\documentclass[12pt,leqno,fleqn,epsfig]{article}
\usepackage{amssymb, epsfig, amsmath, amsthm}
\usepackage{mathrsfs}       

\newcommand{\R}{\mathbb{R}}
\newcommand{\PP}{\mathbb{P}}

\newcommand{\Z}{\mathbb{Z}}

\newcommand{\N}{\mathbb{N}}

\newcommand{\osc}{\operatorname*{osc}}

\newcommand{\trace}{\operatorname{trace}}

\newcommand{\supp}{\operatorname*{supp}}

\newcommand{\const}{\operatorname*{const}}

\newcommand{\bb}{\begin{equation}}
\newcommand{\ee}{\end{equation}}
\newcommand{\bq}{\begin{eqnarray}}
\newcommand{\eq}{\end{eqnarray}}
\newcommand{\bqn}{\begin{eqnarray*}}
\newcommand{\eqn}{\end{eqnarray*}}
\newcommand{\supl}{\sup\limits}
\newcommand{\var}{\varepsilon}

\newcommand{\intl}{\int\limits}

\newcommand{\Beweisende}{\rule{0.2cm}{0.2cm}}

\newcommand{\D}{\displaystyle}

\newcommand{\intmw}{{\int\hspace{-830000sp}-\!\!}}

\newcounter{secnum}

\newtheorem{thm}{Theorem}[section]
\newtheorem{thm1}{Theorem}[]
\newtheorem{cor}[thm]{Corollary}
\newtheorem{lem}[thm]{Lemma}

\theoremstyle{definition}

\newtheorem{defin}[thm]{Definition}
\newtheorem{rem}[thm]{Remark}

\title{The Euler equations  in  a critical case of the generalized Campanato space} 
 
\author{Dongho Chae$^{(*)}$ and J\"{o}rg Wolf$^{(\dagger)}$\\
\ \\
 Department of Mathematics$^{(*), (\dagger)}$ \\
Chung-Ang University\\
Dongjak-gu Heukseok-ro 84\\
Seoul 06974, Republic of Korea\\
and \\
School of Mathematics$^{(*)}$ \\
Korea Institute for Advanced Study\\
Dongdaemun-gu   Hoegi-ro 85 \\
Seoul 02455, Republic of Korea\\
$^{(*)}$e-mail: dchae@cau.ac.kr \\
$^{(\dagger)}$e-mail: jwolf2603@cau.ac.kr}

\date{}
\begin{document}
\maketitle
\begin{abstract}
In this paper we prove  local in time well-posedness for the incompressible Euler equations  in  $\Bbb R^n$ for the initial data  in $\mathscr {L}^{ 1}_{ 1(1)}(\mathbb {R}^{n}) $, 
which corresponds to a critical case of the  generalized Campanato spaces 
$ \mathscr {L}^{ s}_{ q(N)}(\mathbb {R}^{n})$. The space  is studied  extensively in our companion paper\cite{trans}, and in the critical case 
we have embeddings  $ B^{1}_{\infty, 1} (\Bbb R^n)  \hookrightarrow \mathscr {L}^{ 1}_{ 1(1)}(\mathbb {R}^{n}) \hookrightarrow C^{0, 1} (\Bbb R^n)$, 
where  $B^{1}_{\infty, 1} (\Bbb R^n)$ and $ C^{0, 1} (\Bbb R^n)$ are the Besov space and the Lipschitz space respectively. 
In particular  $\mathscr {L}^{ 1}_{ 1(1)}(\mathbb {R}^{n}) $ contains non-$C^1(\Bbb R^n)$ functions as well as linearly growing functions at spatial infinity.
We  can also construct a class of  simple initial velocity  belonging to $ \mathscr {L}^{ 1}_{ 1(1)}(\mathbb {R}^{n})$,  for which  the solution to the Euler equations
blows up in finite time.
  \noindent{\bf
keywords:} Euler equation, generalized Capanato space, local well-posedness

\end{abstract}

\tableofcontents

\section{Introduction}
\label{sec:-1}
\setcounter{secnum}{\value{section} \setcounter{equation}{0}
\renewcommand{\theequation}{\mbox{\arabic{secnum}.\arabic{equation}}}}

Let $ 0< T < +\infty$ and $ Q_T = \R^{n} \times (0,T)$ with $n\in \Bbb N, n\geq 2$. We consider the homogeneous incompressible Euler equations 
\begin{equation}
\begin{cases}
\partial _t v + (v\cdot \nabla) v = -\nabla p \quad  \text{ in}\quad  Q_T,
\\[0.1cm]
\nabla \cdot v =0\quad  \text{ in}\quad  Q_T,\end{cases}
\label{euler}
\end{equation}
equipped with the initial condition 
\begin{equation}
v= v_0\quad  \text{ on}\quad  \R^{n}\times \{ 0\},
\label{initial}
\end{equation}
where $v=(v_1, \cdots, v_n)=v(x,t)$ represents  the velocity of the fluid flows, and $p=p(x,t)$ denotes the scalar pressure. 
The system of Euler equations is of fundamental importance in the mathematical fluid mechanics(see e.g. books\cite{maj, lio, che} or survey paper\cite{con}). Therefore, many authors studied the local well-posedness/ill-posedness of the Cauchy problem \eqref{euler}-\eqref{initial} in various function spaces\cite{kat1, kat2, yud, vis1, vis2, cha, pak, lem, bour1, bour2, bar, lem, 
che, che1, ches, mis}. In particular it is shown that the system \eqref{euler}-\eqref{initial} is locally well-posed in the critical Besov space $ B^{1}_{\infty, 1} (\Bbb R^n)$\cite{pak, lem}, but
 ill-posed  in the Lipschitz space $C^{0, 1} (\Bbb R^n)$\cite{bour1}.
 
  \vspace{0.2cm} 
Our aim in this paper is to show the local well-posedness in a critical generalized Campanato space,  which is embedded into $ C^{ 0,1}(\R^{n})$,  but 
larger than the Besov space $ B^1_{\infty, 1} (\Bbb R^n)$. Furthermore, our  function space include  linearly growing  functions at infinity as well as  bounded functions.  Furthermore it also contains non-$C^1(\Bbb R^n)$ functions as shown in our companion paper\cite{trans}.\\
At a first glance one may think it is impossible to get result of local well-posedness in such function spaces  due to the following example.
\begin{equation}\label{ex1}
v(x, t) = \frac{T_{ \ast}}{T_{ \ast}-t} (x_1,- x_2)^{ \top},\quad (x,t)\in \R^{2}\times (0,T_{ \ast}), 
\end{equation} 
which solves \eqref{euler} with $ v_0 (x)= (x_1,- x_2)^{ \top}$ and $ p(x,t)=
- \frac{1}{2} \frac{T_{ \ast}}{(T_{ \ast}-t)^2} (x_1^2-x_2^2 + T_{ \ast} (x_1^2+ x_2^2)) $. 
Since $ T_{ \ast}>0$ can be chosen arbitrarily small independent of the size of $v_0$, the Euler equations with linear growing initial  data is in general ill-posed. 

We observe that in the above solution one has freedom to choose  the pressure with quadratic growth depending on both the time derivative of $ v$ and the convection term.

  \vspace{0.2cm} 
 In order to avoid such pathological case  we shall restrict our class of solutions by imposing extra condition on choice of  the pressure.
More specifically will introduce a pressure operator $ \Pi = \Pi(v,v)$ such that possible linear growing solutions to \eqref{euler} with $ \nabla p = \nabla \Pi $ are determined uniquely.

\vspace{0.2cm}
{\it In this paper we call a pair $ (v, p)$ a solution to the Euler equations if $ (v,p)\in L^\infty(0, T; L^2_{ loc}( \R^{n})) \times L^\infty(0, T; L^2_{ loc}( \R^{n}))$,   both $ \nabla v$ and $D^2 p$ are bounded in $ Q_T$, and \eqref{euler} holds  a.e. in $ Q_T$. } 

\vspace{0.3cm}  
We start our discussion with the following notion  of equivalent solutions.

\begin{defin}
\label{def1.1}
1. Two solutions  $ (v_1, p_1)$ and $(v_2, p_2)$ are called {\it equivalent} to each other $  (v_1, p_1)\sim(v_2, p_2)$,  if there 
exists $ \xi \in C^{ 1,1}([0,T]; \Bbb R^n)$  such that for almost every $ (x,t)\in Q_T$
\[
v_2(x, t) = v_1(x+ \xi (t), t)- \dot{\xi }(t),\quad   \nabla p_2(x, t) = \nabla p_1(x+ \xi (t), t)- \ddot{\xi}(t) .
\]

2. A solution $ (v,p)$ is called {\it centered},  if 
\begin{equation}
v(0, t)= 0\quad \forall\,t\in [0,T]. 
\label{centered}
\end{equation}

3. We say  a solution $ (v,p)$  {\it decays} as $ | x| \rightarrow +\infty$,  if there exists $ (u, q)\sim (v,p)$ such that 
$  \frac{1}{r^n}\intl_{B(r)}   | u(x, t)| dx \rightarrow 0 $ as $ r\rightarrow +\infty$ for all $ t\in (0,T)$, where $B(r)$ denotes the ball with radius $r$, with its center at the origin.

\end{defin}

\begin{rem}
\label{rem1.1}  1. Clearly, the relation $ \sim$ between two solutions to the Euler equations defines an equivalence relation. 
Given a solution $ (v,p)$ to \eqref{euler} the set $ [(v, p)]$  containing all solutions to \eqref{euler} which are equivalent 
to $ (v, p)$ forms the unique  equivalence class,  which in particular contains $(v,p) $.  
Furthermore, each equivalence class $ [(v, p)]$ contains a centered solution. Indeed,   
we may  find  a solution $ \xi \in C^{ 1,1}([0,T]; \Bbb R^n)$ to the ordinary differential equations
\[
\dot{\xi }(t)= v(\xi(t), t ),\quad  t\in (0,T).  
\]
Setting 
 \[
V(x, t) = v(x+ \xi (t), t)- \dot{\xi }(t),\quad   P(x, t) =  p(x+ \xi (t), t)- \ddot{\xi }(t) \cdot x,\quad  t\in (0, T),  
\]
it is obvious  that $ (V, P)$ is centered and $ (V, P) \sim (v,p)$. 

\vspace{0.5cm}  
2. As an example of non-equivalent  solutions   in $\Bbb R^2$ we consider $ (v, p)$ and $ (u, q)$ both  satisfying  the  same initial condition \eqref{initial}, and  defined by
\begin{align*}
&\hspace*{-1.5cm}v(x,t) = (x_1+x_2, x_1-x_2),\quad -\nabla  p(x,t) = (2x_1, 2x_2),   
\\
&\hspace*{-1.5cm}u(x,t) = (x_1+e ^t x_2, e^t x_1- x_2), \quad  -\nabla q(x.t) = ((e^{ 2t}+1)x_1+e ^tx_2, ( e^{ 2t}-1)x_2 + e^tx_1).
\end{align*}   
This example also shows that we cannot expect uniqueness in the class of solutions with linear growth at infinity without restriction of the pressure as mentioned above.

\vspace{0.3cm}
3. Let $ (v, p)$ be a solution to \eqref{euler} the  fixed properties above. Suppose that $ v(t)\in L^2(\R^{n})$ for all $ t\in [0,T]$, then it holds that $ \| v(t)\|_2 = 
\| v(0)\|_2$ for all $ t\in (0, T)$. Indeed, by interpolating between $v\in L^\infty(0, T; L^2 (\Bbb R^n))$ and $\nabla v\in L^\infty (Q_T) $ one has $v\in 
L^3 (0, T; L^{\frac{3n}{n-1}} (\Bbb R^n))$, in which class we can perform integration by part in the convection term and the pressure term to make them vanish, and finally to get the desired energy conservation.


\end{rem}

\vspace{0.3cm}
Let us introduce the spaces we will use throughout the paper. Let $ N\in \N \cup  \{0\}:=\Bbb N_0$.  By $ \mathcal{P}_N$  ($ \dot{\mathcal{P}} _N$ respectively),  denotes the space of all polynomial (all homogenous polynomials respectively)
of degree less or equal $ N$.  We equip the space $ \mathcal{P}_N$ with the norm $ \| P\|_{ (p)}= \| P\|_{ L^p(B(1))}$. Note that 
since $ \dim (\mathcal{P}_N)<+\infty$ all norms $ \| \cdot \|_{ (p)}, 1 \le p \le  \infty$,  are equivalent.

\vspace{0.5cm}  
 Let $ f\in L^2_{ loc}(\R^{n}), 1 \le p \le +\infty$. For  $ x_0\in  \R^{n}$ and $ 0< r<\infty$ we define 
the oscillation 
\[
\osc_{ p, N} (f; x_0, r):= | B(r)|^{ - \frac{1}{p}} \inf_{ P\in \mathcal{P}_N}  \| f-P\|_{ L^p(B(x_0, r))}. 
\] 

Then, we define for $ 1 \le q, p \le +\infty$ and $ s\in [0, N+1)$  the spaces
$$
 \mathscr{L}^{ s}_{q (p, N)}(\R^{n})= \Big\{ f\in L^p_{ loc}(\R^{n}) \,\Big|\,|f|_{  \mathscr{L}^{ s}_{q (p, N)}} 
  := \supl_{ x_0\in \R^{n}} \left( \sum_{j\in \Bbb Z} \Big(2^{ -sj}\osc_{ p, N} (f; x_0, 2^{j})\Big)^q\right)^{\frac1q} <+\infty\Big\},
$$
Furthermore, by $ \mathscr{L}^{k, s}_{q (p, N)}(\R^{n}) $, $k\in \N$,  we denote the space of all $ f\in W^{k,\, p}_{ loc}(\R^{n})$ such that $ D^k f\in  \mathscr{L}^{ s}_{q (p, N)}(\R^{n})$.  The space $ \mathscr{L}^{k, s}_{q (p, N)}(\R^{n})$ will be equiped with the  norm   
\[
\| f\|_{ \mathscr{L}^{k, s}_{q (p, N)} }= | D^k f|_{ \mathscr{L}^{  s}_{q (p, N)} } + \| f\|_{ L^p(B(1))},\quad  f\in 
{ \mathscr{L}^{k, s}_{q (p, N)}}(\R^{n}). 
\]

Note that the oscillation introduced above is attained  by a unique 
polynomial $ P_{ \ast}\in \mathcal{P}_N$. 
Below we recall basic properties   on this space. According  to the characterization theorem of the Triebel-Lizorkin spaces  in terms of  oscillation(cf. \cite[Theorem, Chap. 1.7.3]{tri}), we have
\begin{align*}
\hspace{-.5in}\begin{cases}
{\D f\in F^s_{ r, q}(\Bbb R^n) \quad  \Leftrightarrow   \quad  \|f\|_{ L^{ \min\{r,q\}}}+
\quad  \bigg\|\left( \sum_{j=-\infty} ^0\Big(2^{ -sj}\osc_{ p, N} (f; \cdot , 2^{j})\Big)^q\right)^{\frac1q}\bigg\|_{ L^r} <+\infty. }
\\[0.5cm]
0< r<+\infty, 0< q \le \infty, \quad s> \Big(\frac{1}{r} - \frac{1}{p}\Big)_{ +},\quad  s> \Big(\frac{1}{q} - \frac{1}{p}\Big)_{ +},
\end{cases}
\end{align*}
and we could regard the spaces $  \mathscr{L}^{s}_{q (p, N)}(\R^{n} ) $ as an extension of the  limit case of 
$F_{ r,q}^s(\R^{n} ) $
as $ r\rightarrow +\infty$. 
In case $ q=+\infty$ and $ s >0$ we get the usual Campanato spaces with the isomorphism relation(cf. \cite{ca, gia})
\[
\mathscr{L}^{n+ps, p}_N(\R^{n}) \cong  \mathscr{L}^{s}_{\infty (p, N)}(\R^{n}).
\] 
 Furthermore, in the case $ N=0, s=0$ and $ q= \infty$ we get the space of bounded 
mean oscillation, i.e.,
\[
\mathscr{L}^{0}_{\infty (p, 0)}(\R^{n})\cong BMO. 
\]
In case $ N=-1$ and $ s\in (-\frac{n}{p}, 0) $ the above space coincides with the usual Morrey space $ {\cal M}^{ n+ ps}(\R^{n} )$. 
Our aim in this paper is to prove the local well-posedness of the Euler equations in  the critical space 
$ \mathscr{L}^{1}_{1 (p, 1)}( \R^{n})$. We recall the following embedding relations(see \cite{trans}).
\begin{equation}
B^{ 1+ \frac{n}{r}}_{ r, 1} \hookrightarrow  B^{1}_{\infty, 1} (\Bbb R^n)   \hookrightarrow  \mathscr{L}^{1}_{1 (p, 1)}( \R^{n})
 \hookrightarrow C^{0, 1}(\R^{n}). 
\label{1.8}
\end{equation}
Accordingly, 
\begin{equation}
\| \nabla u\|_{\infty}  \le c\| u\|_{\mathscr{L }^{1}_{1 (p, 1)}}.
\label{1.11}
\end{equation}

Furthermore,  for every $ f\in  \mathscr{L}^{k}_{1 (p, k)}( \R^{n}), k\in \N_0 = \N \cup \{0\}$, there exists a unique $ \dot{P}^k_{ \infty}(f)\in  \dot{\mathcal{P}} _k  $,  such that for all $ x_0 \in \R^n$
\[
f \quad  \text{converges asymptotically to} \quad  \dot{P}^k_{ \infty}(f)\quad  \text{as }\quad  |x| \rightarrow +\infty.
\]
 The exact meaning of this asymptotic limit is stated in Theorem\,\ref{thm5.4} (see also  \cite[Section 2]{trans}).

\vspace{0.5cm}  
We also  introduce  the following critical homogenous space 
\begin{equation}\label{center}
\mathscr{\dot{L} }^{1}_{1 (p, 1)}(\R^{n})= \Big\{u\in \mathscr{L}^{1}_{1 (p, 1)}( \R^{n}) \,\Big|\, u(0)=0\Big\}.
\end{equation}

The space $ \mathscr{\dot{L} }^{1}_{1 (p, 1)}( \R^{n})$ will be equipped with the homogenous norm 
 \[
\| u\|_{ \mathscr{\dot{L} }^{1}_{1 (p, 1)}} = | u|_{\mathscr{L}^{1}_{1 (p, 1)} } + |P_\infty^0( \nabla u)|,\quad  u\in \mathscr{\dot{L} }^{1}_{1 (p, 1)}( \R^{n}). 
\]
We recall that $ u\in \mathscr{\dot{L} }^{1}_{1 (p, 1)}( \R^{n})$ implies $\nabla u\in 
 \mathscr{\dot{L} }^{0}_{1 (p,0)}( \R^{n})$.

 By $ \mathscr{\dot{L} }^{1}_{1 (1), \sigma }( \R^{n})$ we denote the subspace of all $ u\in \mathscr{\dot{L} }^{1}_{1 (p, 1)}( \R^{n})$ such that $ \nabla \cdot u=0$ 
 almost everywhere in $ \R^{n}$.  Next, we focus on the pressure $ p$, which satisfies the Poisson equation
\begin{equation}
-\Delta  p = \nabla v: (\nabla v)^T\quad  \text{ in}\quad  \R^{n}.
\label{1.12}
\end{equation} 
In contrast to the decaying case this problem for $ \nabla p$   is not well posed in the space $ \mathscr{L}^{1}_{1 (p, 1)}(\R^{n})$, since if
 $\nabla p\in  \mathscr{L}^{1}_{1 (p, 1)}(\R^{n})$  solves  \eqref{1.12}, then the function $   p+ Q$   for any  $ Q\in \mathcal{P}_1$  also solves it. 
The same problem occurs for the general Poisson equation.  In order to have uniqueness of solution we add an   asymptotic condition 
for $ \nabla p$   as  
$ | x| \rightarrow +\infty$ together  with a condition at one point. 
We have the following

\begin{thm}
\label{thm1.2} For every matrix $ H=\{ H_{ \alpha \beta }\}\in \mathscr{L}^{1}_{1 (p, 1)}(\R^{n}),  a_0\in \R$ and $ Q_\infty\in  \dot{\mathcal{P}} _1 $ 
there exists a unique solution $f\in  \mathscr{L}^{1}_{1 (p, 1)}(\R^{n})$ to  the problem  
\begin{equation}
\begin{cases}
- \Delta f = \nabla \cdot  \nabla \cdot  H\quad  \text{ in }\quad  \R^{n},
\\[0.3cm]
f(0)=a_0,\quad  \dot{P}^1 _\infty(f) =Q_\infty.
\end{cases}
\label{1.13}
\end{equation}
Furthermore, there exists a constant $c=c(n, p)  $ such that 
\begin{equation}
| f|_{ \mathscr{L}^{1}_{1 (p, 1)}} \le c | H|_{ \mathscr{L}^{1}_{1 (p, 1)}},\quad  
\| f\|_{ \mathscr{L}^{1}_{1 (p, 1)}} \le c (| H|_{ \mathscr{L}^{1}_{1 (p, 1)}}+ | a_0|+ \| Q_\infty\|). 
\label{1.14}
\end{equation}
\end{thm}

The proof of Theorem\,\ref{thm1.2} is based on  Theorem\,\ref{thm6.3} given in Section\,3.

\vspace{0.3cm}
Next, we discuss the problem of  defining  the pressure. We first  define an operator  $ \nabla \Pi: 
\mathscr{L}^{1}_{1 (p, 1), \sigma }(\R^{n})\times \mathscr{L}^{1}_{1 (p, 1)}(\R^{n}) \rightarrow \mathscr{\dot{L} }^{1}_{1 (p, 1)}(\R^{n})$ as follows. This definition is based on the following theorem, 
which is an immediate  consequence   of Theorem\,\ref{thm6.4} (see also Remark\,\ref{rem6.5})

\begin{thm}
 \label{thmpressure} Let $ (u, v) \in \mathscr{L}^{1}_{1 (p, 1), \sigma }(\R^{n})\times \mathscr{L}^{1}_{1 (p, 1)}(\R^{n})$. 
 There exists a function $ \pi \in \mathscr{L}^{1,1}_{1 (p, 1)}(\R^{n})$ solving the Poisson equation
 \begin{equation}
\begin{cases}
- \Delta \pi  = \nabla \cdot  \nabla \cdot  (u \otimes v)\quad  \text{ in }\quad  \R^{n},
\\[0.3cm]
 \nabla \pi (0)=0  ,\quad  \dot{P}^1 _\infty(\nabla \pi ) =- \frac{1}{n}P^0_\infty(\nabla u: (\nabla v)^{ \top}) x,
\end{cases}
\label{1.13-p}
\end{equation}
which is unique up to a constant. 
\end{thm}

Now we are ready to introduce the following definition.
\begin{defin}
 \label{defpressure} 
Let $ (u, v) \in \mathscr{L}^{1}_{1 (p, 1), \sigma }(\R^{n})\times \mathscr{L}^{1}_{1 (p, 1)}(\R^{n})$. 
Then by  $ \nabla \Pi (u,v)$ we denote the unique function $ \nabla \pi \in \mathscr{L}^{1}_{1 (p, 1)}(\R^{n})$, 
where $ \pi \in \mathscr{L}^{1,1}_{1 (p, 1)}(\R^{n})$ stands for the solution of  \eqref{1.13-p} 
according to Theorem\,\ref{thmpressure}.
\end{defin}

In particular, in view of \eqref{6.50} and  \eqref{6.51} (cf. Section\,3) it holds
\begin{equation}
\| \nabla \Pi(u, v) \|_{ \mathscr{\dot{L} }^{1}_{1 (p, 1)}} \le c 
 \Big(\| u\|_{ \mathscr{\dot{L} }^{1}_{1 (p, 1)}} \| \nabla v\|_{\infty}+ \| v\|_{ \mathscr{\dot{L} }^{1}_{1 (p, 1)}} \| \nabla u\|_{\infty}\Big). 
\label{1.25}
\end{equation}   

We are now in a position to present our first main result.

\begin{thm1}[{\bf Local well posedness in $ \mathscr{\dot{L} }^{1}_{1 (p, 1)}(\R^{n})$}]  For every $ v_0\in \mathscr{\dot{L} }^{1}_{1 (p,1), \sigma}(\R^{n})$  there exists 
\label{thm1}
\begin{equation}
T_0 \ge \frac{1}{c \| v_0\|_{\mathscr{\dot{L} }^{1}_{1 (p, 1)} }},
\label{1.26}
\end{equation}
and a unique solution $ v\in L^\infty(0, T_0; \mathscr{\dot{L}}^{1}_{1 (1), \sigma }(\R^{n}))$ to \eqref{euler}, \eqref{initial}  in $ Q_{ T_0}$  with pressure $ \pi \in L^\infty(0, T_0; L^2_{ loc}(\R^{n}))$ 
such that $ \nabla \pi \in L^\infty(0, T_0; \mathscr{\dot{L} }^{1}_{1 (p, 1)}(\R^{n}))$ and satisfying 
\begin{equation}
\nabla \pi = \nabla \Pi (v,v).
\label{1.27}
\end{equation} 

\end{thm1}

\begin{rem}
\label{rem1.5}
In case of sublinear  growing solutions  the condition   \eqref{1.27}  is automatically satisfied  for the function  $ \nabla \pi(x,t)= \nabla p(x,t)- \nabla p(0, t ) $, if $ (v,p)$ solves  the Euler equations, using the arguments 
in Section\,4.
\end{rem}

Using the Galilean transform $(x, t) = (y+ t a, t)$, $ a\in \R^{n}$, $ (x, t )\in Q_{ T_0}$,  with $ a= -v_0(0)$,  we obtain the following local well posedness in $ \mathscr{L}^{1}_{1 (p, 1)}(\R^{n}).$

\begin{thm1}[{\bf Local well posedness in $ \mathscr{L}^{1}_{1 (p, 1)}(\R^{n})$}] 
\label{thm2}
For every $ v_0\in \mathscr{L}^{1}_{1 (p, 1)}(\R^{n})$ with $ \nabla \cdot v_0=0$ there exists 
\begin{equation}
T_0 \ge \frac{1}{c \| v_0- v_0(0)\|_{\mathscr{\dot{L} }^{1}_{1 (p, 1)} }},
\label{1.28}
\end{equation}
and a unique solution $ v\in L^\infty(0, T_0; \mathscr{L}^{1}_{1 (p, 1)}(\R^{n}))$ to 
\eqref{euler}-\eqref{initial}  
 in $ Q_{ T_0}$  with pressure $ p \in L^\infty(0, T_0; L^2_{ loc}(\R^{n}))$
such that  $\nabla p \in L^\infty(0, T_0;  \mathscr{L}^{1}_{1 (p, 1)}(\R^{n}))$ and for almost all $ t\in (0,T)$
\begin{equation}
\nabla p(\cdot,t )-\nabla p(0, t )  = \nabla \Pi (v(t),v(t)),\quad  \nabla p(  v_0(0) t, t)=0.
\label{1.29}
\end{equation}
\end{thm1}

\begin{rem}
\label{rem1.7}
In fact, our  main result stated in Theorem\,\ref{thm2} improves substantially previous result in \cite{pak} in  both directions,  in the sense of regularity and asymptotic behavior at infinity in space. 
Firstly,  we recall  that  by \eqref{1.8}
\begin{equation}
  B^1_{ \infty, 1}(\R^{n} ) \hookrightarrow  \mathscr{L}^{1}_{1 (p, 1)}(\R^{n})\cap L^\infty(\R^{n} )
  \hookrightarrow  \mathscr{L}^{1}_{1 (p, 1)}(\R^{n}).
\label{emb}
 \end{equation}  
Secondly, according to \cite[p. 85]{tay}, (see also \cite{che}) we have the embedding 
\[
 B^1_{ \infty, 1}(\R^{n} )  \hookrightarrow C^1(\R^{n} ) \cap L^\infty(\R^{n} ). 
\]

On the other hand, there exists a function $ f\in \mathscr{L}^{1}_{1 (p, 1)} (\R^{n} )$, which is not in $ C^1(\R^{n} )$ 
(see in \cite[appendix]{trans}).
Consequently,  $ \mathscr{L}^{1}_{1 (p, 1)} (\R^{n} )$ contains less regular functions then  
$  B^1_{ \infty,1}(\R^{n} ) $.  {\it In particular, this implies that we have local well posedness of the Euler equations for 
initial data not in $ C^1(\R^{n} )$. }

Thirdly, since $ \mathscr{L}^{1}_{1 (p, 1)} (\R^{n} )$ contains linear growing functions, and therefore  polynomials  of degree 
less or equal one,  $ \mathscr{L}^{1}_{1 (p, 1)} (\R^{n} )$ is strictly larger then  $ B^1_{ \infty,1}(\R^{n} )$ in the sense of
asymptotic behavior at spatial infinity.  

\end{rem}

Next,  we introduce the notion of equivalent solutions by using the change of coordinates $ (x,t)= (x+ \xi (t), t)$ for a given function $ \xi \in C^{ 1,1}([0,T] ; \Bbb R^n)$. 

 \begin{defin}
\label{def1.6}  A solution $ (v, p ) \in  L^\infty(0, T_0; \mathscr{L}^{1}_{1 (p, 1)}(\R^{n})\times L^2_{ loc}(\R^{n}))$  
to the Euler equations \eqref{euler} is called {\it eligible} if there exists a centered solution $ (V, P) \sim (v,p)$ with 
\begin{equation}
\nabla P(t) = \nabla \Pi (V(t),V(t))\quad  \forall\,t\in [0,T].
\label{1.32}
\end{equation} 
  
\end{defin}

\begin{rem}
As an example of  non-eligible   solutions  we have  solutions in \eqref{ex1}.  In general, we may look for solutions 
$ v(t) = A(t)x$, where $ A\in \R^{n\times n}$ stands for a matrix $ \trace(A)=0$. In \eqref{euler} replacing $ v$ by $ Ax$ we obtain 
the equation 
\begin{equation}
\partial _t Ax + A^2x = - \nabla \pi \quad  \text{ in}\quad  Q_{ T}.    
\label{1.33a}
\end{equation}
 The compatibility condition \eqref{1.27} yields $\nabla \pi =- \frac{1}{n} \trace (A^2)x$. Inserting this identity into    
 \eqref{1.33a} and applying $ \nabla $ to the resultant equations, we are led to the system of ODE's
\begin{equation}
\dot{A}  + A^2 = \frac{1}{n} \trace(A^2)I  \quad  \text{ in}\quad  [0,T].     
\label{1.33b}
\end{equation} 
This system, under the initial condition $ A(0) = A_0$ has a unique local  solution $ A\in C^1([0,T)); \R^{n\times n})$. 
Applying $ \trace $ to \eqref{1.33b},  we deduce that $ \frac{d}{dt} \trace (A)=0$, which together with  
 $ \trace (A_0)=0$, yields $\trace (A) =0$. This shows that $ v=A x$ is a centered  solution of the Euler equations, and by 
 Theorem\,\ref{thm1} this solution is unique.     
The following examples show the  
global existence and finite time blow up depending on the initial data  in $ \mathscr{\dot{L} }^{1}_{1 (p, 1)}(\R^{n})$.

\vspace{0.3cm}
{\it 1. Global existence in $ n=2$.}  Let $ v_0 (x) = A_0 x$, where $ A_0= {\rm diag} (\lambda_{ 0,1}, \lambda _{ 0,2})$. The condition 
$\nabla \cdot v_0 =0$ implies $ \lambda _{ 0,2} = -\lambda _{ 0,1}$.  Then the centered solution $ v$ to \eqref{euler}  must be given 
by $ v=Ax$, where $ A$ solves \eqref{1.33b}. However, noting that $ \lambda _{ 0,1}^2 = \frac{1}{2} \trace A_0^2$, 
it readily seen that $ A \equiv  A_0$  is a global unique  solution  to \eqref{1.33b}     
and $ v = A_0x$ is the global centered solution to \eqref{euler}.   

\vspace{0.3cm}
{\it 2. Global existence and finite time blow up in $ n=3$.} We begin our discussion with the global existence. 
Let $ v_0= A_0 x$ be given with  $A_0= {\rm diag} (\lambda_{ 0,1}, \lambda _{ 0,3}, \lambda _{ 0,2})$ with 
$ \sum_{i=1}^{3} \lambda _{ 0,i}=0$ such that all eigen values $ \lambda _{ 0,i}$ are differnt.  Then the solution $ A$ to \eqref{1.33b} has the form 
$ A={\rm diag} (\lambda_{ 1}, \lambda _{ 2}, \lambda _{ 3}) $  such that $ \sum_{i=1}^{3} \lambda _{i}=0$. 
By $ T_{ \ast}>0$ we denote the maximal time of existence of this solution, i.e. $ \lambda _i$ solve the system of ODE
\begin{equation}
\dot{\lambda} _i + \lambda _i^2 = \frac{1}{3} | \lambda |^2\quad  \text{ in}\quad  (0, T_{ \ast}),\quad  i=1,2,3. 
\label{1.33c}
\end{equation}

We claim $ T_{ \ast}= +\infty$. 
To see this, first we verify that that all eigen values $ \lambda _{ i}(t)$ are different for all $ t\in (0,T_{ \ast})$. In fact, in view of \eqref{1.33c},  the function  
$\mu = \lambda _i - \lambda _j$ for $ i \neq j$  solves the ODE $ \dot{\mu }  + (\lambda _i+\lambda _j)\mu=0$ in 
$ (0,T_\ast)$.  In case $ \mu (t)=0$ for some 
$ t\in (0, T_{ \ast})$ it follows that  $ \mu \equiv 0$, which contradicts to $ \mu (0) \neq 0$.  We now define the differences 
and sum
\begin{align*}
\mu _1 &= \lambda_2- \lambda _3,\quad \mu _2 = \lambda_3- \lambda _1, \quad  \mu _3 = \lambda_1- \lambda _2,
\\
\nu _1 &= \lambda_2+ \lambda _3,\quad \nu _2 = \lambda_3+ \lambda _1, \quad  \nu _3 = \lambda_1+ \lambda _2. 
\end{align*}
Then \eqref{1.33c} yields 
 \begin{equation}
\dot{\mu} _i + \nu _i\mu _i =0\quad  \text{ in}\quad  (0, T_{ \ast}),\quad  i=1,2,3. 
\label{1.33d}
\end{equation} 
Solving this equations, we get 
\[
\mu _i(t) = \mu_i (0) e^{ - \nu _i(t)},\quad  t\in (0, T_{ \ast}). 
\] 
Verifying that $ \sum_{i=1}^{3} \nu _i \equiv  0$,  we obtain $ \prod_{ i=1}^3 \mu_i \equiv c_0:=  \prod_{ i=1}^3 \mu_{ 0,i}$. 
Accordingly, 
\begin{equation}
\mu _3 = \frac{c_0}{\mu _1 \mu _2}. 
\label{1.33e}
\end{equation} 
We now define $ \alpha = \mu _1+\mu _2= - \mu _3$, and $ \beta = \mu _1- \mu _2$. We calculate, 
\[
\nu _3= \lambda _1+ \lambda _2 = \mu _1- \mu _2 + 2 \lambda _3 =\mu _1- \mu _2 - 2 \nu _3.  
\]
Thus,  
\[
 \nu _3= \frac{1}{3}  (\mu _1- \mu _2) = \frac{\beta }{3}. 
\]
Inserting this identity into \eqref{1.33d} for $ i=3$, we get 
\[
\dot{\alpha} + \frac{\beta }{3} \alpha =0\quad  \Rightarrow \quad  \alpha = \alpha (0) e^{ - \frac{1}{3}\beta }.
\] 
On the other hand, in view of \eqref{1.33e} we infer 
\[
\beta ^2 = \mu_1^2 - 2\mu _1 \mu _2 + \mu _2^2 = \alpha ^2 -  4 \mu _1\mu _2 = \alpha ^2 -  4c_0 \mu _3^{ -1}
=\alpha ^2 + 4c_0 \alpha ^{ -1}.  
\]
Inserting $  \alpha = \alpha (0) e^{ - \frac{1}{3}\beta }$, this yields
\begin{equation}
\beta ^2 = \alpha (0) e^{ - \frac{2}{3}\beta } + 4c_0\alpha (0)^{ -1} e^{ \frac{1}{3}\beta }\quad  \text{ in}\quad  
(0, T_{ \ast}).  
\label{1.33f}
\end{equation}
This shows that $ \beta $ is bounded, which also implies that $ \alpha $ is bounded. Hence, $ \lambda _i, i=1,2,3$ are bounded.
Whence, $ T_{ \ast}=+\infty$. 
  
\vspace{0.5cm}  
Next,  we give an example of finite time blow up. As we have seen above this is only possible if two eigenvalues are equal. 
Thus, we may assume that $ \lambda = \lambda _1=\lambda _2 >0$ and $ \lambda _2 = -2\lambda $. Then, in view of \eqref{1.33c} $ \lambda $
solves the Riccati equation 
\begin{equation}
\dot {\lambda} = \frac{1}{3} \lambda ^2\quad  \text{ in}\quad  (0,T_{ \ast}),
\label{1.33g}
\end{equation}
which has the unique solution 
\[
\lambda(t) = \frac{3\lambda_0}{\lambda (0)t - 3},\quad t\in (0,T_{ \ast}),\quad  T_{ \ast}= \frac{3}{\lambda (0)}. 
\]  
\end{rem}

For the case of  initial data with sub linear growth we get the following third main result which can be directly compared with the known results in Besov spaces 

\begin{thm1}[Local well posedness in $ \mathscr{L}^{1}_{1 (p, 1)}\cap BMO$]
\label{thm3}
 For every $ v_0\in \mathscr{L}^{1}_{1 (p, 1), \sigma }\cap BMO$,  there exists 
\begin{equation}
T_0 \ge \frac{1}{c| v|_{\mathscr{L}^{1}_{1 (p, 1)}}}, 
\label{1.34}
\end{equation}
and a unique   solution $ (v, p ) \in 
L^\infty(0, T_0; (\mathscr{L}^{1}_{1 (p, 1)}\cap BMO)\times 
BMO)$ to \eqref{euler}-\eqref{initial}  in $ Q_{ T_0}$ such that $ (p)_{ 0, 1}=0$.  Such solution is also eligible.   

\vspace{0.2cm} 
In case $ v_0 \in L^\infty(\R^{n})$  the above solution is bounded. 
 
\end{thm1}

We also are able to generalize the Baele-Kato-Majda condition\cite{bea} to the non decaying case as follows.   

\begin{thm1}
\label{thm4} 
Let $v_0\in \mathscr{L}^{1+\delta }_{q (p, 1), \sigma }(\R^{n}) \cap \mathscr{L}^1_{ 1 (2, 1)}(\R^{n}), 1< p< +\infty, 1 \le  q \le +\infty, \delta \in (0,1),$ fulfilling
 \begin{equation}
\forall \varepsilon >0\quad \exists k\in \N\quad  \text{such that}\quad 
\sup_{x_0\in \R^{n} }  \sum_{j=k}^{\infty} 2^{ -j}\osc_{2,1}(v_0; x_0, 2^j ) \le \varepsilon. 
\label{1.34a}
\end{equation}
 Let $  v\in   L^\infty_{ loc}([0, T_{ \ast}); \mathscr{L}^1_{ 1 (2, 1)}(\R^{n}))$ be an eligible  solution to \eqref{euler} 
according to Theorem\,\ref{thm2}. 
Furthermore, assume that 
\begin{equation}
 \intl_{0}^{T_{ \ast}} | \omega (\tau )|_{ BMO} + | P^0_{ \infty}(\nabla v(\tau ))| d\tau <+\infty. 
\label{1.37}
\end{equation}
Then, $ v \in  L^\infty(0, T_{ \ast}; (\mathscr{L}^{ 1+\delta }_{ q (p, 1)}\cap \mathscr{L}^{ 1 }_{ 1 (2, 1)})(\R^{n} )  )$, and the 
solution can be extended to $ [0, T_{ \ast}+ \eta ]$ for some $ \eta >0$. 
\end{thm1}

\begin{rem}

1. We wish to emphasize  that in  the case of sublinear growing  initial data  the condition  \eqref{1.34a} is obviously  satisfied. 
Furthermore, in that case  it holds $  P^0_{ \infty}(\nabla v_0)=0$, and as shown in Section\,7 this implies $  P^0_{ \infty}(\nabla v(\tau ))=0$ 
for all $ \tau \in [0, T_{ \ast})$. Hence,   
\eqref{1.37} reduces to  Kozono-Taniuchi's condition in \cite{koz}
\begin{equation}
 \intl_{0}^{T_{ \ast}} | \omega (\tau )|_{ BMO}  d\tau <+\infty, 
\label{1.38}
\end{equation}
which is a refined version of the Beale-Kato-Majda criterion\cite{bea}.
  \end{rem}

\begin{rem}
The examples of solutions  $ (v,p) \in C([0, T_{ \ast}); \mathscr{L}^1_{ 1 (p, 1)}(\R^{n}))$ in Remark 1.9 show that even if $\omega(t)=0$ for all $t\in [0, T_*)$,   the solution can 
blow up at $t=T_*$,    and it holds $ \int_0 ^{T_*}  | P^0_{ \infty}(\nabla v(t))| dt=+\infty$,    which implies the necessity  of the second integrand of  \eqref{1.37} in the case of solutions having linear growth at inifinity.
\end{rem}

\section{Notations and preliminariy lemmas}
\label{sec:-10}
\setcounter{secnum}{\value{section} \setcounter{equation}{0}
\renewcommand{\theequation}{\mbox{\arabic{secnum}.\arabic{equation}}}}

Let $ X=\{ X_j\}_{ j\in \Z}$ be a sequence of non-negative real numbers. 
 We define $S_{\alpha, q} : X=\{ X_j\}_{j\in \Bbb Z}  \mapsto Y=\{ Y_j \}_{j\in \Bbb Z}$, where
\[
Y_j=S_{ \alpha, q}(X)_j = 2^{ j\alpha }\Big(\sum_{i=j}^{\infty} (2^{ -i\alpha } X_i)^q\Big)^{ \frac{1}{q}},\quad  j\in \Z. 
\]
From the above definition it follows that 
\begin{equation}
\|S_{0, q}(X)\|_{\ell^\infty}\leq \| X\|_{ \ell^q},\quad  \forall\,X \in \ell^q. 
\label{10.1}
\end{equation}

Given $ X=\{ X_j\}_{j\in \Bbb Z} ,  Y=\{ Y_j \}_{j\in \Bbb Z}$,  we denote $X\leq Y$ if  $X_j\leq Y_j$ for all $j\in \Bbb Z$.
Throughout this paper, we frequently make use of the following lemma, which could be regarded as a generalization of the result in \cite{bou}.
 
\begin{lem}
\label{lem10.1}  For all $ \beta < \alpha $  and $ 0< p \le  q \le +\infty$ it holds 
\begin{equation}
 S_{ \beta , q} (S_{ \alpha, p}(X)) \le  \frac{1}{1- 2^{ -( \alpha -\beta) }}S_{\beta , q} (X). 
\label{10.3}
\end{equation}

\end{lem}
For the proof see in \cite[Section\,2]{trans}.

\vspace{0.3cm}
In what follows we provide important properties of the space $  \mathscr{L}^{ k,s}_{  q(p, N)}(\R^{n})$ 
such as embedding properties, equivalent norms. 
First, let us recall the definition of the generalized mean 
for distributions $ f\in {\mathscr S}'$, where $ {\mathscr S}$ denotes the usual Schwarz 
{class of rapidly decaying functions}.  For $ f\in {\mathscr S}'$ 
and $ \varphi \in {\mathscr S}$ we dfine the convolution
\[
f \ast \varphi (x)=  \langle f, \varphi (x- \cdot )\rangle,\quad  x \in \R^{n},
\]
where $<\cdot, \cdot>$ denotes the dual pairing.
Below we use the notation $\N_0= \N\cup\{0\}$.
Then, $ f \ast \varphi \in  C^\infty(\R^{n} )$ and for every  multi index $ \alpha \in \N_0^n$ it holds 
\[
D^\alpha (f \ast \varphi) = f \ast (D^\alpha\varphi)= (D^\alpha f) \ast \varphi. 
\]

Given  $ x_0\in  \R^{n}, 0< r< +\infty$ and  $ f\in {\mathscr S}'$ we define the mean
\[
[f]^\alpha _{ x, r} = f \ast D^\alpha \varphi _r(x).
\]  
where $ \varphi _{r}(y) = r^{ -n}\varphi(r^{ -1}(y))$, and $ \varphi \in C^{\infty}_{ c}(B(1))$ stands for the standard mollifier, such that $ \intl_{ \R^{n} } \varphi dx =1$. Note that in case $ f\in L^1_{ loc}(\R^{n} )$ we get 
\[
[f]^0_{ x, r} =\intl_{ \R^{n}} f(x-y) \varphi _r(y)  dy  = 
\intl_{ B(x, r)} f(y) \varphi _{ x, r}(-y)  dy,
\]
where $ \varphi_{ x, r}= \varphi _r(\cdot +x)$. Furthermore,  from the above definition it follows that 
 \begin{equation}
[f]^{ \alpha }_{ x, r} = (D^\alpha f) \ast \varphi _{r}(x) = [D^\alpha f]^0_{ x, r}.  
 \label{5.1a}
 \end{equation}
For $ f\in L^1_{ loc}(\R^{n} )$ and  $ \alpha \in \N_0^n$ we immediately get 
\begin{equation}
 [f]^{ \alpha }_{ x, r} \le c r^{ -|\alpha |-n}  \|f\|_{ L^1(B(x, r))}\quad \forall x\in \R^{n}, r>0 . 
\label{5.1b}
 \end{equation} 

 By the following lemma we introduce the mean polynomial $ P^N_{ x_0, r}(f)$ together with its properties. The proof of this and all other  lemmas of this section can be found in \cite[Section\,3]{trans}. 
\begin{lem}
\label{lem5.4}
Let $x_0\in \R^{n}, 0<r<+\infty$ and $ N\in \N_0$. For every $ f\in \mathscr{S}'$ there exists a unique polynomial 
$ P^N_{ x_0, r}(f) \in \mathcal{P}_N$ such that 
\begin{equation}
[f- P^N_{ x_0, r}(f)]_{ x_0, r}^\alpha =0\quad  \forall\,| \alpha | \le N. 
\label{5.1}
\end{equation}
In addition, the mapping  $ P^N_{x_0, r }: L^p(B(x_0, r)) \rightarrow \mathcal{P}_N, 1 \le p \le +\infty$,  defines a projection, i.e.
\begin{align}
&\qquad  P^N_{ x_0, r}(Q) = Q\quad  \forall\,Q \in \mathcal{P}_N,
\label{5.3}
\\
&\| P^N_{ x_0, r} (f)\|_{ L^p(B(x_0, 4r))} \le c\| P^N_{ x_0, r} (f)\|_{ L^p(B(x_0, r))} \le  c\| P^N_{ 0, 1}\|_{ p}\| f\|_{ L^p(B(x_0, r))}.
\label{5.4}
\end{align}
where 
\begin{align}
 \| P^N_{ 0, 1}\|_p = \sup_{\substack{g\in L^p(B(1)) \\  g \neq 0}}\frac{\|P_{ 0,1}^N(g)\|_{ L^p(B(1))}}{\|g\|_{ L^p(B(1))}} 
 =\sup_{\substack{g\in L^p(B(x_0, r)) \\ g \neq 0}}
 \frac{\|P_{ x_0,r}^N(g)\|_{ L^p(B(x_0, r))}}{\|g\|_{ L^p(B(x_0, r))}}. 
 \label{5.4a}
\end{align}

Furthermore, 
for all  
 $ f\in W^{p,\, j}(B(x_0, r))$, $ 1 \le p < +\infty,  1 \le j \le N+1$,  it holds 
 \begin{equation}
 \| f- P^N_{ x_0, r}(f)\|_{ L^p(B(x_0, r))} \le c r^j\sum_{| \alpha |=j}\| D^\alpha  f- D^\alpha   P^N_{ x_0, r}(f)  \|_{ L^p(B(x_0, r))}. 
 \label{5.5}
 \end{equation}

\end{lem}

 \begin{rem}
 \label{rem1.6} From   \eqref{5.5}  with $ j = N+1$ we get the generalized Poincar\'e inequality 
 \begin{equation}
 \begin{cases}
 \| f- P^N_{ x_0, r}(f)\|_{ L^p(B(x_0, r))} \le  c r^{ N+1}
 \| D^{ N+1}  f\|_{L^ p( B(x_0, r))}
  \\[0.3cm]
  \forall\,f \in W^{N+1,\, p}(B(x_0, r)).  
 \end{cases}
 \label{5.8}
 \end{equation}
 
 \end{rem}
 
 \begin{cor}
 \label{cor5.7}
 For all $ x_0 \in \R^{n}, 0< r < +\infty$, $ N\in \N_0$, and $ 1 \le p < +\infty$ it holds 
  \begin{equation}
 \| f - P^N_{ x_0, r}(f)\|_{ L^p(B(x_0, r))} \le c \inf_{ Q \in \mathcal{P}_N} \| f - Q\|_{ L^p(B(x_0, r))}
  =cr^{ \frac{n}{p}} \osc_{ p,N} (f; x_0, r). 
 \label{5.9}
 \end{equation}

 \end{cor}

\vspace{0.3cm}
In our discussion below and in the sequel of the paper it will be convenient to work with smooth functions.  
Using the standard mollifier we get the following estimate in  $ \mathscr{L}^{ k,s}_{q (p, N) }(\R^{n})$ for the mollification.

\begin{lem}
\label{lem5.3}
Let $ \varepsilon >0$. 
Given  $ f \in \mathscr{S}'$, we define the mollification  
\[
f_\var (x) =  [f]_{ x, \varepsilon }^{ 0}= f \ast \varphi _{ \varepsilon }(x), \quad  x\in \R^{n}.   
\]
1. For all $ f\in \mathscr{L}^{ k,s}_{q (p, N) }(\R^{n})$,  and all  $ \var >0$ it holds 
\begin{equation}
| f_\var |_{ \mathscr{L}^{ k,s}_{q (p, N) } } \le c| f|_{\mathscr{L}^{ k,s}_{q (p, N) } }. 
\label{5.10}
\end{equation} 
 
 2. Let $ f\in L^p_{ loc}(\R^{n})$ such that for all $ 0<\var<1 $,
\begin{equation}
| f_\var |_{ \mathscr{L}^{ k,s}_{q (p, N) } } \le c_0,
\label{5.11}
\end{equation}  
then $ f\in \mathscr{L}^{ k,s}_{q (p, N) }(\R^{n})$ and it holds $ | f |_{\mathscr{L}^{ k,s}_{q (p, N) }} \le c_0.$
\end{lem}

Next, we  provide the following embedding properties.  
First, let us introduce  the definition of the  projection 
to the space of  homogenous polynomial $\dot{P}^N_{ x_0,r}:  \mathscr{S}' \rightarrow   \dot{\mathcal{P}} _N$ 
defined by 
\[
\dot{P}_{ x_0, r}^N(f) (x) = \sum_{| \alpha |=N}  \frac{1}{ \alpha !} [ f]^{ \alpha }_{ x_0, r} x^\alpha,\quad  x\in \R^{n}.
\]
Clearly, for all  $ f\in \mathscr{S}'$ it holds 
\begin{equation}
D^\alpha \dot{P}_{ x_0, r}^N(f) =\dot{P}_{ x_0, r}^{ N-| \alpha |}(D^{ \alpha } f)\quad  \forall\,| \alpha | \le k.
\label{5.12a}
\end{equation}

\begin{thm}
\label{thm5.4} 
1. For every $ N\in \N_0$ the following embedding holds true.
\begin{equation}
\begin{cases}
 \mathscr{L}^{N}_{1 (p,  N)}( \R^{n}) \hookrightarrow  C^{N-1, 1}(\R^{n})\quad  &\text{ if}\quad N \ge 1
\\[0.3cm]
 \mathscr{L}^{0}_{1 (p, 0)}( \R^{n}) \hookrightarrow  L^\infty(\R^{n})\quad &\text{ if}\quad N = 0.
\end{cases}
\label{5.13}
\end{equation}

2. For every $ f\in  \mathscr{L}^{N}_{1 (p,  N)}( \R^{n})$ there exists a unique $ \dot{P} ^N_{ \infty}\in \dot{\mathcal{P}} _N$ such that 
for all $ x_0 \in \R^n$
\[
\lim_{r \to \infty}  \dot{P} _{ x_0, r}(f) \rightarrow \dot{P}^N_{ \infty}(f)\quad  \text{ in }\quad \mathcal{P}_N. 
\]
Furthermore, $\dot{P}^N_{ \infty}:  \mathscr{L}^{N}_{1 (p,  N)}( \R^{n}) \rightarrow \dot{\mathcal{P}} _N$ is a projection with the property 
\begin{equation}
D^\alpha \dot{P}^N_{ \infty}(f) = \dot{P}^{ N-| \alpha |}_\infty (D^\alpha f) \quad  \forall\,| \alpha | \le N. 
\label{5.14}
\end{equation}
3. For all $ g, f\in \mathscr{L}^{1}_{1 (p,  1)}( \R^{n})$ it holds
 \begin{align}
 \dot{P}^1_{ \infty}(g \partial _k f) &= 
 \dot{P}^1_{ \infty}(g) \partial _k\dot{P}^1_{ \infty} (f) = \dot{P}^1_{ \infty}(g) \dot{P}^0_{ \infty} (\partial _k f), \quad 
 k=1, \ldots, n.
 \label{5.15}
\end{align}
In addition, for $ g\in C^{ 0,1}(\R^{n}; \R^{n}  )$, and for all 
$ f\in \mathscr{L}^{0}_{1 (p,  0)}( \R^{n})$  it holds 
\begin{align}
  \dot{P}^0_{ \infty}(g\partial _k  f) := \lim_{r \to \infty}  P^0_{ 0, r}(g\partial _k   f) =0,\quad k=1, \ldots, n, 
\label{5.15b}
\end{align}
where $ g\partial _k   f = \partial _k(gf)- \partial _kg f\in \mathscr{S}'$. 

\vspace{0.3cm}
4. For all $ v\in \mathscr{L}^{1}_{1 (p,  1)}( \R^{n}; \R^{n} )$ with $ \nabla \cdot  v=0$ almost everywhere in $ \R^{n}$  
and $ f\in \mathscr{L}^{1}_{1 (p,  1)}( \R^{n})$ it holds 
\begin{align}
 \dot{P}^0_{ \infty}(\nabla v \cdot \nabla f) &= 
 \dot{P}^0_{ \infty}(\nabla v)\cdot \dot{P}^0_{ \infty}(\nabla f). 
\label{5.15a}
\end{align}

\end{thm}

\vspace{0.3cm}
 Next, we have the following norm equivalence, which is similar to the properties of the usual Campanato spaces.  
 
 \begin{lem}
 \label{lem5.9}  Let  $1 \le p < +\infty, 1 \le q \le +\infty$, and $ N, N'\in \N_0, N< N',  s \in [- \frac{n}{p}, N+1)$. 
  If $ f\in \mathscr{L}^{k,s}_{q  (p, N')}(\R^{n})$, and satisfies   
 \begin{equation}
   \lim_{m \to \infty}\dot{P}^{ L}_{ 0, 2^m}(D^k f) = 0\quad \forall\,L=N+1, \ldots, N'.
 \label{5.23}
 \end{equation} 
then $ f\in \mathscr{L}^{k,s}_{q (p, N)}(\R^{n})$ and it holds, 
 \begin{equation}
  | f|_{  \mathscr{L}^{ k,s}_{q(p, N')}} \le | f|_{  \mathscr{L}^{ k,s}_{ q(p,N)}} \le c  | f|_{  \mathscr{L}^{ k,s}_{q(p, N')}}.
 \label{5.24}
 \end{equation}

 \end{lem} 
 
\begin{rem}
\label{rem5.7} 
For all   $ f\in \mathscr{L}^{s}_{q (p, N)}(\R^{n}), 1 \le p < +\infty, 1 \le q \le +\infty, s\in [- \frac{n}{p}, N+1)$,  the condition  \eqref{5.23} is fulfilled, and therefore 
\eqref{5.24} holds for all $ f\in \mathscr{L}^{s}_{q (p, N)}(\R^{n})$ under the assumptions on $ p,q,s, N$ and $ N'$ 
of Lemma\, \ref{lem5.9}.  To verify this fact we observe  that   for  $ f\in \mathscr{L}^{s}_{q (p, N)}(\R^{n})$\begin{equation}
\sup_{ m\in \Z} 2^{ -Nm} \osc_{ p, N}(f, 0, 2^m) \le |f|_{ \mathscr{L}^{s}_{q (p, N)}}. 
\label{5.29}
\end{equation} 
Then for $ L \in \N$, $ L>N$,  we estimate for all multi index $ \alpha $ with $ |\alpha |=L$
\begin{align*}
  |D^\alpha \dot{P}_{0, 2^m}^L(f)| &= |D^{ \alpha }\dot{P}_{0, 2^m}^L((f- P^N_{ 0, 2^m})) | 
  \le c2^{ -Lm}\osc_{ p, N}(f, 0, 2^m) 
\\
&\le c 2^{m (N-L)} |f|_{ \mathscr{L}^{s}_{q (p, N)}}
\rightarrow 0\quad   \text{as}\quad m \rightarrow +\infty.   
\end{align*}
Hence, \eqref{5.23} is fulfilled.  
\end{rem}

\begin{rem}
 \label{rem5.7b} 
 In case $ q=\infty$,  since  $ \mathscr{L}^{s}_{\infty (p, N)}(\R^{n} ) $ coincides with the 
 usual Campanato space,  and Lemma\,\ref{lem5.9}  reduces to the well known result(cf. \cite[p. 75]{gia}).  
 
\end{rem}

\vspace{0.3cm}
We also have the  following growth properties  of functions in $\mathscr{L}^{s}_{q (p, N)}(\R^{n} ) $ as $ |x| 
\rightarrow +\infty$ (see \cite{trans}). 

\begin{lem}
 \label{growth}
 Let $ N\in \N_0$. Let  $f\in  \mathscr{L}^{s}_{q (p, N)}(\R^{n} ) , 1 \le q \le +\infty, 1 \le p < +\infty, 
 s \in  [N, N+1)$. 
 
 1. In case $ s\in (N, N+1)$ it holds
 \begin{equation}
  |f(x)| \le c (1+|x|^s) \|f\|_{ \mathscr{L}^{s}_{q (p, N)}}\quad \forall x\in \R^{n}. 
 \label{growth1}
  \end{equation} 
  
  2. In case $ s=N$ it holds
 \begin{equation}
  |f(x)| \le c (1+\log(1+|x|)^{ \frac{1}{q'}}|x|^{ N}) \|f\|_{ \mathscr{L}^{N}_{q (p, N)}}\quad \forall x\in \R^{n}. 
 \label{growth11}
  \end{equation} 
  Here  $q'= \frac{q}{q-1}$, and the constant $ c=\const>0$, depends on $ q,p,s, N$ and $ n$. 
 
\end{lem}

\section{Calder\'on-Zygmund estimate involving $\mathscr{L}^{ s}_{ q ( p, N)} $ norm}
\label{sec:-6}
\setcounter{secnum}{\value{section} \setcounter{equation}{0}
\renewcommand{\theequation}{\mbox{\arabic{secnum}.\arabic{equation}}}}

In this section we establish the  Calder\'on-Zygmund type estimate  for our spaces.  For this purpose let us introduce the partition 
of unity, which will be used in what follows. We set $ U_j = B(2^{ j+1})  \setminus \overline{B(2^{ j-1})} , j\in \Z$. Clearly, 
$ \{ U_j\}$ is a local finite covering of $ \R^{n}  \setminus \{ 0\}$. By $\{ \psi _j\}$ we denote a corresponding partition 
of unity of radial symmetric functions $ \psi _j \in C^\infty_c(U_j)$, such that $ 0 \le \psi _j \le 1$ in $ U_j$, $ | D^k \psi _j| 
\le c_k 2^{ -k j}$ in $ U_j$ and $ \sum_{j\in \Z} \psi _j=1$ on $ \R^{n}  \setminus \{ 0\}$.  We have the following.  

\begin{lem}
\label{lem6.1}
Let $ K\in C^{ 2}( \R^{n}  \setminus \{ 0\})$ be a Calder\'on-Zygmund kernel, i.e. 
\begin{itemize}
  \item[(i)]  $ K(x) \sim t ^{ n}K(t x)$ for all $ x\in \R^{n}  \setminus \{ 0\}, t >0$.
  \item[(ii)] $ \intl_{\partial B(1)} K(x)  dS =0$. 
\end{itemize} 
By  $ \{ \psi _j\}$ we denote a partition of unity introduced above. Let $ m, k \in \Z, m < k$. Define 
\[
T_{m}^k (h)(x) = \sum_{i=m}^{k} \intl_{ \R^{n}} h(x-y) K(y) \psi _i (y) dy,\quad x\in \R^{n},\quad   h \in L^1_{ loc}( \R^{n}),  
\]
then for all $1< p< +\infty, 1 \le q \le +\infty, s \in [0, N+1), N\in \N_0  $ the operator $ T_{m}^{ k}$ is uniformly  bounded in  
$ \mathscr{L}^{ s}_{ q (p, N)}( \R^{n})$, i.e. it holds 
\begin{equation}
| T_{m}^k(h)|_{ \mathscr{L}^{ s}_{ q (p, N)}} \le c_{ N,  s , q, n} | h|_{\mathscr{L}^{ s}_{ q (p, N)}}\quad \forall\,h\in \mathscr{L}^{ s}_{ q (p, N)}( \R^{n}). 
\label{6.1}
\end{equation}
\end{lem}

{\bf Proof}:
Let $K : \R^{n}  \setminus \{ 0\}$ be a Calder\'on-Zygmund kernel. Let $ h\in \mathscr{L}^{ s}_{ q (p, N)}(\R^{n})$. 
Given $ m, M \in \Z, m < M$ we set 
\[
f_{ m}^k(x) = T_m^k(h)(x)= \sum_{i = m}^{k} \intl_{ \R^{n}}  h(x-y) K(y) \psi _i (y) dy,\quad  x\in \R^{n}.   
\]

Let $ x_0 \in \R^{n}$ be arbitrarily chosen. Let $ j\in \Z$ be fixed. Our aim will be to evaluate the 
$ \osc_{ p,N}(f^k_m; x_0, 2^{ j})$. 
We decompose  $ f_{ m}^k$ into the sum $g_{ m}^k + G_{ m}^k$ by means of 
\begin{align*}
g_{ m}^k(x)  &= \sum_{\substack{i = m \\ i \le j}}^{k} \intl_{ \R^{n}}  h(x-y)  K(y) \psi _i (y) dy,
\\
G_{ m}^k(x)  &= \sum_{\substack{i = m \\ i > j}}^{k} \intl_{ \R^{n}}  h(x-y) K(y) \psi _i (y) dy,\quad  x\in \R^{n}.   
\end{align*}

Let $ x \in B(x_0, 2^j)$ arbitrarily chosen, but fixed. Defining  $ Q\in \mathcal{P}_N$ by means of 
\[
Q(x)=\sum_{\substack{i = m \\ i \le j}}^{k} \intl_{ \R^{n}}  P^{ N}_{ x_0, 2^{ j+2}}(h)(\cdot -y) K(y) \psi _i (y) dy,\quad  x \in \R^{n},
\]
 it can be checked easily that     
\begin{align*}
g_{ m}^k(x) - Q(x) &= \sum_{\substack{i = m \\ i \le j}}^{k} \intl_{B(2^{ j+1})}  (h(x-y)- P^{ N}_{ x_0, 2^{ j+2}}(h)(x-y)) K(y) \psi _i (y) dy
\\
&=\sum_{\substack{i = m \\ i \le j}}^{k} \intl_{B(x, 2^{ j+1})}  (h(y)- P^{ N}_{ x_0, 2^{ j+2}}(h)(y)) K(x-y) \psi _i (x-y) dy.
\end{align*}

Clearly, $ B(x, 2^{ j+1}) \subset B(x_0, 2^{ j+2})$ and $ \supp(\psi _i(x-\cdot ))  \subset B(x, 2^{ j+1})$ for all $ i \le j$.  
This shows that 
\begin{align*}
g_{ m}^k(x) - Q(x) 
&=\sum_{\substack{i = m \\ i \le j}}^{k}\intl_{ \R^{n} }  \chi _{ B(x_0, 2^{ j+2})}(y)(h(y)- P^{ N}_{ x_0, 2^{ j+2}}(h)(y)) K(x-y) \psi _i (x-y) dy
\\
&=\sum_{\substack{i = m \\ i \le j}}^{k} \intl_{ \R^{n}}  
\Big(\chi _{ B(x_0, 2^{ j+2})}(h- P^{ N}_{ x_0, 2^{ j+2}}(h))\Big)(x-y)) K(y) \psi _i (y) dy. 
\end{align*}
By virtue of the well known Cader\'on-Zygmund inequality in $ L^p$ we find 
\begin{equation}
\| g_{ m}^k(x)-Q\|_{ L^p(B(x_0, 2^j))} \le c \| h- P^{ N}_{ x_0, 2^{ j+2}}(h)\|_{L^p(B(x_0, 2^{ j+2})) }
\le c 2^{ j\frac{n}{p}}\osc_{ p,N}(h; x_0, 2^{ j+2}).  
\label{6.2}
\end{equation}
Next, we estimate $  G_{ m}^k$. Let $ \alpha \in \N_0^n$ be any multi index with $ | \alpha |= N+1$. Clearly, 
\begin{align*}
D^\alpha  G_{ m}^k(x) =  \sum_{\substack{i = m \\ i > j}}^{k} \intl_{ \R^{n}}  (h(x-y) -P^{ N}_{ x_0, 2^{ i+2}}(h)(x-y)) D^\alpha (K(y)  \psi _i (y)) dy.
\end{align*}
 Let  $ i \in \{j+1, \ldots, k\}$. Noting that $ B(x, 2^{ i+1}) \subset  B(x_0, 2^{ i+1}+ 2^{ j}) \subset B(x_0, 2^{ i+2})$, and employing 
Jensen's inequality,   
we estimate 
\begin{align*}
&\bigg|\intl_{ \R^{n}}  (h(x-y) -P^{ N}_{ x_0, 2^{ i+2}}(h)(x-y)) D^\alpha (K(y)  \psi _i (y)) dy\bigg|
\\
&\quad \le c\intl_{ B(2^{ i+1})}  | h(x-y) -P^{ N}_{ x_0, 2^{ i+2}}(h)(x-y)| 2^{ -i (n +N+1)} dy
\\
&\quad = c\intl_{ B(x, 2^{ i+1})}  | h(y) -P^{ N}_{ x_0, 2^{ i+2}}(h)(y)| 2^{ -i (n +N+1)} dy
\\
&\quad \le  c\intl_{ B(x_0, 2^{ i+2})}  | h(y) -P^{ N}_{ x_0, 2^{ i+2}}(h)(y)| 2^{ -i (n +N+1)} dy
\\
&\quad \le  c 2^{  -(i+2) (\frac{n}{p} + N+1)} \bigg(\intl_{ B(x_0, 2^{ i+2})}  | h(y) -P^{ N}_{ x_0, 2^{ i+2}}(h)(y)|^p dx\bigg)^{ \frac{1}{p}}
\\
& \quad \le  c 2^{  -(i+2) (N+1)} \osc_{ p, N}(h; x_0, 2^{ i+2}). 
\end{align*}
Summing over $ i= j+1$ to $ k$ to both sides of the above inequality and multiplying the result by $ 2^{ j(N+1)}$, we get  
\begin{equation}
2^{ j (N+1)}\| D^{ N+1} G_{ m}^k(x) \|_{ L^\infty(B(x_0, 2^j))}  \le cS_{ N+1,1}(\osc_{ p, N}(h; x_0))_j. 
\label{6.3}
\end{equation}
Thanks to Poincar\'e's inequality \eqref{6.3} implies 
\begin{align}
\|  G_{ m}^k- P^{ N}_{ x_0, 2^j}(G_{ m}^k) \|_{ L^p(B(x_0, 2^j))}  
&\le c 2^{ j ( \frac{n}{p} +N+1) } \| D^{ N+1} G_{m}^k(x) \|_{ L^\infty(B(x_0, 2^j))}
\cr
&\le c 2^{ j \frac{n}{p} }S_{ N+1,1}(\osc_{ p,N}(h; x_0))_j. 
\label{6.4}
\end{align}
Furthermore, noting that $ P_{ x_0, 2^j}^{ N}(f^k_m) = P_{ x_0, 2^j}^{ N}(g^k_{ m})+ P_{ x_0, 2^j}^{ N}(G^k_{m})$, 
we infer 
\begin{align*}
&\|  f_{ m}^k- P^{ N}_{ x_0, 2^j}(f_{ m}^k) \|_{ L^p(B(x_0, 2^j))}  
\\
&\quad \le \|  g_{ m}^k- P^{ N}_{ x_0, 2^j}(g_{ m}^k) \|_{ L^p(B(x_0, 2^j))}   + \|  G_{ m}^k- P^{ N}_{ x_0, 2^j}(G_{ m}^k) \|_{ L^p(B(x_0, 2^j))}  
\\
&\quad \le c\| g_{ m}^k- Q\|_{ L^p(B(x_0, 2^j))}   + \|  G_{ m}^k- P^{ N}_{ x_0, 2^j}(G_{ m}^k) \|_{ L^p(B(x_0, 2^j))}.  
\end{align*}
Combining this inequality with \eqref{6.2} and \eqref{6.4},  we obtain 
\begin{align}
&\osc_{ p, N}(f_m^k; x_0, 2^{ j})
\cr
& \le c\osc_{p,  N}(h; x_0, 2^{ j+2})  + cS_{ N+1,1}(\osc_{ p, N}(h; x_0))_j
\cr
& \le  c S_{ N+1,1}(\osc_{p, N}(h; x_0))_j. 
\label{6.5}
\end{align}
We now perform  $ S_{ s, q}$ to the both sides of \eqref{6.5}, and use Lemma\,\ref{lem10.1} with $ {\D X=\osc_{ p, N}(h; x_0)}$, 
with $ p=1, \alpha =N+1, \beta =s $. Then taking the supremum over $ x_0\in \R^{n}$ on both sides, we obtain 
\begin{equation}
| f_{ m}^k|_{ \mathscr{L}^{ s}_{q (p, N)}} \le c | h|_{ \mathscr{L}^{ s}_{q (p, N)}}.   
\label{6.6}
\end{equation}  
Whence, \eqref{6.1}.  \hfill \Beweisende

\begin{lem}
\label{lem6.2}
Let $ N\in \N_0$. Let $ \{ \varphi _k\}$ be a sequence of functions in $ C^{l}_c (\R^{n}), l\in \N, l \ge N+1,$ 
with $ \supp(\varphi_k ) \subset \R^{n}  \setminus B(2^k)$, and 
\begin{equation}
| D^\alpha \varphi_k | \le c | y|^{ -n-| \alpha |}\quad  \text{ in}\quad \R^{n}\quad  \forall\,| \alpha | \le l. 
\label{6.10a}
\end{equation}
Let $1 \le q \le +\infty,  s\in  [0, N+1)  $. Given $f\in  L^p_{ loc}(\R^{n})$ such that 
\begin{equation}
  \sup_{x_0\in \R^{n} } \sum_{j=1}^{\infty} 2^{- jsq } \osc_{ p,N}(f; x_0, 2^j)^q <+\infty.
\label{6.10aa}
 \end{equation} 
Define, 
\begin{align*}
w^k(x) = \intl_{ \R^{n}} f(x-y) \varphi _k(y) dy, \quad  x\in  \R^{n}, \quad  k\in \N.   
\end{align*}
Then for every multi index $ \alpha $ with  $ |\alpha | \in \{N+1, \ldots, l\}$, 
\begin{equation}
D^{ \alpha  }  w^k \rightarrow 0 \quad  \text{{\it uniformly in}}\quad \R^{n}\quad  \text{{\it as}}\quad  k \rightarrow +\infty. 
\label{6.10b}
\end{equation}
\end{lem}
\vspace{0.3cm}

{\bf Proof}: Let $ \{ \psi _j\}$ denote the partition of unity introduced in the beginning of this section. 
Since $ \supp(\varphi _k) \subset \R^{n}  \setminus B(2^k)$ we get 
\begin{align*}
w^k(x) = \sum_{j=k-1}^{\infty}\intl_{ \R^{n}} f(x-y) \varphi _k(y) \psi _j(y) dy, \quad  x\in  \R^{n}, \quad  k\in \N.   
\end{align*}
Let $ \alpha \in \N_0^n$ be any multi index with $ | \alpha |\in \{N+1, \ldots, l\}$. We calculate 
\begin{align*}
D^\alpha  w^k(x) =  \sum_{j = k-1}^{\infty} \intl_{ \R^{n}}  (f(x-y) -P^{ N}_{ x_0, 2^{ j+2}}(x-y)) D^\alpha (\varphi _k(y)  \psi _i (y)) dy,\quad x\in \R^{n}. 
\end{align*}
Let  $  j\in \Z,  j \ge k-1 $. Fix $ x_0\in \R^{n}$.  Noting that $ B(x, 2^{ j+1}) \subset  B(x_0, 2^{ j+1}+ 2^{ k}) \subset B(x_0, 2^{ j+2})$  for all $ x\in B(x_0, 2^k)$, observing  \eqref{6.10a}, and arguing as in the proof 
of Lemma\,\ref{lem6.1}, we find for all $ x\in B(x_0, 2^k)$, 
\begin{align*}
&\bigg|\intl_{ \R^{n}}  (f(x-y) -P^{ N}_{ x_0, 2^{ j+2}}(x-y)) D^\alpha (\varphi _k(y)  \psi _j (y)) dy\bigg|
\\
&\quad \le c\intl_{ B(2^{ j+1})}  | f(x-y) -P^{ N}_{ x_0, 2^{ j+2}}(x-y)| 2^{ -j (n +|\alpha |)} dy
\\
& \quad \le  c 2^{  -j |\alpha |} \osc_{p, N}(f; x_0, 2^{ j+2}). 
\end{align*}
This  together with H\"older's inequality yields
\begin{align}
|D^{ \alpha } w^k(x_0)| &\le \| D^{ \alpha } w^k\|_{ L^\infty(B(x_0, 2^k))}  
\cr
&\le c\sum_{j=k-1}^{\infty}
2^{  -j |\alpha |} \osc_{ p, N}(f; x_0, 2^{ j+2})
 \le  c\sum_{j=k}^{\infty} 2^{  -j (|\alpha |-s)} 2^{ -js}\osc_{p, N}(f; x_0, 2^{ j})
\cr
&\le  c \Big(\sum_{j=k}^{\infty} 2^{-j (|\alpha |-s) q'}\Big)^{ \frac{1}{q'}} 
\Big(\sup_{x\in \R^{n} } \sum_{j=1}^{\infty} 2^{- jsq } \osc_{ p,N}(f; x, 2^j)^q \Big)^{ \frac{1}{q}}. 
\label{6.10c}
\end{align}
Observing  \eqref{6.10aa}, the right-hand side tends to zero as $ k \rightarrow +\infty$ uniformly in $ x_0\in \R^{n}$ we get the claim.  
\hfill \Beweisende

\vspace{0.3cm}
Next, we apply  Lemma\,\ref{6.1} to the Laplace equation 
\begin{equation}
- \Delta f = \sum_{\alpha ,\beta =1}^{n} \partial _\alpha  \partial _\beta  H_{ \alpha \beta }\quad  \text{ in}\quad  \R^{n}. 
\label{6.7}
\end{equation}
Let $ H\in L^2_{ loc}(\R^{n}; \R^{n^2})$. A function $ f\in L^2_{ loc}(\R^{n})$ is called a {\it very weak solution} to \eqref{6.7} if the following identity holds for all $ \phi \in C^{\infty}_c(\R^{n})$   
\begin{equation}
- \intl_{ \R^{n}} f \Delta \phi   dx  = \intl_{ \R^{n}} H:  \nabla ^2 \phi   dx.  
\label{6.8}
\end{equation}

Below we shall  also make use of the following.

\begin{lem}
 \label{lem6.2a}
 Let $ \{f_k\}$ be a bounded sequence in $ \mathscr{L}^{ s}_{q (p, N)}(\R^{n} )$. Suppose there exists $ f\in L^p_{ loc}(\R^{n} )$ such that  
 \begin{equation}
  f_k \rightarrow f\quad  \text{in}\quad L^p_{ loc}(\R^{n} )\quad   \text{as}\quad k \rightarrow +\infty. 
 \label{6.8a}
  \end{equation} 
 Then $ f\in \mathscr{L}^{ s}_{q (p, N)}(\R^{n} )$ and it holds
 \begin{equation}
|f|_{ \mathscr{L}^{ s}_{q (p, N)}} \le \sup_{ k\in \N} |f_k|_{ \mathscr{L}^{ s}_{q (p, N)}}.   
 \label{6.8b}
  \end{equation} 
\end{lem}
{\bf Proof}: Let $ m,l\in \Z, m< l$ be arbitrarily chosen. By means of  \eqref{6.8a} we get for all $ x_0\in \R^{n} $
\begin{align*}
 \sum_{j=m}^{l} \osc_{ p,N} (f; x_0, 2^j) =  \lim_{k \to \infty}\sum_{j=m}^{l} \osc_{ p,N} (f_k; x_0, 2^j)   
\le \sup_{ k\in \N} |f_k|_{ \mathscr{L}^{ s}_{q (p, N)}}. 
\end{align*}
Passing $ m \rightarrow -\infty$ and $ l \rightarrow +\infty$ and taking the supremum over all $ x_0 \in \R^{n} $, we obtain the claim  \eqref{6.8b}.  \hfill \Beweisende

\vspace{0.2cm}
We have following 

\begin{thm}
\label{thm6.3}
Let $ N\in \{ 0,1\}$, $ s\in [0, N+1)$. 
For  every $ H\in \mathscr{L}^{ s}_{q (p, N)}(\R^{n}; \R^{n^2})$ there exists a unique very weak solution 
$f\in \mathscr{L}^{ s}_{q (p, N)}(\R^{n}) $ to \eqref{6.7} such that 
\begin{equation}
P^{ N}_{ 0,1}(f)=0. 
\label{6.9}
\end{equation} 
In particular, the following estimate holds true 
\begin{equation}
| f|_{ \mathscr{L}^{ s}_{q (p, N)}} \le  c| H|_{ \mathscr{L}^{ s}_{q (p, N)}},
\label{6.10}
\end{equation} 
where $ c= \const>0$, depending only on $ n, q, N$ and $ s$. 

\end{thm} 

{\bf Proof}: By $ K_{ \alpha \beta }$ we denote the kernel $ \partial _{ \alpha }\partial _\beta  \Gamma $, where $ \Gamma $ stands for the Newtonian  
potential in $ \R^{n}$, i.e.
\[
\Gamma (x)=\begin{cases}
\dfrac{1}{n | B(1)| | x|^{ n-2}}\quad  &\text{ if}\quad  n \ge 3
\\[0.3cm]
\dfrac{1}{2\pi }\log | x|\quad  &\text{ if}\quad  n =2. 
\end{cases}
\]
It is readily seen that  $ K_{ \alpha \beta }$ is a Calder\'on-Zygmund kernel. Let $H\in \mathscr{L}^{ s}_{q (p, N)}(\R^{n}; \R^{n^2})$. 
We define, for $ m, k\in \Z, m <k$,  
\[
f_{m}^k(x) =  \sum_{i=m}^{k} \intl_{ \R^{n}} H(x-y):  K(y) \psi _i (y) dy,\quad x\in \R^{n}. 
\]
According to Lemma\,\ref{lem6.1} it holds 
\begin{equation}
| f_{ m}^k|_{ \mathscr{L}^{ s}_{q (p, N)}} \le  c| H|_{ \mathscr{L}^{ s}_{q (p, N)}}, 
\label{6.11}
\end{equation}  
where the constant $ c>0$ does not depend on $ m, k\in \Z$. 

\vspace{0.3cm}
1. Assume $ H \in C^\infty(\R^{n})$. Using integration by parts,  and noting that 
\[
 \sum_{i=m}^{m+2} \intl_{B(2^{ m+3})}  \partial _\beta \Gamma (y) \partial _\alpha  \psi _i (y) dy=0, 
\] 
we get 
 \begin{align*}
f_{m }^k(x) &=  \sum_{i=m}^{k} \intl_{ \R^{n}} \partial_{ \alpha } H_{ \alpha \beta }(x-y) \partial _{ \beta }
\Gamma(y) \psi _i (y) dy
 \\
 &\qquad - \sum_{i=k-2}^{k} \intl_{B(2^{ k+1})} H_{ \alpha \beta }(x-y)   \partial _\beta \Gamma (y) \partial _\alpha  \psi _i (y) dy
 \\
 &\qquad  -
\sum_{i=m}^{m+2} \intl_{B(2^{ m+3})}  (H_{ \alpha \beta }(x-y)  -H_{ \alpha \beta }(x))    \partial _\beta \Gamma (y) \partial _\alpha  \psi _i (y) dy. 
 \end{align*}
 Clearly,  from the above  identity 
 we deduce that  $ f_m^k(x) \rightarrow f^k(x)$ as $ m \rightarrow -\infty$ for all $ x\in \R^{n}$, where 
 \begin{align*}
f^k(x) &=   \sum_{i=-\infty}^{k} \intl_{ \R^{n}} \partial_{ \alpha } H_{ \alpha \beta }(x-y) \partial _{ \beta }\Gamma (y) \psi _i (y) dy
  \\
 &\qquad  -
\sum_{i=k-2}^{k} \intl_{B(2^{ k+1})} H_{ \alpha \beta }(x-y)    \partial _\beta \Gamma (y) \partial _\alpha  \psi _i (y) dy. 
 \end{align*} 
 By Lebesgue's theorem of dominated 
 convergence we see that $ f_m^k \rightarrow f^k $ in $ L^p_{ loc}(\R^{n})$ as $ m \rightarrow -\infty$. Applying integration by 
 parts, we find that 
 \begin{align*}
f^k(x) &=   \intl_{ \R^{n}} \partial_{ \alpha } \partial _{ \beta }H_{ \alpha \beta }(x-y)\Gamma(y)\chi _k(y)  dy
  \\
  &\qquad -\sum_{j=k-2}^{k} \intl_{ \R^{n}}  H_{ \alpha \beta }(x-y) \Gamma(y) \partial _{ \alpha }\partial _\beta \psi _j (y) dy
  \\
 &\qquad  -
\sum_{j=k-2}^{k} \intl_{ \R^{n}} H_{ \alpha \beta }(x-y) ( \partial _\beta \Gamma(y) \partial _\alpha  \psi _j (y)+ 
\partial _\alpha  \Gamma(y) \partial _\beta   \psi _j (y) )dy ,
 \end{align*}  
 where 
 \[
 \chi _k(y))=\begin{cases}
\sum_{j=-\infty}^{k} \psi _j(y) \quad   &\text{if}\quad y\in \R^{n} \setminus \{0\},
\\[0.3cm]
1\quad  &\text{if}\quad y=0. 
\end{cases}
\]

 In addition,  \eqref{6.8b} together with \eqref{6.11} implies 
 \begin{equation}
| f^k|_{ \mathscr{L}^{ s}_{q (p, N)}} \le  c| H|_{ \mathscr{L}^{ s}_{q (p, N)}}. 
\label{6.12}
\end{equation}

Set $ g^k:= f^k- P^{ N}_{ 0,1}(f^k)$. We easily  get for $ j\in \N$
\begin{align*}
g^k &= f^k - P^{ N}_{ 0, 2^j}(f^k) + \sum_{i=1}^{j}(P^{ N}_{ 0, 2^i}(f^k)- P^{ N}_{ 0, 2^{ i-1}}(f^k) )
\\
& = f^k - P^{ N}_{ 0, 2^j}(f^k) + \sum_{i=1}^{j}P^{ N}_{ 0, 2^i}(f^k- P^{ N}_{ 0, 2^{ i-1}}(f^k) ). 
\end{align*}
Thus, by the triangle inequality along with \eqref{5.4} we find 
\begin{align*}
\| g^k\|_{ L^p(B(2^j))} &\le c 2^{ \frac{n}{p}j } \osc_{ p, N}(f^k, 0, 2^j)
+  \sum_{i=1}^{j} \| P^{ N}_{ 0, 2^i}(f^k- P^{ N}_{ 0, 2^{ i-1}}(f^k) )\|_{ L^p(B(2^j))}
\\
& \le c 2^{ \frac{n}{p}j } \osc_{p, N}(f^k, 0, 2^j) 
+ c \sum_{i=1}^{j}  2^{ \frac{n}{p} i}\osc_{ p, N}(f^k, 0, 2^i ). 
\end{align*}
Applying H\"older's inequality, and using \eqref{6.12}, we find   
\[
\| g^k\|_{ L^p(B(2^j))} \le c 2^{js+  j\frac{n}{p} }  | f^k|_{ \mathscr{L}^{ s}_{q (p, N)}} \le 
 c2^{js+  j\frac{n}{p} }| H|_{ \mathscr{L}^{ s}_{q (p, N)}}. 
\]
Thus, $ \{ g^k\}$ is bounded in $ L^p_{ loc}(\R^{n})$. Noting that for all $ x_0\in \R^{n}$  it 
holds ${\D  \osc_{ p, N}(g^k, x_0) = \osc_{ p, N}(f^k, x_0)}$, owing to \eqref{6.12},  we infer 
 \begin{equation}
| g^k|_{ \mathscr{L}^{ s}_{q (p, N)}} \le  c| H|_{ \mathscr{L}^{ s}_{q (p, N)}}. 
\label{6.13}
\end{equation}    
Thus, by the compact embedding $ \mathscr{L}^{ s}_{q (p, N)}(\R^{n} ) \hookrightarrow  L^p_{ loc}(\R^{n})$, 
eventually  passing to a subsequence, we get a function 
$ f\in L^p_{ loc}(\R^{n})$ such that 
\begin{equation}
 g^k \rightarrow f  \quad  \text{{\it in}}\quad  L^p_{ loc}(\R^{n})\quad  \text{{\it as}}\quad  m \rightarrow +\infty.  
\label{6.17}
\end{equation}
This together with  \eqref{6.8b} and \eqref{6.13} shows that $ g\in  \mathscr{L}^{ s}_{q (p, N)}(\R^{n})$, and satisfies the inequality 
\begin{equation}
| g|_{ \mathscr{L}^{ s}_{q (p, N)}} \le  c| H|_{ \mathscr{L}^{ s}_{q (p, N)}}. 
\label{6.18}
\end{equation}

Setting 
\begin{align*}
\varphi _k &= - \sum_{j=k-2}^{k} \Big(\Gamma(y) \partial _{ \alpha }\partial _\beta \psi _j (y)+ 
\partial _{ \alpha } \Gamma(y) \partial _\beta \psi _j(y) +
\partial _{ \beta } \Gamma(y) \partial _{ \alpha } \psi _j (y)\Big),
\end{align*}
we may write $ f^k = f^k_1+ f^k_2$, where 
\begin{align*}
f^k_1(x) &=   \intl_{ \R^{n}} \partial_{ \alpha } \partial _{ \beta }H_{ \alpha \beta }(x-y)\Gamma(y) \chi _k(y)  dy
\\
f^k_2(x) &=   \intl_{ \R^{n}} H_{ \alpha \beta }(x-y)\varphi _k(y)  dy,\quad  x\in \R^{n}. 
\end{align*}
Clearly,  $ \varphi _k \in C^{ N+1}_c(\R^{n})$ with $ \supp(\varphi _k) \subset \R^{n}  \setminus B(2^{ k-3})$ and satisfying 
condition \eqref{6.10a} of Lemma\,\ref{lem6.2}. Thus, thanks to Lemma\,\ref{lem6.2}
\begin{equation}
 D^{ N+1} f^k_2 \rightarrow 0  \quad  \text{{\it uniformly in}}\quad  \R^{n}\quad  \text{{\it as}}\quad  k \rightarrow +\infty. 
\label{6.20}
\end{equation} 
Let $ \phi \in C^{\infty}_c(\R^{n})$ be arbitrarily chosen. 
Employing \eqref{6.20}, recalling that $ N \le 1$, we immediately  verify that  
\begin{equation}
\lim_{ k\to \infty}\intl_{ \R^{n}} f^k_2(x)   \Delta \phi(x)    dx  =0.
\label{6.21}
\end{equation}
Using Fubini's theorem,  and applying integration by parts, we calculate 
\begin{align*}
&\intl_{ \R^{n}} f^k_1(x)   \Delta \phi(x)    dx  
\\
&= \intl_{ \R^{n}}  \intl_{ \R^{n}} \partial _\alpha \partial _\beta H(y) \Gamma(x-y) \chi_{ k}(x-y)      \Delta \phi(x)   dy dx  
\\
&= \intl_{ \R^{n}}  \intl_{ \R^{n}} \partial _\alpha \partial _\beta H(y) \Gamma(x-y) \chi_{ k}(x-y)      \Delta \phi(x)    dx dy 
\\
&= -\intl_{ \R^{n}}  \intl_{ \R^{n}} \partial _\alpha \partial _\beta H(y)        \phi(y)   dy 
\\
&\quad  +2\intl_{ \R^{n}}  \intl_{ \R^{n}} \partial _\alpha \partial _\beta H(y) \nabla \Gamma(x-y)\cdot  \nabla  \chi_{ k}(x-y)      \phi(x)     dy dx 
\\
&\quad  +\intl_{ \R^{n}}  \intl_{ \R^{n}} \partial _\alpha \partial _\beta H(y)  \Gamma(x-y) \Delta \chi_{ k}(x-y)     \phi(x)    dy dx 
\\
&= \intl_{ \R^{n}}  \intl_{ \R^{n}} \partial _\alpha \partial _\beta H(y)       \phi(y)   dy +I_k+ II_k. 
\end{align*}
Again applying  integration by parts, we infer 
\begin{align*}
I_k &= 2\intl_{ \R^{n}}  \intl_{ \R^{n}}  H(y)   \partial _\alpha \partial _\beta   (\nabla \Gamma(x-y)\cdot \nabla \chi_{ k}(x-y))     \phi(x)    dy dx 
\\
& =2(-1)^{ N-1}\intl_{ \R^{n}}   \partial_\alpha \partial _\beta   \intl_{ \R^{n}}  H(x- y)      
\nabla \Gamma(y)\cdot \nabla    \chi_{ k}(y) dy    \phi(x)dx. 
\end{align*}
Using  Lemma\,\ref{lem6.2}, we get  $ I_k=o(1)$. By a similar reasoning we see that 
$ II_k = o(1)$. By means of this  properties together with  $ \chi _k \rightarrow 1$ uniformly on each ball as $ k \rightarrow +\infty$, 
along with \eqref{6.21},  we deduce that 
\begin{equation}
\lim_{k \to \infty} \intl_{ \R^{n}} f^k(x)    \Delta \phi(x)    dx  = -
\intl_{ \R^{n}}  \intl_{ \R^{n}} \partial _\alpha \partial _\beta H(y)       \phi(y)   dy. 
\label{6.22}
\end{equation}
On the other hand, recalling the definition of $ g^k$,  we see that 
\[
-\intl_{ \R^{n}} f^k(x)    \Delta \phi(x)    dx =- \intl_{ \R^{n}} g^k(x)    \Delta \phi(x)    dx. 
\]
By the aid of \eqref{6.17} letting  $ k \rightarrow +\infty$ on the right-hand side, and using  \eqref{6.22}, we obtain the identity 
\begin{equation}
-\intl_{ \R^{n}} f(x)    \Delta \phi(x)    dx = \intl_{ \R^{n}}  \intl_{ \R^{n}} \partial _\alpha \partial _\beta H(y)       \phi(y)   dy. 
\label{6.23}
\end{equation}
Accordingly,  $f$ is a very weak solution to \eqref{6.8}.  Recalling that $ P^N_{ 0,1}(g_k)=0$ for all $ k\in \N$ thanks to \eqref{6.17} 
it holds   $P^N_{ 0,1}(f)=0$. 

\vspace{0.2cm}
Now, let $ H\in   \mathscr{L}^{ s}_{q (p, N)}(\R^{n})$ be arbitrarily chosen. By $ H_\var $ for $ \var >0$ we denote the standard 
mollification of $H$.  According to Lemma\,\ref{5.3} it satisfies
\begin{equation}
| H_\var |_{ \mathscr{L}^{ s}_{q (p, N)}} \le  c| H|_{ \mathscr{L}^{ s}_{q (p, N)}}. 
\label{6.25}
\end{equation} 
From the previous step we get a solution $ f_\var \in \mathscr{L}^{ s}_{q (p, N)}(\R^{n})$ to  \eqref{6.8} such that $ P^N_{ 0,1}(f_\var )=0$. 
According to \eqref{6.18},  the following  a priori estimate holds
\begin{equation}
| f_\var |_{ \mathscr{L}^{ s}_{q (p, N)}} \le  c| H|_{ \mathscr{L}^{ s}_{q (p, N)}}. 
\label{6.26}
\end{equation} 
Set $g_\varepsilon = f_\varepsilon - P^{ N}_{ 0,1}(f_\varepsilon )$. As above we  verify  that $ \{ g_\var \}$ is bounded in $ L^p_{ loc}(\R^{n})$ 
and in $ \mathscr{L}^{ s}_{q (p, N)}(\R^{n})$. By similar argument to the above  we get a function $ f\in \mathscr{L}^{ s}_{q (p, N)}(\R^{n})$
 together with a sequence $  \var _k \searrow 0$  as $ k \rightarrow +\infty$ such that 
 \begin{equation}
 f_{ \var _k} \rightarrow f  \quad  \text{{\it in}}\quad  L^p_{ loc}(\R^{n})\quad  \text{{\it as}}\quad  k \rightarrow +\infty.  
\label{6.27}
\end{equation}
In addition, it holds
\begin{equation}
| f|_{ \mathscr{L}^{ s}_{q (p, N)}} \le  c| H|_{ \mathscr{L}^{ s}_{q (p, N)}}. 
\label{6.28}
\end{equation}    
Let $ \phi \in C^{\infty}_c(\R^{n})$. 
Since $ f_{ \var _k} $ solves \eqref{6.8} with $ H_{ \var _k}$ in place of $ H$,  we infer that the following identity holds true 
\[
 -\intl_{ \R^{n}} f_{ \var _k}   \Delta \phi   dx =   
\intl_{ \R^{n}}  H_{ \var _k } :  D ^2 \phi   dx.
\]
Letting $ k \rightarrow +\infty$ on both sides of the above identity,  and making use of \eqref{6.27}, we are led to 
\[
- \intl_{ \R^{n}} f\Delta \phi   dx =\intl_{ \R^{n}}  H :  D^2 \phi   dx.
\]
This shows that $ f$ is a very weak solution to \eqref{6.8} satisfying \eqref{6.9}. 

\vspace{0.3cm}
{\it Uniqueness}.  Let $ \overline{f} $ be another very weak solution to \eqref{6.8} satisfying \eqref{6.9}.  Then 
$ f- \overline{f} \in \mathscr{L}^{ s}_{ q (p, N)}(\R^{n})$ and $ f- \overline{f}$ is  harmonic. By the 
virtue of the Caccioppoli inequality for harmonic functions  we get 
\begin{align*}
&\bigg(\intmw_{B(  2^j)} | D^{ N+1} (f-\overline{f} )|^p  dx\bigg)^{ \frac{1}{p}} 
\\
&\qquad \le c2^{ - j (N+1) }    
\bigg(\intmw_{B(2^{ j+1})} | (f-\overline{f} ) - P^{ N}_{ 0, 2^{ j+1}}(f-\overline{f} )|^p  dx\bigg)^{ \frac{1}{p}}
\\
& \qquad \le c2^{ - j (N+1) }   \osc_{ p, N} (f-\overline{f} ,0 , 2^{ j+1}) 
\le c2^{ - j (N+1-s) }  | f-\overline{f}  |_{  \mathscr{L}^{ s}_{ q (p, N)}}.  
\end{align*}
Since the right-hand side tends to zero as $ j \rightarrow +\infty$,  we
deduce   that $ D^{ N+1} (f-\overline{f} )=0$. Hence, $ f-\overline{f} \in \mathcal{P}_{N}$. Observing \eqref{6.9}, it follows that 
$f-\overline{f} =0 $. This completes the proof of the uniqueness. 

 \hfill \Beweisende

We are now in a position to prove Theorem\,\ref{thm1.2}.

\vspace{0.3cm}
{\bf Proof of Theroem\,\ref{thm1.2}}:  
Let $ H\in \mathscr{L}^{ 1}_{ 1 (p,1)}(\R^{n})$. Thanks to Theorem\,\ref{thm6.3} there exists 
a unique very weak  solution $ g \in \mathscr{L}^{ 1}_{ 1 (p,1)}(\R^{n})$ to \eqref{6.7}  such that 
$ P^1_{ 0, 1}(g)=0$.  Let $ a\in \R^{n}$ and $ Q_\infty \in \mathcal{\dot{P} }_1$. We define 
\[
f(x) = g(x)-g(0)+a -\dot{P}^1_{ \infty}(g)(x) + Q_{\infty}(x),\quad  x\in  \R^{n}.   
\]
Clearly, $f\in \mathscr{L}^{ 1}_{ 1 (p,1)}(\R^{n})$ is a very weak solution to \eqref{6.7} satisfying $ f(0)=a$ and 
$ \dot{P}^1_{ \infty}(f)=Q_{\infty}$.  Assume $ \overline{f} $ is another very weak solution to \eqref{6.7} satisfying 
$ \overline{f} (0)=a$ and $ \dot{P}^1_{ \infty}(\overline{f} )=Q_{\infty}$. Then by Weyl's lemma $ f- \overline{f} $ is harmonic. 
Using Liouville's theorem we see that $ f-\overline{f} \in \mathcal{P}_1$. Owing to $ f(0)-\overline{f}(0)=0$ and 
$ \nabla (f-\overline{f} )= \nabla \dot{P}^1_{ \infty}(f- \overline{f})=0$ we get $ f=\overline{f} $.   \hfill \Beweisende

\vspace{0.3cm}
{\it Definition of the Helmholtz-Leray projection. } Let $ u\in \mathscr{L}^{ 1}_{ 1 (p,1)}(\R^{n})$. 
Applying   Theorem\,\ref{thm1.2} with $ a=0$ and $ Q_\infty= - \frac{1}{n}P ^0_\infty(\nabla \cdot u) x_\alpha , 
\alpha =1, \ldots,n$, we get  a unique very weak solution 
 $w= Q^{ \sharp} (u)\in 
 \mathscr{L}^{ 1}_{ 1 (p,1)}(\R^{n})$ to the equation 
 \begin{equation}
 \begin{cases}
 -\Delta w_\alpha   = - \partial _\alpha \nabla \cdot u\quad  \text{ in}\quad  \R^{n}, \quad  \alpha =1, \ldots, n. 
 \\[0.3cm]
 w_\alpha (0)=0,\quad  P ^0_\infty(w_\alpha ) = - \frac{1}{n}P ^0_\infty(\nabla \cdot u) x_\alpha 
 \end{cases}
  \label{6.29}
 \end{equation}
 Then we define $ \PP : \mathscr{L}^{ 1}_{ 1 (p,1)}(\R^{n}) \rightarrow \mathscr{L}^{ 1}_{ 1 (p,1)}(\R^{n})$
 by means of 
\[
\PP  u = u - Q ^{ \sharp} (u),\quad  u\in \mathscr{L}^{ 1}_{ 1 (p,1)}(\R^{n}). 
\]

 Thanks, to Theorem\,\ref{thm1.2}, both $ Q ^{ \sharp}$ and $ \PP $ are bounded operators.  
 It is readily seen that $  \intl_{ \R^{n}}  \PP  u\cdot \nabla \phi dx =0$ 
 for all $ \phi \in C^{\infty}_c(\R^{n})$, which shows that 
$ \nabla \cdot \PP  u=0$. On the other hand, if $ u\in \mathscr{L}^{ 1}_{ 1 (p,1)}(\R^{n})$ with $ \nabla \cdot u=0$ then 
$\Delta  Q ^{ \sharp} (u)=0$ and $ \dot{P}^1_\infty  (Q ^{ \sharp} (u))=0$, which implies that  $Q ^{ \sharp} (u)$ is constant. Observing \eqref{6.29}$ _2$, we 
infer $ Q ^{ \sharp} (u)=0$ and 
\[
Q ^{ \sharp} (u) = \dot{P}^1_\infty  (Q ^{ \sharp} (u)) =0. 
\] 
Accordingly, $ \PP(u)=u$, which shows that $ \PP: \mathscr{L}^{ 1}_{ 1 (p,1)}(\R^{n}) \rightarrow 
\mathscr{L}^{ 1}_{ 1 (p,1)}(\R^{n})$ defines a projection onto divergence free fields.

 \vspace{0.3cm}  
 In what follows we consider the equation \eqref{6.8} for matrices $ H=u  \otimes    v$.  
 We first prove the following lemma,  
 which covers the situation $ \nabla  \cdot u =0$  and  $v=(h, 0, \ldots, 0) $.  
 
 \begin{thm}
 \label{thm6.4}
Let $ N\in \{ 0,1\}$, $ 1< p< +\infty, 1 \le q \le +\infty, s\in [0, N+1)$. 
 Let $ u, v\in \mathscr{L}^{s}_{ q (p, N)}(\R^{n})\cap  \mathscr{L}^{ 1}_{ 1 (p, 1)}(\R^{n})$,   such that 
\begin{equation}
\nabla \cdot u=0\quad  \text{a.e.  in}\quad  \R^{n}. 
\label{6.31}
\end{equation}
Then for every $ l\in \{1, \ldots, n\}$  there exists a unique  solution 
$ f\in \mathscr{L}^{s}_{ q (p, N)}(\R^{n})\cap  \mathscr{L}^{1}_{ 1 (p, 1)}(\R^{n})$  
to the equation 
 \begin{equation}
 -\Delta f = \partial _l\nabla \cdot  ( (v\cdot \nabla) u)\quad  \text{ in}\quad  \R^{n},
 \label{6.32}
 \end{equation}
such that 
\begin{equation}
\dot{P} ^1_{\infty}(f)= -\frac{1}{n}P ^0_\infty(\nabla u : (\nabla v)^{ \top}) x_l, \quad  f(0)=0. 
\label{6.33}
\end{equation}
In addition, it holds
\begin{align}
\begin{cases}
| f|_{  \mathscr{L}^{ 1}_{ 1 (p,1)}}
  \le  c\Big(\| \nabla u\|_{ \infty}
| v|_{ \mathscr{L}^{ 1}_{ 1 (p, 1)}}
+ \| \nabla v\|_{ \infty}
| u|_{ \mathscr{L}^{ 1}_{ 1 (p, 1)}}\Big),
\\[0.3cm]
| f|_{  \mathscr{L}^{s}_{ q (p, N)}}
  \le  c\Big(\| \nabla u\|_{ \infty}
| v|_{ \mathscr{L}^{s}_{ q (p, N)}}
+ \| \nabla v\|_{ \infty}
| u|_{ \mathscr{L}^{s}_{ q (p, N)}}\Big).
\end{cases}
\label{6.34}
\end{align}
In particular, $ \nabla  f\in L^\infty(\R^{n})$,  and it holds 
\begin{equation}
\| \nabla  f\|_{ \infty}
  \le  c\Big\{\| \nabla u\|_{ \infty}
| v|_{ \mathscr{L}^{ 1}_{ 1 (p, N)}}
+ \| \nabla v\|_{ \infty}
| u|_{ \mathscr{L}^{ 1}_{ 1 (p, N)}}\Big\}. 
\label{6.35b}
\end{equation}
In addition, given  $ 1 < r < p$, the following inequality holds for all $ j\in \Z$
 \begin{align}
(\osc_{r, N}( f; x_0))_j
&\le c\Big( |\nabla u|_{ BMO}  + 
\sup_{ i \ge j} |P^0_{ x_0, 2^i}(\nabla u)|\Big) S_{ N+1,1}(\osc_{p,  N}(v; x_0))_j 
\cr
&\qquad + c
\Big(|\nabla v|_{ BMO}  + 
\sup_{ i \ge j} |P^0_{ x_0, 2^i}(\nabla v)| \Big)S_{ N+1,1}(\osc_{ p, N}(u; x_0))_j.
\label{6.41c}
\end{align}

 \end{thm}
 
 {\bf Proof}: 
Let $ l\in \{ 1, \ldots, n\}$.  We define,  for $ m, k\in \Z, m < k$, 
 \[
f_{ m}^{ k}(x) = \sum_{j=m}^{k} \intl_{ \R^{n}}  \partial _\beta u_\alpha (x-y) 
 v_\beta  (x-y)  \partial _l \partial _\alpha   (\Gamma(y) \psi _j(y)) dy, \quad x \in \R^{n}.   
\] 
 Let $ j\in \Z$ be fixed. We decompose $ f_{ m }^k$ into the sum $ g_{ m,}^k+G_{m}^k$, where 
  \begin{align*}
g_{ m }^k(x) &= \sum_{\substack{i=m \\ i \le j}}^{k}  \intl_{ \R^{n}}  \partial _\beta   u_\alpha (x-y)  v_\beta  (x-y)  \partial _l \partial _\alpha   (\Gamma(y) \psi _j(y)) dy, 
   \\
  G_{ m }^k(x) &= \sum_{\substack{i=m \\ i > j}}^{k}  \intl_{ \R^{n}}  \partial _\beta   u_\alpha (x-y)  v_\beta  (x-y)  \partial _l \partial _\alpha   (\Gamma(y) \psi _j(y)) dy, \quad y \in \R^{n}.   
  \end{align*}
Setting  
 \[
Q (x)=\sum_{\substack{i=m \\ i \le j}}^{k} \intl_{ \R^{n}}   P^{ 0}_{ x_0, 2^{ j+2}}(\partial _\beta   u_\alpha ) 
P^{ N}_{ x_0, 2^{ j+2}}( v_\beta )(x -y) \partial _l\partial _\alpha   (\Gamma(y) \psi _i (y) )dy \in \mathcal{P}_{ N},
\]
we  see that 
\begin{align*}
&g_{ m}^k(x)  - Q (x) 
\\
&=\sum_{\substack{i=m \\ i \le j}}^{k} \intl_{ \R^{n}} 
\Big( \partial _\beta  u_\alpha   v_\beta   -
P^{ 0}_{ x_0, 2^{ j+2}}(\partial _\beta   u_\alpha ) 
P^{ N}_{ x_0, 2^{ j+2}}( v_\beta )\Big) (x-y) \partial _l \partial _\alpha  (\Gamma(y) \psi _i(y)) dy
\\
& =
\sum_{\substack{i=m \\ i \le j}}^{k} \intl_{ \R^{n}} 
\partial _\beta   u_\alpha 
( v_\beta   -P^{ N}_{ x_0, 2^{ j+2}}(  v_\beta ))(x-y)
 \partial _l \partial _\alpha   (\Gamma(y) \psi _i(y)) dy
\\ 
&+\sum_{\substack{i=m \\ i \le j}}^{k} \intl_{ \R^{n}} 
( \partial _\beta u_\alpha   -  P^{0}_{ x_0, 2^{ j+2}}(\partial _\beta    u_\alpha )) 
P^{ N}_{ x_0, 2^{ j+2}}(  v_\beta )(x -y)\partial _l \partial _\alpha   (\Gamma(y) \psi _i(y)) dy
\\
&= J_1(x)+ J_2(x).
\end{align*}
(i) {\it Estimation of $ J_1$}: Arguing as in the proof of Lemma\,\ref{lem6.1}, using Calder\'on-Zygmund inequality and 
H\"older's inequality,  
we find for all $ 1< r \le p$
\begin{align*}
\| J_1\|_{ L^r(B(x_0, 2^j))}
&\le  c\Big\| \nabla    u \cdot  
( v   -P^{ N}_{ x_0, 2^{ j+2}}(  v ))\Big\|_{ L^r(B(x_0, 2^{ j+2}))}  
\\
&\le c \| \nabla u\|_{ L^{ \frac{pr}{p-r}}(B(x_0, 2^{ j+2}))} \|  v - P^{ N}_{ x_0, 2^{ j+2}} (v) \|_{ L^p(B(x_0, 2^{ j+2}))}. 
\end{align*}

(ii) {\it Estimation of $ J_2$}:  Applying integration by parts and recalling that $ \nabla \cdot u=0$, we infer 
\begin{align*}
&J_2(x)
=-\sum_{\substack{i=m \\ i \le j}}^{k} \intl_{ \R^{n}} 
\Big(  u_\alpha   -   P^{ N}_{ x_0, 2^{ j+2}}(  u_\alpha) 
P^{ 0}_{ x_0, 2^{ j+2}}( \nabla \cdot  v )\Big)(x-y)\partial _l \partial _\alpha   (\Gamma(y) \psi _i(y)) dy
\\
& \qquad \qquad +\sum_{\substack{i=m \\ i \le j}}^{k} \intl_{ \R^{n}} 
(    u_\alpha   -   P^{ N}_{ x_0, 2^{ j+2}}(  u_\alpha) )
 P^{ 0}_{ x_0, 2^{ j+2}}(\partial _\alpha   v_\beta )(x-y) \partial _l \partial _{ \beta }   (\Gamma(y) \psi _i(y)) dy. 
\end{align*} 
Once more applying Calder\'on-Zygmund's inequality using Poincar\'e's inequality, arguing as above, we obtain   
for $ 1 < r \le p$
\begin{align*}
\| J_2\|_{ L^r(B(x_0, 2^j))}
\le c  \|  u - P^{ N}_{ x_0, 2^{ j+2}} (u) \|_{ L^p(B(x_0, 2^{ j+2}))} \| \nabla v\|_{L^{ \frac{pr}{p-r}}(B(x_0, 2^{ j+2}))}. 
\end{align*}

Employing the two estimates for $ J_1$ and $ J_2$, we get  
\begin{align}
&\|g _{ m }^k- Q \|_{ L^r(B(x_0, 2^j))} 
\cr
&\le c \| \nabla u\|_{ L^{ \frac{pr}{p-r}}(B(x_0, 2^{ j+2}))}\|  v - P^{ N}_{ x_0, 2^{ j+2}} (v) \|_{ L^p(B(x_0, 2^{ j+2}))} 
\cr
&\qquad + c 
\|  u - P^{ N}_{ x_0, 2^{ j+2}} (u) \|_{ L^p(B(x_0, 2^{ j+2}))} \| \nabla v\|_{ L^{ \frac{pr}{p-r}}(B(x_0, 2^{ j+2}))}. 
\label{6.36}
\end{align}
 By the aid of \eqref{6.36} we deduce the following estimate for the oscillation of $ g_m^k$.
\begin{align}
& \osc_{r, N}(g _{ m}^k, x_0, 2^j)
\cr
&\le c \| \nabla u\|_{L^{ \frac{pr}{p-r}}(B(x_0, 2^{ j+2}))}  \osc_{ p, N}(  v, x_0, 2^{ j+2})) 
 + c  \osc_{ p, N}( u, x_0, 2^{ j+2})  \| \nabla v\|_{ L^{ \frac{pr}{p-r}}(B(x_0, 2^{ j+2}))}. 
\label{6.36a}
\end{align}

Next, we estimate $  G_{ m}^k$. Let $ \sigma  \in \N_0^n$ be any multi index with $ | \sigma |= N+1$. Clearly, 
\begin{align*}
& D^\sigma   G_{ m}^k(x) 
\\
&=  \sum_{\substack{i=m \\ i > j}}^{k} \intl_{ \R^{n}} 
\Big(\partial _\beta  u_\alpha   v_\beta  - P^{ 0}_{ x_0, 2^{ i+2}}(\partial _\beta  u_\alpha )P^{ N}_{ x_0, 2^{ i+2}}( v_\beta  )\Big)(x-y)D^\sigma \partial _l\partial_ \alpha  (\Gamma(y)  \psi _i (y)) dy
\\
&= \sum_{\substack{i=m \\ i > j}}^{k} \intl_{ \R^{n}} 
\partial _\beta  u_\alpha (x-y)  ( v_\beta  (x-y) -P^{ N}_{ x_0, 2^{ j+2}}(  v_\beta )(x -y))
D^\sigma  \partial _l\partial_ \beta (\Gamma(y)  \psi _i (y)) dy
\\
& \quad +   \sum_{\substack{i=m \\ i > j}}^{k} \intl_{ \R^{n}} ( \partial _\beta  u_\alpha (x-y)   -P^{ 0}_{ x_0, 2^{ j+2}}(\partial _\beta  u_\alpha )) 
P^{ N}_{ x_0, 2^{ j+2}}(  v_\beta )(x -y)
D^\sigma  \partial _l\partial_ \alpha  (\Gamma(y)  \psi _i (y)) dy. 
\end{align*}
Let  $ j \le  i \le k$. Noting that $ B(x, 2^{ i+1}) \subset  B(x_0, 2^{ i+1}+ 2^{ j}) \subset B(x_0, 2^{ i+2})$, and employing 
Jensen's inequality, we infer  
\begin{align*}
&\bigg|\intl_{ \R^{n}} 
\partial _\beta  u_\alpha (x-y)  ( v_\beta  (x-y) -P^{ N}_{ x_0, 2^{ j+2}}(  v_\beta )(x -y))
D^\sigma  \partial _l\partial_ \beta (\Gamma(y)  \psi _i (y)) dy\bigg|
\\
&\quad \le c\intl_{ B(2^{ i+1})}  |\nabla u(x-y)|\,| v(x-y) -P^{ N}_{ x_0, 2^{ i+2}}(v)(x-y)|  2^{ -i (n +N+1)} dy
\\
&\quad = c\intl_{ B(x, 2^{ i+1})}  |\nabla u(y)|\,| v(y) -P^{N}_{ x_0, 2^{ i+2}}(v)(y)| 2^{ -i (n +N+1)} dy
\\
&\quad \le  c \| \nabla u\|_{ L^{ p'}(B(x_0, 2^{ i+2}))}2^{ -i (n +N+1)}\bigg(\intl_{ B(x_0, 2^{ i+2})}\,  | v-P^{ N}_{ x_0, 2^{ i+2}}(v)|^p  dy\bigg)^{ \frac{1}{p}}
\\
& \quad \le  c 2^{ - i\frac{n}{p'}}\| \nabla u\|_{ L^{ p'}(B(x_0, 2^{ i+2}))}2^{  -i(N+1)} \osc_{  p, N}( v; x_0, 2^{ i+2}). 
\end{align*}
Similarly, using integration by parts along with \eqref{6.31}, we get  
\begin{align*}
&\bigg|\intl_{ \R^{n}} ( \partial _\beta  u_\alpha (x-y)   -P^{ 0}_{ x_0, 2^{ j+2}}(\partial _\beta  u_\alpha )) 
P^{ N}_{ x_0, 2^{ j+2}}(  v_\beta )(x -y)
D^\sigma  \partial _k\partial_ \alpha  (\Gamma(y)  \psi _i (y)) dy \bigg|
\\
& \quad \le  c 2^{ - i\frac{n}{p'}}\| \nabla v\|_{ L^{ p'}(B(x_0, 2^{ i+2}))}2^{  -i(N+1)} \osc_{  p, N}( u; x_0, 2^{ i+2}). 
\end{align*}
This yields
\begin{align}
&\|  D^{ N} G_{ m}^k \|_{ L^\infty(B(x_0, 2^j))}   
\cr
 &\le c  \sum_{i=j}^{\infty}2^{ - i\frac{n}{p'}}\| \nabla u\|_{ L^{ p'}(B(x_0, 2^{ i+2}))}2^{  -i(N+1)} \osc_{  p, N}( v; x_0, 2^{ i+2})
 \cr
 &\qquad + c\sum_{i=j}^{\infty}2^{ - i\frac{n}{p'}}\| \nabla v\|_{ L^{ p'}(B(x_0, 2^{ i+2}))}2^{  -i(N+1)} \osc_{  p, N}( u; x_0, 2^{ i+2}). 
\label{6.40a}
\end{align}
Noting that for all $ i \ge j$
\begin{align*}
&2^{ - i\frac{n}{p'}}\| \nabla u\|_{ L^{ p'}(B(x_0, 2^{ i+2}))} \le c |\nabla u|_{ BMO}  + 
c\sup_{ i \ge j} |P^0_{ x_0, 2^i}(\nabla u)|
\\
&2^{ - i\frac{n}{p'}}\| \nabla v\|_{ L^{ p'}(B(x_0, 2^{ i+2}))} \le c |\nabla v|_{ BMO}  + 
c\sup_{ i \ge j} |P^0_{ x_0, 2^i}(\nabla v)|  
\end{align*}
we infer from   \eqref{6.40a} 
\begin{align}
&2^{ j(N+1)}\|  D^{ N} G_{ m}^k \|_{ L^\infty(B(x_0, 2^j))}  
\cr
& \le c\Big( |\nabla u|_{ BMO}  + 
\sup_{ i \ge j} |P^0_{ x_0, 2^i}(\nabla u)|\Big) S_{ N+1,1}(\osc_{p,  N}(v; x_0))_j 
\cr
&\qquad + c
\Big(|\nabla v|_{ BMO}  + 
\sup_{ i \ge j} |P^0_{ x_0, 2^i}(\nabla v)| \Big)S_{ N+1,1}(\osc_{ p, N}(u; x_0))_j. 
\label{6.40}
\end{align}

With the  help of Poincar\'e's inequality \eqref{5.8} along with \eqref{6.40} we find 
\begin{align}
&\osc_{r, N}( G_{m}^k; x_0))_j
\cr
& \le 2^{ - j \frac{n}{r}}\|  G_{m}^k- P^{N}_{ x_0, 2^j}(G_{m}^k) \|_{ L^r(B(x_0, 2^j))}  
\cr
&\le c 2^{ j(N+1)  } \| D^{ N+1} G_m^k \|_{ L^\infty(B(x_0, 2^j))}
\cr
&\le c\Big( |\nabla u|_{ BMO}  + 
\sup_{ i \ge j} |P^0_{ x_0, 2^i}(\nabla u)|\Big) S_{ N+1,1}(\osc_{p,  N}(v; x_0))_j 
\cr
&\qquad + c
\Big(|\nabla v|_{ BMO}  + 
\sup_{ i \ge j} |P^0_{ x_0, 2^i}(\nabla v)| \Big)S_{ N+1,1}(\osc_{ p, N}(u; x_0))_j.
\label{6.41}
\end{align}
Combining this inequality with \eqref{6.36a} and \eqref{6.41},  and noting that for all $ j\in \Z$
\begin{align*}
\osc_{ p,N}(u; x_0, 2^{ j+2}) &\le c S_{ N+1,1}(\osc_{ p,N}( u; x_0))_j,\quad  
\\
\osc_{ p,N}(v; x_0, 2^{ j+2}) &\le c S_{ N+1,1}(\osc_{ p,N}( v; x_0))_j,
\end{align*}
we obtain in case $ r=p$
\begin{align}
 &\osc_{ p,N}(f^k_m; x_0) 
\cr
&\quad \le c\| \nabla u\|_{ \infty} S_{ N+1,1}(\osc_{ p,N}(v; x_0)) + c\| \nabla v\|_{ \infty} S_{ N+1,1}(\osc_{p, N}(u; x_0)). 
\label{6.41a}
\end{align}
and in case $ 1<r<p$
\begin{align}
&(\osc_{r, N}( f_{m}^k; x_0))_j
\cr
&\le c\Big( |\nabla u|_{ BMO}  + 
\sup_{ i \ge j} |P^0_{ x_0, 2^i}(\nabla u)|\Big) S_{ N+1,1}(\osc_{p,  N}(v; x_0))_j 
\cr
&\qquad + c
\Big(|\nabla v|_{ BMO}  + 
\sup_{ i \ge j} |P^0_{ x_0, 2^i}(\nabla v)| \Big)S_{ N+1,1}(\osc_{ p, N}(u; x_0))_j.
\label{6.41b}
\end{align}

We are now in a position to apply  Lemma\,\ref{lem10.1} with $ \alpha = N+1, \beta = s, p=1$. Performing $ S_{ s, q}$ 
to both sides of \eqref{6.41a}, we get  
 \begin{equation}
 |  f_m^k|_{ \mathscr{L}^{ s}_{ q (p,N)}} \le c \Big(\| \nabla u\|_{ \infty} 
  | \nabla v|_{ \mathscr{L}^{ s}_{ q (p,N)}} + \| \nabla v\|_{ \infty} 
  | \nabla u|_{ \mathscr{L}^{ s}_{ q (p,N)}}\Big). 
 \label{6.42}
 \end{equation}
 In particular, for $ N=1, s=1$ and $ q=1$ it follows that  
 \begin{equation}
 |  f^k_m|_{ \mathscr{L}^{ 1}_{ 1 (p,1)}} \le c \Big(\| \nabla u\|_{ \infty} 
  | \nabla v|_{ \mathscr{L}^{ 1}_{ 1 (p, 1)}} + \| \nabla v\|_{ \infty} 
  | \nabla u|_{ \mathscr{L}^{ 1}_{ 1 (p, 1)}}\Big). 
 \label{6.43}
 \end{equation}
  
  \vspace{0.3cm}
 Set $ g_m^k= f_m^k - P^1_{ 0, 1}(f_m^k)$.  Clearly, $ \{ g_m^k\}$ is bounded in $  \mathscr{L}^{ s}_{ q (p,N)}\cap \mathscr{L}^{1}_{ 1 (p,1)}(\R^{n}) \hookrightarrow C^{ 0, 1}(\R^{n})$.
  Arguing as in the proof of Theorem\,\ref{thm6.3}, we get a  very weak solution 
  $ g\in  \mathscr{L}^{ s}_{ q (p,N)}\cap \mathscr{L}^{1}_{ 1 (p,1)}(\R^{n})$  with  $ \nabla  g = L^\infty(\R^{n})$
 to  
  \[
 -\Delta g =  \partial _k\nabla \cdot ((v\cdot \nabla) u)\quad    \text{ in}\quad   \R^{n}. 
\]
We make the ansatz: $ f= g + Ax +b$. Clearly, for all $ A\in \R^{n\times n}$ and $ b\in \R^{n}$, 
$  f\in \mathscr{L}^{s}_{ q (p,N)}\cap \mathscr{L}^{1}_{ 1 (p,1)}(\R^{n})$  is a very weak solution 
to \eqref{6.32}.  The condition $ f(0)=0$ implies $ b=-g(0)$, while the first condition in \eqref{6.33} implies 
\[
P^0_{ \infty}(\nabla u:(\nabla v)^{ \top} )x_l = \dot{P}_{ \infty} ^1(f)= \dot{P}_{ \infty} ^1(g + Ax +b)
= \dot{P}_{ \infty} ^1(g) + Ax. 
\] 
Setting $ Ax= P^0_{ \infty}(\nabla u:(\nabla v)^{ \top} )x_l- \dot{P}_{ \infty} ^1(g)$, the function $ f$ fulfills \eqref{6.33}. Arguing as in the proof of Theorem\,\ref{thm6.3}, the estimate  \eqref{6.34} follows from  \eqref{6.42} and  \eqref{6.43}. Furthermore, the estimate  \eqref{6.41c} follows from  \eqref{6.41b} after passing $ k \rightarrow +\infty$ and $ m \rightarrow -\infty$.

\vspace{0.3cm}
{\it Uniqueness}. Let $ \overline{f} \in \mathscr{L}^{ s}_{ q (p,N)}\cap \mathscr{L}^{ 1}_{ 1 (p,1)}(\R^{n})$ be a second solution 
which satisfies \eqref{6.33}. Clearly, $f- \overline{f}  $ is harmonic. Arguing as in the proof of Theorem\,\ref{thm1.2}, we conclude 
$ f- \overline{f} \in \mathcal{\dot{P}}_1 $, and by \eqref{6.33} $ f= \overline{f} $.

\hfill \Beweisende 
 
\begin{rem}
\label{rem6.5a}
Noting that $ {\D  \osc_{ p,N}(f^k_m; x_0, 2^j)  \rightarrow  \osc_{p, N}( f; x_0, 2^j)} $ as $ k \rightarrow +\infty$ and 
$ m \rightarrow -\infty$, we infer from \eqref{6.41a} the estimate  
\begin{align}
\osc_{ p,N}( f; x_0) 
\le c\| \nabla u\|_{ \infty} S_{ N+1, 1}(\osc_{p,N}(v; x_0)) + c\| \nabla v\|_{ \infty} S_{ N+1,1}
(\osc_{ p,N}(u; x_0)). \label{6.44a}
\end{align}
\end{rem} 
 
 \begin{rem}
 \label{rem6.5} By $ \mathscr{L}^{ 1}_{ 1 (p, 1), \sigma } (\R^{n})$ we define the space of all $ u\in
  \mathscr{L}^{ 1}_{ 1 (p,1)} (\R^{n})$ such that $ \nabla \cdot u=0$ in $ \R^{n}$.  Then by Theorem\,\ref{thm6.4} we are able to construct 
 the pressure $ \pi\in  \mathscr{L}^{1, 1}_{1 (p, 1), \sigma } (\R^{n}) $, by the relation $ \nabla \pi = f$, where 
 $ f\in  \mathscr{L}^{ 1}_{ 1 (p,1)} (\R^{n})$ is the unique very weak solution to \eqref{6.32}, \eqref{6.33}.  In fact, from \eqref{6.32} it follows that $ \nabla \times f$ is harmonic and bounded. 
Thus, by Liouville's theorem for harmonic functions we see that $\nabla \times f $ is constant. 
On the other hand, by means of \eqref{6.33} we find 
\[
\nabla \times f= P ^0_\infty  (\nabla \times f) =
\nabla \times\dot{P}^1_\infty  ( f) = \nabla \times (P ^0_\infty(\nabla u : (\nabla v)^{ \top})x) =0.  
\]   
Thus, $ f= \nabla \pi $ for a unique  $ \pi \in \mathscr{L}^{1, 1}_{1 (p, 1)} (\R^{n})$ fulfilling  $ [\pi]^0 _{0,1 }=0$.
Setting $ \nabla \Pi (u, v)= \nabla \pi $, defines a    linear mapping 
 $ \nabla \Pi: \mathscr{L}^{ 1}_{ 1 (p, 1), \sigma } (\R^{n})\times 
 \mathscr{L}^{ 1}_{ 1 (p,1)} (\R^{n}) \rightarrow  \mathscr{\dot{L} }^{1}_{ 1 (p, 1)} (\R^{n})$. 
Accordingly, it holds 
 \begin{align}
&- \Delta \pi  = \nabla u: (\nabla v)^{ \top}\quad  \text{ in }\quad  \R^{n}, 
 \label{6.48}
 \\
 & \dot{P}^1_\infty(\nabla \pi )  = - \frac{1}{n } P^0_\infty( \nabla  u: (\nabla  v)^{ \top}) x. 
 \label{6.49}
 \end{align}
In addition, from \eqref{6.34} we deduce that $ \nabla \Pi $ is bounded. More precisely,   
\begin{equation}
| \nabla \Pi (u, v)|_{ \mathscr{L}^{1}_{1 (p, 1)} }
\le c \Big(\| \nabla u\|_{ \infty} |v |_{ \mathscr{L}^{1}_{1 (p, 1)}}
+ \| \nabla v\|_{ \infty} |u |_{ \mathscr{L}^{1}_{1 (p, 1)}} \Big), 
\label{6.50}
\end{equation}
and in view of \eqref{6.33}
\begin{equation}
| \dot{P} ^1_{\infty}(\nabla \Pi (u, v))| \le c|  P ^0_{\infty}(\nabla u)|\, |  P ^0_{\infty}(\nabla v)|.   
\label{6.51}
\end{equation} 

In case $ u, v\in  \mathscr{L}^{ s}_{ q (p, N)}(\R^{n})$ we get from \eqref{6.34}$ _2$ 
\begin{align}
| \nabla \Pi (u, v)|_{  \mathscr{L}^{ s}_{ q (p, N)}}
  \le  c\Big\{\| \nabla u\|_{ \infty}
| v|_{ \mathscr{L}^{ s}_{ q (p, N)}}
+ \| \nabla v\|_{ \infty}
| u|_{ \mathscr{L}^{ s}_{ q (p, N)}}\Big\}.
\label{6.52}
\end{align}
 In addition,  given $ 1 < r < p$, in view of \eqref{6.41a} for all $ j\in \Z$ it holds 
\begin{align}
&(\osc_{r, N}( \nabla \Pi (u, v); x_0))_j
\cr
&\le c\Big( |\nabla u|_{ BMO}  + 
\sup_{ i \ge j} |P^0_{ x_0, 2^i}(\nabla u)|\Big) S_{ N+1,1}(\osc_{p,  N}(v; x_0))_j 
\cr
&\qquad + c
\Big(|\nabla v|_{ BMO}  + 
\sup_{ i \ge j} |P^0_{ x_0, 2^i}(\nabla v)| \Big)S_{ N+1,1}(\osc_{ p, N}(u; x_0))_j.
\label{6.53}
\end{align}

 \end{rem}

\vspace{0.3cm}
We also get the following pressure estimate for the case of sublinear growth less then $ \frac{1}{2}$. The following theorem 
will be used in the proof of Theorem\,\ref{thm3}. 

\begin{thm}
\label{thm6.8}  
Let $1 \le q \le +\infty$ and $ s\in [0, \frac{1}{2})$. Let $ v\in \mathscr{L}^{ s}_{ q (p, 0), \sigma }(\R^{n})\cap C^{ 0,1}(\R^{n})$. 
 Let $ 1 < r< (1-s)p$ with $ r \ge \frac{p}{2}$.   There exists a unique very weak solution 
$ \pi \in \mathscr{L}^{ 2- (1-s) \frac{p}{r}}_{ q (r,0), \sigma }(\R^{n})$ to the equation 
\begin{equation}
-\Delta \pi  = \nabla \cdot \nabla \cdot (v \otimes v)\quad  \text{ in}\quad  \R^{n}
\label{6.57a}
\end{equation} 
with $ P^0_{ 0, 1}(\pi )=0$. 
Furthermore, for all $ x_0\in \R^{n}$ and $ j\in \Z$ it holds 
\begin{equation}
\osc_{p, 0}( \pi ; x_0, 2^j) \le c \| \nabla v\|_{ \infty}^{ 2- \frac{p}{r}}  S_{1,1}
\Big(\{2^{ (2- \frac{p}{r}) i}\osc_{p, 0}(v; x_0, 2^i)^{ \frac{p}{r}}\} \Big)_j.
\label{6.57}
\end{equation} 

\end{thm}
 
 {\bf Proof}: Let $ 1< r< (1-s)p$ with $ r \ge \frac{p}{2}$.  
    Let $ m,k\in \Z$ with $ m< k$. Define, 
 \[
f_{ m}^k(x) = \sum_{i=m}^{k} \intl_{ \R^{n}} 
 v_\alpha(x-y)  v_\beta(x-y)
 \partial _\alpha  \partial _\beta    (\Gamma(y) \psi _i(y)) dy,\quad  x\in \R^{n}. 
\]
 Let $ j\in \{m, \ldots , k-1\}$ be fixed. Our aim is to estimate 
 $ \osc_{ r, 0}(f_m^k; x_0, 2^j)$. We decompose $ f_{ m}^k$ into the sum $ g_m^k + G_m^k$, where 
\begin{align*}
g_{ m}^k(x) &= \sum_{i=m}^{j} \intl_{ \R^{n}} 
 v_\alpha(x-y)  v_\beta(x-y)
 \partial _\alpha  \partial _\beta    (\Gamma(y) \psi _i(y)) dy,
\\
G_{ m}^k(x) &= \sum_{i=j+1}^{k} \intl_{ \R^{n}} 
 v_\alpha(x-y)  v_\beta(x-y)
 \partial _\alpha  \partial _\beta    (\Gamma(y) \psi _i(y)) dy,\quad  x\in \R^{n}.
\end{align*}
 Recalling that $ \nabla \cdot v=0$, we obtain 
   \begin{align*}
g_{ m}^k(x) & =
\sum_{i=m}^{j} \intl_{ \R^{n}} 
(v_\alpha(x-y) - P^0_{ x_0, 2^{ j+2}}) (v_\beta(x-y)- P^0_{ x_0, 2^{ j+2}})
 \partial _\alpha  \partial _\beta    (\Gamma(y) \psi _i(y)) dy.
\end{align*}
Arguing as in the proof of Theorem\,\ref{lem6.1}, by the aid of Calder\'on-Zygmund's estimate along with H\"older's inequality we find 
\begin{align*}
\| g_m^k\|_{ L^r(B(x_0, 2^j))}
&\le  c\| v- P^0_{ x_0, 2^{ j+2}}(v)\|^{ 2}_{ L^{ 2r}(B(x_0, 2^{ j+2}))}
\\
&\le c 2^{ j (2- \frac{p}{r})}\| \nabla v\|_{ \infty}^{ 2- \frac{p}{r}}\| v- P^0_{ x_0, 2^{ j+2}}(v)\|^{ \frac{p}{r}}_{ L^{ p}(B(x_0, 2^{ j+2}))}. 
\end{align*}
This shows that 
\begin{equation}
\osc_{ r, 0} (g_{ m}^k; x_0, 2^j)  
\le   c 2^{ j (2- \frac{p}{r})}\| \nabla v\|_{ \infty}^{ 2- \frac{p}{r}}\osc_{p, 0} (v; x_0, 2^{ j+2})^{ \frac{2}{p}}. 
\label{6.58}
\end{equation}
Arguing as in Lemma\,\ref{lem6.1}, we obtain 
 \begin{align*}
&2^{ j}\|  \nabla  G_{ m}^k \|_{ L^\infty(B(x_0, 2^j))}   \le c\| \nabla v\|_{ \infty}^{2- \frac{p}{r}}
S_{  1,1}(\{2^{ (2- \frac{p}{r})i}\osc_{p, 0}(v; x_0, 2^i)^{ \frac{p}{r}}\})_j. 
\end{align*}
Using Poincar\'e's inequality, this yields
 \begin{align}
 \osc_{ r, 0} (G_{ m}^k; x_0, 2^j)  
 \le  c \| \nabla v\|_{ \infty}^{ 2- \frac{p}{r}}
S_{  1,1}(\{2^{ (2- \frac{p}{r})i}\osc_{p, 0}(v; x_0, 2^i)^{ \frac{p}{r}}\})_j. 
\label{6.59}
\end{align}

 Combining \eqref{6.58} and \eqref{6.59}, and letting $ m \rightarrow -\infty$, we arrive at 
\begin{align}
 \osc_{ r, 0} (f^k; x_0, 2^j)  
&\le  c \| \nabla v\|_{ \infty}^{ 2- \frac{p}{r}}
S_{  1,1}(\{2^{ (2- \frac{p}{r})i}\osc_{p, 0}(v; x_0, 2^i)^{ \frac{p}{r}}\})_j 
\cr
&\le c\| \nabla v\|_{ \infty}^{ 2- \frac{p}{r}} | v|_{ \mathscr{L}^s_{ q(p, 0)}}^{ \frac{p}{r}-1}
S_{1,1}(\{2^{i (2- (1-s)\frac{p}{r})-is}\osc_{p, 0}(v; x_0, 2^i)\})_j. 
\label{6.60}
\end{align}
Recalling that $ r < p(1-s)$, we have $ 2-  (1-s)\frac{p}{r}  <1$.  Applying $ S_{ 2- (1-s) \frac{p}{r}, q}$ to both sides of \eqref{6.60},  thanks to Lemma\,\ref{lem10.1} with $ \alpha =1$ and $ \beta = 2- (1-s) \frac{p}{r} $ we get 
\begin{align*}
S_{ 2- (1-s) \frac{p}{r}, q}(\osc_{ p, 0} (f^k; x_0)  )_j &\le c\| \nabla v\|_{ \infty}^{ 2- \frac{p}{r}} | v|_{ \mathscr{L}^s_{ q(p, 0)}}^{ \frac{p}{r}-1}
S_{2-  (1-s)\frac{p}{r} ,q}(\{2^{ i(2- (1-s)\frac{p}{r}) -is}\osc_{p, 0}(v; x_0, 2^i)\})_j
\\
&= c 2^{j (2-(1-s) \frac{p}{r}) - js}\| \nabla v\|_{ \infty}^{ 2- \frac{p}{r}} | v|_{ \mathscr{L}^s_{ q(p, 0)}}^{ \frac{p}{r}-1}S_{s,q}(\osc_{p, 0}(v; x_0))_j
\\
& \le c2^{j (2- (1-s)\frac{p}{r})}\| \nabla v\|_{ \infty}^{ 2- \frac{p}{r}} 
| v|_{ \mathscr{L}^s_{ q(p, 0)}}^{ \frac{p}{r}}. 
\end{align*}
Multiplying both sides by $ 2^{j (2- (1-s)\frac{p}{r})}$, and taking the supremum over all $ x_0 \in \R^{n}$, we arrive at  
\begin{align}
 | f^k|_{  \mathscr{L}^{ 2- (1-s) \frac{p}{r}}_{ q(r, 0)}}
\le   c\| \nabla v\|_{ \infty}^{ 2- \frac{p}{r}} 
| v|_{ \mathscr{L}^s_{ q(p, 0)}}^{ \frac{p}{r}}. 
 \label{6.61}
\end{align}
Set $ g^k = f^k - P^0_{ 0,1}(f^k)$. In view of \eqref{6.61} $ \{ g^k\}$ is bounded in 
$ \mathscr{L}^{ 2- (1-s) \frac{p}{r}}_{ q(r, 0)}(\R^{n})$.  
Thus, by the compact embedding $ \mathscr{L}^{ 2- (1-s) \frac{p}{r}}_{ q(r, 0)}(\R^{n}) \hookrightarrow  L^r_{ loc}(\R^{n})$, eventually passing 
to a sub sequence, we get $ \pi \in \mathscr{L}^{ 2- (1-s) \frac{p}{r}}_{ q(r, 0)}(\R^{n})$ such that 
\begin{equation}
 g^k \rightarrow \pi   \quad  \text{{\it in}}\quad  L^r_{ loc}(\R^{n})\quad  \text{{\it as}}\quad  k \rightarrow +\infty. 
\label{6.62}
\end{equation}
Arguing as in the proof of Theorem\,\ref{thm6.4}, it can be checked that $ \Delta \pi  = - \nabla \cdot \nabla \cdot (v \otimes v)$ in the 
sense of distributions.  Using \eqref{6.62}, we immediately get \eqref{6.57} from the first inequality in \eqref{6.60}  and $ P^0_{ 0,1}(\pi)=0$ from 
$ P^0_{ 0,1}(g^k)=0$.  

\vspace{0.2cm}
{\it Uniqueness}. Assume there is a second very weak  solution $ \overline{\pi }\in \mathscr{L}^{ 2- (1-s) \frac{p}{r}}_{ q(r, 0)}(\R^{n})$ to \eqref{6.57a}. Then by Weyl's Lemma $ \pi - \overline{\pi } $ is harmonic. Thus, by Liouvill's theorem of harmonic functions it follows that 
$ \pi - \overline{\pi }$ is constant. Taking into account the condition $ P^0_{ 0,1}(\pi - \overline{\pi })=0$ we obtain 
$ \pi = \overline{\pi } $.  \hfill \Beweisende     

\vspace{0.3cm}
 A careful inspection of the proof of Theorem\,\ref{thm6.8} shows that we may remove the condition $ v\in C^{ 0,1}(\R^{n} )$ in case $ 2< p< +\infty$ and $ r= \frac{p}{2}$. Thus, we have the following 

\begin{cor}
 \label{cor6.9}
 Let $1 \le q \le +\infty, 2 < p <+\infty$ and $ s\in [0, \frac{1}{2})$. Let $ v\in \mathscr{L}^{ s}_{ q (p, 0), \sigma }(\R^{n})$.  There exists a unique very weak solution 
$ \pi \in \mathscr{L}^{2s}_{q ( \frac{p}{2},0), \sigma }(\R^{n})$ to  \eqref{6.57a} 
with $ P^0_{ 0, 1}(\pi )=0$. 
Furthermore, it holds 
\begin{equation}
|\pi |_{ \mathscr{L}^{2s}_{q ( \frac{p}{2},0) }}\le c |v|_{ \mathscr{L}^{s}_{q ( p,0) }}^2. 
\label{6.65}
\end{equation} 
In particular, if $ v\in BMO$, then  $ \pi \in BMO$ and it holds 
\begin{equation}
 |\pi |_{ BMO} \le c |v|_{ BMO}^2. 
\label{6.66}
 \end{equation} 

\end{cor}

\section{Proof of Theorem\,\ref{thm1}}
\label{sec:-2}
\setcounter{secnum}{\value{section} \setcounter{equation}{0}
\renewcommand{\theequation}{\mbox{\arabic{secnum}.\arabic{equation}}}}

Before turning to the  proof of Theorem\,\ref{thm1} we state  the following local  energy inequality for the transport equation,  which is proved in 
\cite{trans}. 

\begin{lem}
\label{lem2.1}  Given $v \in  L^1(0, T; C^{ 0,1}(\R^{n}))$ 
and $ g\in L^1(0,T; L^p_{ loc}(\R^{n} ))$,  let $ f\in C([0, T]; W^{1,\,p} _{ loc}(\R^{n}))$ be a weak  solution to the   transport equation 
\begin{equation}
\partial _t f + v\cdot \nabla f = g\quad  \text{ in}\quad  Q_T. 
\label{2.2}
\end{equation}
Let $ N \in \N_0$.  Then 
the following inequality  holds  for all $ t\in [0,T]$
 \begin{align*}
& \osc_{ p, \max\{2N-1, N\}} \Big(f(t); x_0,  \frac{r}{2}\Big)
\le   c \osc_{ p, N} (f(0); x_0,  r) +c r^{ -1} \intl_{0}^{t} \|   v(\tau )\|_{ L^\infty(B(x_0, r))}  
\osc_{ p, N} (f(\tau ); x_0,  2r)d\tau  
 \cr
&\qquad  +c \intl_{0}^{t}   \|  \nabla \cdot  v(\tau )\|_{ L^\infty(B(x_0, r))} \osc_{ p, N} (f(\tau ); x_0,  2r)d\tau  
\end{align*}
\begin{align}
& \qquad + \delta_{ N0} c  \intl_{0}^{t}  \osc_{ p, N} (v(\tau ); x_0,  r) 
\| \nabla P^{ N}_{ x_0, r}(f(\tau ))\|_{ L^\infty(B(x_0, r))} d\tau 
\cr
&\qquad +  c\intl_{0}^{t} \osc_{ p, N} (g(\tau ); x_0,  r)  d\tau,
  \label{2.4}
 \end{align}       
where $ \delta _{ N0}=0$ if $ N=0$ and $ 1$ otherwise. 
\end{lem}

\begin{rem}
 \label{rem2.2} Given
 $v\in  L^1(0, T; C^{ 0,1}(\R^{n}; \R^{n} ))$,  and 
$ \pi  \in L^1(0, T; W^{1,\,2}_{ loc}(\R^{n}; \R^{n}  ))$, 
let $ f\in L^\infty(0, T; C^{ 0,1}(\R^{n}; \R^{n} ))$ with $ \nabla \cdot f=0$ be a weak solution to the system 
 \begin{equation}
\partial _t f + (v\cdot \nabla) f = -\nabla \pi \quad  \text{ in}\quad  Q_T. 
\label{2.2s}
\end{equation}

Then,   \eqref{2.4} can be replaced by 
 \begin{align}
& \osc_{ 2, 1} \Big(f(t); x_0,  \frac{r}{2}\Big)
\le   c \osc_{2, 1} (f(0); x_0,  r) +c r^{ -1} \intl_{0}^{t} \|   v(\tau )\|_{ L^\infty(B(x_0, r))}  
\osc_{ 2, 1} (f(\tau ); x_0,  2r)d\tau  
 \cr
&\quad  +c \intl_{0}^{t}   \|  \nabla \cdot  v(\tau )\|_{ L^\infty(B(x_0, r))} \osc_{2, 1} (f(\tau ); x_0,  2r)d\tau  
\cr
& \quad +  c  \intl_{0}^{t}  \osc_{ 2, 1} (v(\tau ); x_0,  r) 
 |\nabla P^{ 1}_{ x_0, r}(f(\tau )| d\tau 
\cr
&\quad +  c\intl_{0}^{t} \osc_{ \frac{2n}{n+2}, 1} (\nabla \pi (\tau ); x_0,  r)  d\tau.
  \label{2.4a}
 \end{align}

\end{rem}

\vspace{0.3cm}
{\bf Proof of Theorem\,\ref{thm1}.}  The proof of  Theorem\,\ref{thm1} is based on a fixed point argument using 
Banach's fixed point theorem.  Let $ v_0\in \mathscr{\dot{L} }^{1}_{1 (p, 1)}(\R^{n})$ be arbitrarily chosen.  Let $ T_0 = \frac{1}{c \| v_0\|_{ \mathscr{\dot{L} }^{1}_{1 (p, 1)}}}$, 
with a constant $ c>0$ which will be specified below. 
We construct an operator 
$ \mathcal{T}: L^\infty(0, T_0; \mathscr{\dot{L} }^{1}_{1 (p, 1)}(\R^{n})) \rightarrow L^\infty(0, T_0; \mathscr{\dot{L} }^{1}_{1 (p, 1)}(\R^{n}))$ as follows. Given 
$ u\in L^\infty(0, T_0; \mathscr{\dot{L} }^{1}_{1 (p, 1)}(\R^{n}))$ by $ \mathcal{T}(u):=v\in L^\infty(0, T_0; \mathscr{\dot{L} }^{1}_{1 (p, 1)}(\R^{n}))$   we denote the unique solution to the model 
problem 
\begin{equation}
\begin{cases}
\partial _t v + (\PP(u) \cdot \nabla) v = - \nabla \Pi (\PP(u),u)\quad  \text{ in}\quad  Q_{ T_0}
\\[0.3cm]
v= v_0\quad  \text{ on}\quad  \R^{n}\times \{ 0\}.
\end{cases}
\label{2.9}
\end{equation}
Here $ \nabla \Pi (\PP(u),u)\in L^\infty(0, T_0; \mathscr{\dot{L} }^{1}_{1 (p, 1)}(\R^{n}))$, and $ \pi(\tau ) := \Pi (\PP(u(\tau )),u(\tau ))\in 
\mathscr{L}^{ 1,1}_{ 1 (1)}(\R^{n})$ 
stands for the unique solution  to the Poisson equation 
\begin{equation}
\begin{cases}
- \Delta \pi(\tau )  = \nabla \PP(u(\tau )) : (\nabla u(\tau ))^{ \top}\quad  \text{ in}\quad \R^{n}
\\[0.3cm]
\nabla \pi (0)=0,\quad  \dot{P} ^1_{ \infty}(\nabla  \pi(\tau )  ) = - \frac{1}{n}P^0_{ \infty}(\nabla \PP(u(\tau ) : (\nabla (u(\tau ))^{ \top}) x
\end{cases}
\label{2.10}
\end{equation}
with $ ( \pi (\tau ) )_{ B(1)}=0$.  According to \eqref{6.50}, \eqref{6.51} (cf. Remark\,\ref{rem6.5}) the following  estimate holds true 
for all $ \tau \in [0,T_0]$
\begin{align}
\| \nabla \Pi (\PP(u(\tau )), u(\tau ))\|_{ \mathscr{\dot{L} }^{1}_{1 (p, 1)} } &\le c\| u(\tau )\|^2_{ \mathscr{\dot{L} }^{1}_{1 (p, 1)} }. 
\label{2.11a}
\end{align}
Furthermore, we wish to remark that  \cite[Theorem\,1.2]{trans}  ensures the existence and uniqueness of the solution 
$ v =  \mathcal{T}(u)$. 

\vspace{0.3cm}
Let $ x_0\in \R^{n}$.   Let $ 0<t \le  T_0$ be arbitrarily chosen, but fixed. By $ \xi \in C^{ 1,1}([0,T])$ 
we denote the unique solution to the ODE
\begin{equation}
\dot{\xi}  (\tau ) = u(x_0+ \xi (\tau ), \tau ),\quad  \tau \in [0, t], \quad  \xi (t)= 0.  
\label{2.12}
\end{equation}
We set
\begin{align*}
V(x,\tau ) &=  v  (x+ \xi (\tau ), \tau ),\quad 
\\
U(x,\tau ) &= u(x+ \xi (\tau ), \tau )-  \dot{\xi}  (\tau ), 
\\
 g(x,\tau  ) &= - \nabla \Pi (\PP(u(\tau )),u(\tau ))(x+ \xi (\tau ), \tau),\quad  (x,\tau )\in Q_t. 
\end{align*}
It is readily seen that $ V$ solves the transport equation 
\begin{equation}
\partial _t V + (U\cdot \nabla) V = g  \quad  \text{ in }\quad  Q_t.
\label{2.13}
\end{equation} 
Furthermore, it holds 
\begin{equation}
U(x_0, \tau )=0\quad \forall\,\tau \in (0, t],\quad  V(t)=v(t). 
\label{2.14}
\end{equation}

Let $ j\in \Z$ be arbitrarily chosen. Observing \eqref{2.4} with $ V$ ($ U$ respectively) in place of $ v$ ($ u$ respectively), 
$ r= 2^{ j+1}$, and $ N=1$, using \eqref{2.14},    we obtain 
 \begin{align}
& \osc_{p, 1} \Big(v(t); x_0,  2^{ j}\Big)
\le   c \osc_{p, 1} (V(0); x_0,  2^{ j+1})   
 \cr
&\quad  +c \intl_{0}^{t}   \|  \nabla  u(\tau )\|_{ \infty} \osc_{p, 1} (V(\tau ); x_0,  2^{ j+2})d\tau  
\cr
& \quad + c  \intl_{0}^{t}  \osc_{p, 1} (U(\tau ); x_0,  2^{ j+1}) \| \nabla v(\tau )\|_{\infty} d\tau 
 +  c\intl_{0}^{t} \osc_{p, 1} (g(\tau ); x_0,  2^{ j+1})  d\tau.
  \label{2.15}
 \end{align}       
Noting that all functions $ U, V$ and $ g$ belong to $ L^\infty(0, T; \mathscr{L}^1_{ 1 (p, 1)}(\R^{n}))$, 
 we may multiply both sides by $ 2^{ -j}$,   apply 
the sum over $ j\in \Z$ on  both sides and take the supremum over $ x_0\in \R^{n}$. This yields 
 \begin{align}
&| v(t)|_{ \mathscr{L}^1_{ 1 (p, 1)}}
\le   | V(0 )|_{ \mathscr{L}^1_{ 1 (p, 1)}}
 \cr
&\quad  +c \intl_{0}^{t}   \|  \nabla  u(\tau )\|_{ \infty} | V(\tau )|_{ \mathscr{L}^1_{ 1 (p, 1)}}d\tau  
\cr
& \quad + c  \intl_{0}^{t}  | U(\tau )|_{ \mathscr{L}^1_{ 1 (p, 1)}} \| \nabla v(\tau )\|_{\infty} d\tau 
 +  c\intl_{0}^{t}| g(\tau )|_{ \mathscr{L}^1_{ 1 (p, 1)}} d\tau.
  \label{2.16}
 \end{align}       
 Obviously,   for all $ \tau \in (0,t)$,
 \[
| V(\tau )|_{ \mathscr{L}^1_{ 1 (p, 1)}}=|v(\tau )|_{ \mathscr{L}^1_{ 1 (p, 1)}}, 
| U(\tau )|_{ \mathscr{L}^1_{ 1 (p, 1)}}=|u(\tau )|_{ \mathscr{L}^1_{ 1 (p, 1)}}. 
\] 
Furthrmore, thanks to  \eqref{2.11a},  we see that for all $ \tau \in (0,t)$, 
\[
| g(\tau )|_{ \mathscr{L}^1_{ 1 (p, 1)}} = 
 |\nabla \Pi (\PP(u(\tau )),u(\tau ))|_{ \mathscr{L}^1_{ 1 (p, 1)}} \le c
  \| u(\tau )\|^2_{ \mathscr{\dot{L} }^{1}_{1 (p, 1)}}. 
\]
Inserting the estimates above into \eqref{2.16} we are led to 
 \begin{align}
&| v(t)|_{ \mathscr{L}^1_{ 1 (p, 1)}}
\le   c\| v_0\|_{\mathscr{\dot{L} }^{1}_{1 (p, 1)}} +c \intl_{0}^{t}  ( \| u(\tau )\|_{\mathscr{\dot{L} }^{1}_{1 (p, 1)}} \| v(\tau )\|_{\mathscr{\dot{L} }^{1}_{1 (p, 1)}}  + 
\| u(\tau )\|^2_{\mathscr{\dot{L} }^{1}_{1 (p, 1)}} )d\tau.   
  \label{2.17}
 \end{align}       

In order to estimate $ | \dot{P} ^1_{ \infty}(v)|$, we argue as follows. Applying $  \dot{P} ^1_{ \infty}$ to both sides 
of the equation \eqref{2.9}, and using \eqref{2.10}$ _2$, along with \eqref{5.15} and \eqref{5.15a},   we  see that $P_v=\dot{P} ^1_{ \infty}(v) $, $ P_u=\dot{P} ^1_{ \infty}(u) $ 
and $ \widetilde{P} _u=\dot{P} ^1_{ \infty}(\PP(u)) $ solve the following transport equation 
\begin{equation}
 \frac{d}{d\tau  }P_v + (\widetilde{P} _u\cdot \nabla) P_v = \frac{1}{n} 
  \nabla  \widetilde{P}_u: (\nabla P_u))^{ \top} x\quad  \text{ in }\quad  Q_{ T_0}. 
\label{2.18}
\end{equation}    
Note that from the definition of $ \PP$ we get 
\[
\widetilde{P}_u  = \dot{P}^1_\infty(u) - \frac{1}{n} P ^0_\infty(\nabla \cdot u  )x
= P_u - \frac{1}{n} (\nabla \cdot P_u) x,
\]
Applying $\nabla  $ to both side of \eqref{2.18}, we see that $ A:= \nabla P_v\in C^1([0,T_0]; \R^{n^2})$ solves the ODE
\begin{equation}
 \frac{d}{d\tau  }A + \nabla \widetilde{P} _u\cdot A = \frac{1}{n} 
  \nabla  \widetilde{P}_u: (\nabla P_u))^{ \top}I \quad  \text{ in }\quad  (0, T_0). 
\label{2.20}
\end{equation} 
Multiplying both sides by $ \frac{A(\tau )}{| A(\tau )|}$, integating the result over $ (0,t)$, $ t\in (0, T_0]$ and  applying integration by parts, we obtain 
\begin{align}
| A(t)| &\le | A(0)| + c \intl_{0}^{t}  | A(\tau )|\, | \nabla \widetilde{P} _u(\tau )| + 
| \nabla  \widetilde{P}_u(\tau ) : (\nabla P_u(\tau )  ))^{ \top}| d\tau
\cr
&\le | A(0)| + c \intl_{0}^{t}  | A(\tau )|\, | P ^0_{ \infty}(\nabla u(\tau ))|  + 
| P ^0_{ \infty}(\nabla u(\tau )  )| ^2 d\tau.  
\label{2.21}
\end{align}
Combining \eqref{2.17} and \eqref{2.21},  we obtain for all $ t\in (0,T_0]$ 
 \begin{align}
& \| v(t)\|_{ \mathscr{\dot{L}}^1_{ 1 (p, 1)}}
\le   c\| v_0\|_{\mathscr{\dot{L} }^{1}_{1 (p, 1)}} +c \intl_{0}^{t}  ( \| u(\tau )\|_{\mathscr{\dot{L} }^{1}_{1 (p, 1)}} \| v(\tau )\|_{\mathscr{\dot{L} }^{1}_{1 (p, 1)}}  + 
\| u(\tau )\|^2_{\mathscr{\dot{L} }^{1}_{1 (p, 1)}} )d\tau.   
  \label{2.22}
 \end{align}       
Using Gronwall's Lemma, we infer from \eqref{2.22} for all $ t\in [0, T_0]$
\begin{align}
& \| v(t)\|_{\mathscr{\dot{L} }^{1}_{1 (p, 1)}}
\cr
&\quad \le   c_0\bigg\{\| v_0\|_{\mathscr{\dot{L} }^{1}_{1 (p, 1)}} + \intl_{0}^{T_0}  \| u(\tau )\|^2_{\mathscr{\dot{L} }^{1}_{1 (p, 1)}}  d\tau\bigg\}  
\exp\bigg(c_0\intl_{0}^{T_0}   \| u(\tau )\|_{\mathscr{\dot{L} }^{1}_{1 (p, 1)}} d\tau\bigg).   
\label{2.22a}
\end{align}
Assume that $ \| u\|^2_{L^\infty(0, T_0; \mathscr{\dot{L} }^{1}_{1 (p, 1)})} \le 2c_0 e\| v_0\|_{ \mathscr{\dot{L} }^{1}_{1 (p, 1)}} $, 
and 
\[
T_0 = \frac{1}{8c_0^2e^2 \| v_0\|_{\mathscr{\dot{L} }^{1}_{1 (p, 1)} }}. 
\]
Then, \eqref{2.22a} gives 
\begin{align*}
 &\|  v \|^2_{L^\infty(0, T_0; \mathscr{\dot{L} }^{1}_{1 (p, 1)})} 
 \\
  &\quad \le   c_0\Big\{\| v_0\|_{\mathscr{\dot{L} }^{1}_{1 (p, 1)}} + T_0 4c_0^2 e^2  \| v_0\|_{\mathscr{\dot{L} }^{1}_{1 (p, 1)}}^2 )\Big\}  
\exp\bigg(2c_0^2 eT_0 \| v_0\|_{ \mathscr{\dot{L} }^{1}_{1 (p, 1)}}\bigg) 
\\
&\quad \le 2 c_0 e\| v_0\|_{\mathscr{\dot{L} }^{1}_{1 (p, 1)}}.    
\end{align*}  
This shows that $ \mathcal{T}|_{ M}: M \rightarrow M$, where
\[
M = \Big\{ u\in L^\infty(0, T_0; \mathscr{\dot{L} }^{1}_{1 (p, 1)}(\R^{n})) \,\Big|\, 
\| u\|_{ L^\infty(0, T_0; \mathscr{\dot{L} }^{1}_{1 (p, 1)})} \le 2c_0e 
\| v_0\|_{\mathscr{\dot{L} }^{1}_{1 (p, 1)}}\Big\}. 
\]

{\it Proof that $ \mathcal{T}|_{ M}$ is a contractive.}
Let $ u_1, u_2 \in M$ be given. Set $ v_i =T(u_i), i=1,2$, and define $ w= v_1-v_2$. Then $ w $ solves the transport equation 
\begin{equation}
\begin{cases}
\partial _t w + (\PP(u_1)\cdot \nabla)  w 
\\
= -  (\PP(u_1-u_2)\cdot  \nabla) v_2 - \nabla \Pi  (\PP(u_1), u_1-u_2)-  \nabla \Pi  (\PP(u_1-u_2), u_2)
\\
\text{ in}\quad  Q_{ T_0}
\\[0.3cm]
w=0\quad  \text{ on}\quad  \R^{n}\times \{ 0\}.
\end{cases}
\label{2.24}
\end{equation}
Arguing as above,  we get the estimate 
\begin{align*}
\| w(t)\|_{ \mathscr{\dot{L} }^{1}_{1 (p, 1)}} &\le  c  \intl_{0}^{t}  \|u_1(\tau )\|_{ \mathscr{\dot{L} }^{1}_{1 (p, 1)}}\| w(\tau )\|_{ \mathscr{\dot{L} }^{1}_{1 (p, 1)}} d\tau 
\\
&\quad  +  \intl_{0}^{t}  \| u_1(\tau )- u_2(\tau )\|_{ \mathscr{\dot{L} }^{1}_{1 (p, 1)}}   \|v_1(\tau )\|_{ \mathscr{\dot{L} }^{1}_{1 (p, 1)}}d\tau 
\\
&\quad +  \intl_{0}^{t}  (\|u_1(\tau )\|_{ \mathscr{\dot{L} }^{1}_{1 (p, 1)}} +\|u_2(\tau )\|_{ \mathscr{\dot{L} }^{1}_{1 (p, 1)}})
 \| u_1(\tau )- u_2(\tau )\|_{ \mathscr{\dot{L} }^{1}_{1 (p, 1)}}  d\tau.  
\end{align*}
Applying Gronwall's lemma, we arrive at 
\begin{align*}
&\| w(t)\|_{ \mathscr{\dot{L} }^{1}_{1 (p, 1)}} 
\\
&\le   c_0\intl_{0}^{T_0}  (\|u_1(\tau )\|_{ \mathscr{\dot{L} }^{1}_{1 (p, 1)}} +\|u_2(\tau )\|_{ \mathscr{\dot{L} }^{1}_{1 (p, 1)}}
+\|v_2(\tau )\|_{ \mathscr{\dot{L} }^{1}_{1 (p, 1)}})
 \| u_1(\tau )- u_2(\tau )\|_{ \mathscr{\dot{L} }^{1}_{1 (p, 1)}}  d\tau \times 
\\
&\quad\times \exp\bigg(c_0\intl_{0}^{T_0}  \|u_1(\tau )\|_{ \mathscr{\dot{L} }^{1}_{1 (p, 1)}} \bigg)
 d\tau 
 \\
 & \le 6 c_0^2 e^2 T_0  \| v_0\|_{  \mathscr{\dot{L} }^{1}_{1 (p, 1)}} \| u_1- u_2\|_{ L^\infty(0, T_0;  \mathscr{\dot{L} }^{1}_{1 (p, 1)})}
 \le \frac{2}{3} \| u_1- u_2\|_{ L^\infty(0, T_0;  \mathscr{\dot{L} }^{1}_{1 (p, 1)})}. 
\end{align*}

By virtue of Banach's fixed point theorem there exists a unique fixed point $ v\in M$ such that $ \mathcal{T}(v)=v$. 

\vspace{0.3cm}
In order to verify that $ v$ is a solution to  \eqref{euler} it only remains to show that $ \nabla \cdot v= 0$ or what is equivalent to
$ \PP(v)=v$.  First note that due to the definition of  $ \mathcal{T}$,  the function $ v\in L^\infty(0, T_0; \mathscr{\dot{L} }^{1}_{1 (p, 1)}(\R^{n}))$ 
solves the transport equation 
\begin{equation}
\begin{cases}
\partial _t v + (\PP(v)\cdot \nabla) v = -\nabla \Pi (\PP(v),v)\quad  \text{ in}\quad  \R^{n},
\\[0.3cm]
v=v_0\quad \text{ on}\quad  \R^{n}\times \{ 0\}. 
\end{cases}
\label{2.25}
\end{equation}
Applying $ \nabla \cdot $ to both sides of \eqref{2.25}, we see that 
\[
\begin{cases}
\partial _t \nabla \cdot v + (\PP(v)\cdot \nabla) \nabla \cdot v + \nabla \PP(v):(\nabla v)^{ \top} = \nabla \PP(v):(\nabla v)^{ \top}\quad  \text{ in}\quad  \R^{n},
\\[0.3cm]
\nabla \cdot v=0\quad \text{ on}\quad  \R^{n}\times \{ 0\}. 
\end{cases}
\]
Accordingly, $ \nabla \cdot v \in L^\infty(Q_{ T_0})$ solves the transport equation with zero data. 
The strong-weak uniqueness \cite[pp. 46]{trans} implies $ \nabla \cdot v=0$.  This completes the proof of Theorem\,\ref{thm1}.  

\hfill \Beweisende 

\section{Proof of Theorem\,\ref{thm2}}
\label{sec:-9}
\setcounter{secnum}{\value{section} \setcounter{equation}{0}
\renewcommand{\theequation}{\mbox{\arabic{secnum}.\arabic{equation}}}}

Let $ v_0 \in \mathscr{L }^{1}_{1 (p, 1)}(\R^{n})$. Then $ v_0 - v_0(0) \in \mathscr{\dot{L} }^{1}_{1 (p, 1)}(\R^{n})$. 
According to Theorem\,\ref{thm1} there exists 
$ T_0 \ge  \frac{1}{c\|v_0- v_0(0)\|_{ \mathscr{\dot{L} }^{1}_{1 (p, 1)}}}$ and a unique solution 
 $ \widetilde{v}\in L^\infty(0, T_0; \mathscr{\dot{L} }^{1}_{1 (p, 1)}(\R^{n}))$ to the Euler equations  \eqref{euler}, 
 \eqref{initial} with $ v_0- v_0(0)$ in place of $v_0$ and pressure 
 $ \widetilde{p}\in L^\infty(0,T_0; L^2_{loc}(\R^{n} ))$ such that $ \nabla \widetilde{p}\in L^\infty(0, T_0; \mathscr{\dot{L} }^{1}_{1 (p, 1)}(\R^{n}))$
 \begin{equation}
\nabla \widetilde{p} = \nabla \Pi(\widetilde{v}, \widetilde{v}).    
 \label{9.1}
  \end{equation} 
 Setting $ v(x,t) = \widetilde{v} (x- t v_0(0), t) + v_0(0)$, we see that $ v\in L^\infty(0, T_0; \mathscr{L }^{1}_{1 (p, 1)}(\R^{n}))$ and solves the Euler equations  \eqref{euler},  \eqref{initial} with pressure 
$ p(x,t)= \widetilde{p}(x- v_0(0)t, t)$.  We now verify that for almost all $ t\in (0,T_0)$
\[
\nabla p(t) - \nabla p(0,t) = \nabla \Pi(v(t), v(t)).   
\]
Clearly, as $ \nabla \cdot v=0$ in $ Q_{ T_0}$ we find that $ \pi = p(t)- \nabla p(0, t)\cdot x$ solves the Poisson equation 
\[
\Delta \pi  = \nabla v(t) : (\nabla v)^{ \top}= \nabla \cdot ((v(t) \cdot \nabla) v(t)).
\]
Obviously, it holds $ \nabla \pi (0)= 0$. It only remains to verify the asymptotics as $ |x| \rightarrow +\infty$. 
By the definition of $ p$ along with  \eqref{9.1},  recalling the definition of $ \nabla \Pi $,  it follows that 
\begin{align*}
\dot{P}^1_\infty(\nabla \pi ) &= \dot{P}^1_\infty(\nabla p (t))= \dot{P}^1_\infty(\nabla \widetilde{p} (t))
= \frac{1}{n} P^0_\infty(\nabla \widetilde{v}(t): (\nabla \widetilde{v}(t))^{ \top})
\\
&= \frac{1}{n} P^0_\infty(\nabla v(t): (\nabla v(t))^{ \top}).
 \end{align*}

\section{Proof of Theorem\,\ref{thm3}}
\label{sec:-?}
\setcounter{secnum}{\value{section} \setcounter{equation}{0}
\renewcommand{\theequation}{\mbox{\arabic{secnum}.\arabic{equation}}}}

Let $ 0<T<+\infty$. We begin our discussion with the following oscillation estimate.  
Let $ v\in L^1(0, T; C^{ 0,1}(\R^{n}))$.  Let $ x_0\in \R^{n}$. We call $ \xi \in C^1([0,T])$ an 
{\it cancelling shift 
in $ x_0$} if $ \xi $ satisfies the following ODE
\begin{equation}
\dot{\xi} (t) =v(x_0+ \xi (t), t)\quad  \forall\,t\in (0,T).   
\label{7.12a}
\end{equation}

\begin{thm}
\label{thm7.1}
Let $1 \le q \le +\infty, N\in \N_0  , s\in [0, N+1)$,  and let $ P\in 
L^\infty(0,T; \mathcal{P}_{ N})$. Let  $ v\in L^{ \infty}(0, T; C^{ 0,1}(\R^{n}))$ be a solution to the  Euler equations   
\begin{equation}
\begin{cases}
\nabla \cdot v =0 \quad  \text{ in}\quad  Q_T,
\\[0.3cm]
\partial _t v + (v\cdot \nabla) v = - \nabla \Pi (v,v)+ P\quad  \text{ in}\quad  Q_T,
\\[0.3cm]
v=v_0 \quad  \text{on}\quad \R^{n}\times \{ 0\},
\end{cases}
\label{m-euler}
\end{equation}
with initial value  $v_0\in   \mathscr{L}^{s}_{ q (p, N)}(\R^{n})$. 
Then $ v \in L^{ \infty}(0, T; \mathscr{L}^{s}_{ q (p, N)}(\R^{n}))$. 

Furthermore, the following estimate holds true for all $ x_0\in \R^{n}$ and $ t\in [0,T]$,  and all characteristics   $ \xi \in C^{ 1,1}([0,T])$ satisfying \eqref{7.12a}
\begin{align}
&\osc_{ p, \max\{2N-1, N \}}(v(t ), x_0 + \xi (t)) 
\cr
&\le cS_{ N+1, 1} (\osc_{ p, N}(v_0, x_0 + \xi (0))) +  c\intl_{0}^{t} 
\| \nabla v(\tau )\|_{ \infty}  S_{ N+1, 1} (\osc_{ p, N}(v(\tau  ), x_0 + \xi (\tau )))  d\tau.
\label{7.12b}
\end{align}
Furthermore, it holds
 \begin{equation}
 | v(t)|_{  \mathscr{L}^{ s}_{q (p, N) }}\le c | v_0|_{  \mathscr{L}^{ s}_{q (p, N) }} \exp \bigg(c\intl_{0}^{t} 
\| \nabla v(\tau ) \|_{ \infty}  d\tau\bigg). 
 \label{7.28a}
 \end{equation}

\end{thm}

{\bf Proof}: 1. First, let us consider the case $ N \ge 1$.  Let $ x_0\in \R^{n}$ be fixed.  Let $ \xi \in C^{ 1,1}([0,T])$ satisfying \eqref{7.12a}.  
We set 
\begin{align*}
& V(x, \tau ) = v(x + \xi (\tau ), \tau ) - \dot{\xi} (\tau ),\quad  \Pi (x, \tau )=  \pi (x + \xi (\tau ), \tau ),
\\
&\qquad \widetilde{P} (x, \tau ) = P(x + \xi (\tau ), \tau ),\quad  (x, \tau )\in Q_T,  
\end{align*}
where $ \nabla \pi = \nabla \Pi (v, v)$. By the definition of $ \Pi $ it follows that 
\begin{equation}
 \nabla \Pi = \nabla \Pi (V, V).
\label{7.12}
\end{equation}

Clearly, $ (V, P)$ solves the transformed  Euler equations
\begin{equation}
\partial _t V + (V\cdot \nabla) V =- \nabla \Pi + \widetilde{P}  - \ddot{\xi}  \quad  \text{ in}\quad Q_T.
\label{7.13}
\end{equation}
 According to  \eqref{2.4} in Lemma\,\ref{lem2.1} with  $ v=u=V, g = -\nabla \Pi  +\widetilde{P} - \ddot{\xi} $  and $ r= 2^{ j+1}$ we get 
 \begin{align}
& \osc_{p, 2N-1\}} (V(t); x_0, 2^{ j})
\le   c \osc_{p, N} (V(0); x_0,  2^{ j+1})    
\cr
& \quad + c \intl_{0}^{t}  \| \nabla v(\tau )\|_{\infty}  \osc_{p, N} (V(\tau ); x_0,  2^{ j+2})d\tau 
+ c\intl_{0}^{t} \osc_{p, N} (\nabla \Pi (\tau ); x_0,  2^{ j+1})
d\tau.
  \label{7.14}
 \end{align}       
Thanks to \eqref{6.44a} (cf. Remark\,\ref{rem6.5a}) with $u=v=V $ it holds 
\begin{align}
& \osc_{p, N}(\nabla \Pi (\tau ); x_0, 2^{ j+1}) 
\le  c\| \nabla v(\tau )\|_{ \infty}  S_{ N+1, 1} ( \osc_{p, N}(V(\tau ); x_0))_j,\quad j\in \Z.
\label{7.14a} 
\end{align}

Inserting \eqref{7.14a} into the last integral on the right-hand side of \eqref{7.14}, 
and operating $ S_{s', 1}$ for some $ \max\{s,1\}< s'< N+1$, using Lemma\,\ref{lem2.1}, we arrive at 
 \begin{align}
& S_{ s', 1} (\osc_{p, 2N-1} (V(t); x_0))_j
\cr
&\le   c S_{ s', 1}(\osc_{p, N} (V(0); x_0) )_j   
+ c \intl_{0}^{t}  \| \nabla v(\tau )\|_{\infty}  S_{ s', 1}(\osc_{p, N} (V(\tau ); x_0))_j d\tau. 
  \label{7.15}
 \end{align}  
 Appealing to \cite[Corollary\,3.10]{trans}, we get 
 \begin{equation}
\osc_{p, N} (V(t); x_0))_j \le cS_{ N+1, 1} (\osc_{p, 2N-1} (V(t); x_0))_j \le cS_{ s', 1} (\osc_{p, 2N-1} (V(t); x_0))_j.  
 \label{7.15a}
  \end{equation} 
Estimating the left-hand side of  \eqref{7.15} by  \eqref{7.15a}, operating $ S_{s'', 1}$ for some $ \max\{s,1\}< s''< N+1$, again appealing to  Lemma\,\ref{lem2.1}, we find  for all $ j\in \Z$
 \begin{align}
& S_{ s'', 1} (\osc_{p, N} (V(t); x_0))_j
\cr
&\le   c S_{ s'', 1}(\osc_{p, N} (V(0); x_0) )_j   
+ c \intl_{0}^{t}  \| \nabla v(\tau )\|_{\infty}  S_{ s'', 1}(\osc_{p, N} (V(\tau ); x_0))_j d\tau. 
  \label{7.24}
 \end{align}

 By the aid of  Gronwall's lemma  we obtain from \eqref{7.24} for all $ t\in [0, T]$
\begin{align}
 \osc_{p, N} (V(t); x_0) &\le S_{ s'', 1} ( \osc_{p, N} (V(t); x_0)) 
 \cr
 &\le cS_{s'', N} ( \osc_{p, N} (V(0); x_0)) \exp \bigg(c\intl_{0}^{t} 
\| \nabla v(\tau ) \|_{ \infty}  d\tau\bigg). 
\label{7.25}
\end{align} 
Let $ t\in (0, T]$ be arbitrarily chosen, but fixed.   Since the constant in the above estimate does 
not depend on the choice of the characteristic, we may choose $ \xi \in C^{ 1,1}([0,T])$ such that $ \xi (t)=0$. Then $V(t)= v(t)$. 
Thus, replacing in \eqref{7.25} $ V(t)$ by $ v(t)$,  operating  $ S_{ s, q}$ to both sides of \eqref{7.25}, multiplying the result by $ 2^{ -js}$,   
and using \eqref{10.3} 
of Lemma\,\ref{lem10.1}, we arrive at  
\begin{align}
\| \{2^{ -sj}\osc_{p, N} (v(t); x_0, 2^j)\}\|_{ \ell^q} \le c | v_0|_{  \mathscr{L}^{ s}_{q (p, N) }} \exp \bigg(c\intl_{0}^{t} 
\| \nabla v(\tau ) \|_{ \infty}  d\tau\bigg). 
\label{7.26}
\end{align}    
 In \eqref{7.26}  taking the supremum over all $ x_0\in \R^{n}$,  we obtain 
 \begin{equation}
 | v(t)|_{  \mathscr{L}^{ s}_{q (p, N) }}\le c | v_0|_{  \mathscr{L}^{ s}_{q (p, N) }} \exp \bigg(c\intl_{0}^{t} 
\| \nabla v(\tau ) \|_{ \infty}  d\tau\bigg). 
 \label{7.27}
 \end{equation}
Whence,   $ v\in L^\infty(0, T; \mathscr{L}^{ s}_{ q (p, N)}(\R^{n}))$. 

\vspace{0.3cm}  
2. The case $ N=0, s\in [0, 1)$.  Noting that $ \mathscr{L}^{ s}_{q (p, 0) }(\R^{n}) \hookrightarrow  \mathscr{L}^{ s}_{q (p, 1) }(\R^{n}) $, 
we get from the case $ N=1$  that $ v\in L^\infty(0, T;  \mathscr{L}^{ s}_{q (p, 1) }(\R^{n}))$ together with the estimate \eqref{7.27}. 
 On the other hand, in view of Lemma\,\ref{lem5.9},  from  $ v_0\in \mathscr{L}^{ s}_{ q (p, 0)}(\R^{n})$ we deduce that 
\[
\dot{P}^1_\infty(v_0)=0.  
\]
Applying, $ \dot{P}^1_\infty$ to \eqref{m-euler}, and  $P \in L^\infty(0,T;\mathcal{P}_0)$, using \eqref{5.15}, \eqref{5.15a} and \eqref{6.49},     we see that 
$ Q(t)=\dot{P}^1_\infty(v(t ))$, solves the equations in $ Q_T$ 
\[
\partial _t Q + (Q\cdot \nabla) Q =  \frac{1}{n} \nabla Q: \nabla Q^{ \top}x,\quad  Q(0)=0. 
\]
Applying $ \nabla $ to the above equation, we see that $ A(t) = \nabla Q(t)$ solves the ODE
\[
\partial _t A + A^2 =  \frac{1}{n} {\rm tr} (A^2)\quad  \text{in}\quad (0,T),\quad  A(0)=0. 
\]
With the help of Gronwall's lemma, we easily get $ A(t)=0$ for all $ t\in [0,T]$. Once more appealing to Lemma\,\ref{lem5.9} 
it follows that $ v(t)\in \mathscr{L}^{ s}_{ q (p, 0)}(\R^{n})$, and  \eqref{7.27} together with \eqref{5.24} implies 
 \begin{equation}
 | v(t)|_{  \mathscr{L}^{ s}_{q (p, 0) }}\le c | v_0|_{  \mathscr{L}^{ s}_{q (p, 0) }} \exp \bigg(c\intl_{0}^{t} 
\| \nabla v(\tau ) \|_{ \infty}  d\tau\bigg). 
 \label{7.28}
 \end{equation}
Whence, $ v\in L^\infty(0,T; \mathscr{L}^{ s}_{q (p, 0) }(\R^{n}))$  and \eqref{7.12b} holds.  \hfill \Beweisende  

\begin{rem}
 \label{rem7.3} 1. Firsly, we wish to remark that  Theorem\,\ref{thm7.1} still holds under weaker assumption 
 $ v\in L^1(0,T; C^{ 0,1}(\R^{n} ))$ together with the assumption $ v_0\in \mathscr{L}^s_{ q (p, N), \sigma }(\R^{n})\cap \mathscr{L}^1_{ 1 (p, 1), \sigma }(\R^{n})$. In fact, from Theorem\,\ref{thm2} we get 
  $ v\in L^\infty_{ loc}([0,T_{ \ast}); \mathscr{L}^1_{ 1 (p, 1), \sigma }(\R^{n}))$ for a maximal time $ T_{ \ast}>0$. 
  Assume $ T_{ \ast} \le T$. Then thanks to  \eqref{7.28a} we get $ v\in 
  L^\infty(0,T_{ \ast}; \mathscr{L}^1_{ 1 (p, 1), \sigma }(\R^{n}))$. In case $ T=T_{ \ast}$ we get the claim. In the other case since $ v(T_{ \ast}) \in  \mathscr{L}^1_{ 1 (p, 1), \sigma }(\R^{n})$ we are in a position to apply again  
 Theorem\,\ref{thm2}, which shows that $ L^\infty(0,T_{ \ast}+\delta ; \mathscr{L}^1_{ 1 (p, 1), \sigma }(\R^{n}))$ 
 for some $ \delta >0$, which clearly contradicts to the definition of $ T_{ \ast}$. Whence, the claim. 
 
 \vspace{0.2cm}
 2. As Corollary of the first remark we get the local well-posedness of the Euler equations in $ \mathscr{L}^s_{ q (p, N), \sigma }(\R^{n})\cap \mathscr{L}^1_{ 1 (p, 1), \sigma }(\R^{n})$ for $ N\in \N_0, 1 <p< +\infty, 1 \le q \le +\infty, s\in [0, N+1)$.

\end{rem}

\vspace{0.3cm}
Next, we provide the  following uniqueness result
\begin{lem}
\label{lem3.1}
Let $ (v, p), (u, q)\in L^\infty(0, T; \mathscr{L}^1_{ 1 (p, 1), \sigma }(\R^{n}))$ are two  solution to  \eqref{euler}, \eqref{initial}. 
Assume that $ \dot{P}^0_{\infty} (\nabla v_0) = 0$, and 
\begin{equation}
P^0_\infty(D ^2 p) = P^0_\infty(D ^2 q) =0\quad  \text{ in}\quad  (0,T). 
\label{3.0a}
\end{equation}
Then $ (v, p)= (u, q)$.  
\end{lem}

{\bf Proof}:  Set $ A(t) = P^0_\infty(\nabla u(t)), t\in (0,T)$. As in the proof of Theorem\,\ref{thm1} we get 
 \[
\dot{A} + A^2 =  P^0_\infty(D ^2 p)=0\quad  \text{ in}\quad  (0,T). 
\]
Owing to $ A(0) = \dot{P}^0_{\infty} (\nabla v_0)=0$ it follows that $ A(t)=0$ for all $ t\in (0,T)$. Whence, 
$ \dot{P}^0_{\infty} (\nabla v)=0$ in $ (0,T)$. Analogously, we see that $ \dot{P}^0_{\infty} (\nabla u)=0$ in $ (0,T)$. 

Next, let $ \xi \in C^{ 1,1}([0,T])$ be the unique solution to the ODE
\[
\dot {\xi}(t)= v(\xi (t), t)\quad  t\in (0,T),\quad  \xi (0)=0.   
\]
Set 
\[
V(x,t)=v(\xi (t), t)-  \dot {\xi}(t), \quad  P(x,t)=p(\xi (t), t),\quad  (x,t ) \in Q_{ T}. 
\]
Clearly, $ (V, P) \sim (v, p)$ and $ (V, P)$ is a centered solution to \eqref{euler}, \eqref{initial}. In addition, 
$ P^0_\infty(D^2 p)=0$ implies  $ P^0_\infty(D^2 P)=0$, and $ P^0_\infty(\nabla  v)=0$ implies 
$  P^0_\infty(\nabla  V)=0$. Hence, 
\[
\dot{P} ^1_\infty( \nabla P ) = 0 = -\frac{1}{n}  P^0_\infty(\nabla  v: (\nabla v)^{ \top}) x. 
\]
Noting that \eqref{euler} and $ V(0, t)=0$ implies $ \nabla P(0, t)=0$ we  infer $P= \Pi (V, V)$. This shows that 
$ (V, P)$ is eligible .  Similar, there exists a unique centered solution $ (U, Q) \sim (u, q)$. Since 
$ P^0_\infty(D^2 q)=0$ and $ P^0_\infty(\nabla  u)=0$ this solution is eligible  too.  According to  
Theorem\,\ref{thm1} this solutions are unique, which gives $ (V, P)= (U, Q)$. Accordingly, $ (v,p)= (u, q)$.  \hfill \Beweisende

\vspace{0.5cm}  
{\bf Proof of Theorem\,\ref{thm3}:} 
Let $ v_0\in \mathscr{L}^1_{ 1 (p, 1), \sigma }(\R^{n})\cap BMO$.  Set $ u_0= 
v_0-v_0(0) \in \mathscr{\dot{L} }^{1}_{1 (p,1), \sigma}(\R^{n})\cap BMO$. 
According to  Theorem\,\ref{thm1} there exists a unique centered  eligible    solution $ (u, \pi )\in 
L^\infty(0, T_0; \mathscr{\dot{L} }^{1}_{1 (p,1), \sigma}(\R^{n})\times L^2_{ loc}(\R^{n}))$ to \eqref{euler}
 with $ u_0$ in place of $ v_0$, where $ T_0>0$ satisfies $ T_0 \ge \frac{1}{c \| u_0\|_{\mathscr{\dot{L} }^1_{ 1 (1)}}}= 
\frac{1}{c | v_0|_{ \mathscr{L}^1_{ 1 (p, 1)}} }. $  By the definition of eligible  centered  solutions 
  to \eqref{euler} it holds 
\begin{equation}
\nabla \pi(t) =  \nabla \Pi (u(t), u(t)) \quad \forall\,t\in [0,T]. 
\label{3.1}
\end{equation}
 As it has been proved in Theorem\,\ref{thm7.1},  $ \dot{P}  ^1_{ \infty}(u_0) =0 $ implies 
  $  \dot{P}  ^1_{ \infty}(u(\tau )) =0 $ for all $ \tau \in (0, T_0)$.  Thus,  it holds 
 \begin{equation}
-\Delta  \pi = \nabla u: (\nabla u)^{ \top}\quad  \text{ in}\quad  \R^{n},\quad  \dot{P}  ^1_{ \infty}(\nabla \pi ) =0 . 
\label{3.2}
\end{equation} 
Let $ x_0 \in \R^{n}$.  In view of   \eqref{7.12b} of Theorem\,\ref{thm7.1} with $ N=0, s=0 $ and $ q=\infty$
we have $ u\in L^\infty(0,T_0; BMO)$ and it  holds   
for all $ t\in (0,T_0)$ 
\begin{align}
 | u(t)|_{BMO}\le c | v_0|_{ BMO} \exp \bigg(c\intl_{0}^{t} 
\| \nabla u(\tau ) \|_{ \infty}  d\tau\bigg). 
 \label{3.7}
 \end{align}
 
Let $ \tau \in (0,T)$. Thanks to Corollary\,\ref{cor6.9} with $ s=0$ and $ q=\infty$ there exists a unique very weak solution 
$ \pi _0(\tau ) \in \mathscr{L}^{ 0}_{ \infty (p,0)}(\R^{n}) \cong BMO$ to 
\[
-\Delta \pi _0(\tau )= \nabla \cdot \nabla \cdot (v(\tau ) \otimes v(\tau ))\quad  \text{ in}\quad  \R^{n}
\]
with $P^0_{ 0,1} (\pi_0 (\tau ))=0$.  This implies that $ P(\tau ):= \pi(\tau ) - \pi _0(\tau ) $ is harmonic. 
By virtue of the  Liouville theorem 
for harmonic functions it follows that $P(\tau ) \in \mathcal{P}_2$.    Noting that $ \dot{P}_\infty^1(\nabla \pi _0(\tau ))=0 $, and   
observing \eqref{3.2}, we see that $\dot{P} _\infty^1(\nabla P(\tau ))=0$. Accordingly, $P(\tau ) \in \mathcal{P}_1$, and $ \nabla P(\tau )$ is constant for all $ \tau \in [0, T_0]$.
Set $ \eta (\tau ) = \nabla P(\tau )$. Define,
\[
\xi (t) =  -\intl_{0}^{t}  \intl_{0}^{\tau }  \eta (s) ds  d\tau - t v_0(0), \quad  t\in [0, T_0].  
\]
Set 
\begin{align*}
v(x,t) &= u(x+ \xi (t), t)- \dot{\xi}  (t),\quad 
\\
p (x,t) &= \pi (x+ \xi (t), t)- P(x+ \xi (t), t)=\pi_0 (x+ \xi (t), t),\quad  (x,t)\in Q_{ T_0}. 
\end{align*}
Clearly, $ (v, p)\in L ^\infty(0, T; \mathscr{L }^{1}_{1 (p, 1), \sigma}(\R^{n} )\cap 
BMO  \times BMO$ solves the Euler equations \eqref{euler}, \eqref{initial}.  Eventually, replacing $ p$ by $ p- P_{ 0,1}^0(p)$ we may assume that $ P^0_{ 0,1}(p)=0$ is satisfied.   Furthermore, thanks to  \eqref{6.66} it holds 
\begin{equation}
\|p(\tau )\|_{ BMO} \le c \|v(\tau )\|^2_{ BMO}\quad \forall \tau \in (0,T_0).  
\label{3.7a}
 \end{equation} 

\vspace{0.2cm}
{\it Uniqueness}. Let $ (\overline{v}, \overline{p} )\in L ^\infty(0, T; \mathscr{L }^{1}_{1 (p, 1), \sigma}(\R^{n} )\cap 
BMO  \times BMO$ 
be another solution to \eqref{euler}, \eqref{initial}, with $ P^0_{ 0,1}(\overline{p} )=0$.  From  Theorem\,\ref{thm6.3} we get  $ \nabla \overline{p } (t)
\in \mathscr{\dot{L} }^{1}_{1 (p, 1)}(\R^{n}) $. Clearly, 
 $ \lim_{m \to \infty}\dot{P} ^2_{ 0, 2^m}(\overline{p} (t))=0$. This shows that $ P^0_\infty(D^2 \overline{p}(t) )=0$ 
 for all $ t\in [0,T]$. Analogously, $ P^0_\infty(D^2 p(t) )=0$  for all $ t\in [0,T]$. Applying 
 Lemma\,\ref{growth} 
 with $ s<1$ we see that $ P^0_\infty(\nabla v_0)=0$. We are now in a position to apply Lemma\,\ref{lem3.1}, which yields
 $ (v, p)= (\overline{v}, \overline{p}  )$.

\vspace{0.3cm}  
It only remains to prove that $ v_0\in L^\infty(\R^{n})$ implies that  $ v\in L^\infty(Q_{ T_0})$. 
 Since $ v\in L^\infty(0,T_0; C^{ 0,1}(\R^{n}))$, to verify the claim  it will be sufficient to prove 
\begin{equation}
|P^0_{ x_0,1}(v(t))  |\le c \bigg\{\| v_0\|_{ \infty} +  \intl_{0}^{t}   | v(\tau )|^2_{  BMO} d\tau \bigg\}\exp  \intl_{0}^{t} \| \nabla v(\tau )\|  d\tau\quad  \forall\,x_0\in \R^{n}.  
\label{3.8}
\end{equation}
In fact, applying $ P^0_{ x_0,1}$ to both sides of \eqref{euler}, we get the identity  
\begin{align}
\frac{d}{d \tau }  P^0_{ x_0,1}(v)  &=- P^0_{ x_0,1}(v\cdot \nabla )v)- P^0_{ x_0,1}(\nabla p).  
 \cr
&=- P^0_{ x_0,1}(v) \cdot P^0_{ x_0,1} (\nabla v)-  \nabla \cdot P^0_{ x_0,1}((v-P^0_{ x_0,1}(v)) \otimes  (v-P^0_{ x_0,1}(v)))
\cr
&\qquad - \nabla P^1_{ x_0,1}( p- P^0_{ x_0,1}(p) ).  
\label{3.9}
\end{align}
The first term can be estimated by $ | P^0_{ x_0,1}(v(\tau ))| \| \nabla v(\tau )\|_{ \infty}$, while the remaining two terms are bounded by
$ | v(\tau )|_{ BMO}^2+ | p(\tau )|_{ BMO} \le c| v(\tau )|_{ BMO}^2 $,  where we have used  \eqref{3.7a}.   Thus, 
\[
\frac{d}{d \tau }  P^0_{ x_0,1}(v) \le P^0_{ x_0,1}(v)  \| \nabla v\|_{ \infty}+  c| v|^2_{ BMO}. 
\] 
Using Gronwall's lemma, we get \eqref{3.8}. We now easily estimate 
\[
| v(x_0,\tau )| \le | v(x_0,\tau ) - P^0_{ x_0,1}(v(\tau ))|  +| P^0_{ x_0,1}(v(\tau ))| 
\le c \| \nabla v(\tau )\|_{ \infty} + | P^0_{ x_0,1}(v(\tau ))|. 
\]   
Together with \eqref{3.8} we see that $ v\in L^\infty(Q_{ T_0})$.  
This completes the proof of Theorem\,\ref{thm3}. 

 \hfill \Beweisende  

\section{Proof of Theorem\,\ref{thm4}}
\label{sec:-?}
\setcounter{secnum}{\value{section} \setcounter{equation}{0}
\renewcommand{\theequation}{\mbox{\arabic{secnum}.\arabic{equation}}}}

The proof of Theorem\,\ref{thm4} will be carried out, using logarithmic Sobolev  type inequality similarly to the decaying case.
We provide such inequality for the space $ \mathscr{L}^{ 1+\delta }_{ 1 (p, 1)}(\R^{n})$. 
 
 \begin{lem}[Logarithmic inequality]
 \label{lem7.3}
 Let  
 $ u \in  \mathscr{L}^{ 1+\delta }_{ 1 (p, 1)}(\R^{n}) \cap BMO_1,  \delta \in (0,1)$. Then for all $ x_0\in \R^{n}$ and all $ k\in \Z$ it holds 
 \begin{equation}
  \sum_{j=-\infty}^{ k } 2^{ -j}\osc_{p, 1} (u; x_0, 2^j)
 \le c 2^{ \delta k}+ c| \nabla u|_{ BMO}\log (1+ | u|_{ \mathscr{L}^{ 1+\delta }_{ 1 (p, 1)}}). 
 \label{7.38}
 \end{equation}
 In particular, for all $ m,k\in \N$ with $m<k$ it holds 
 \begin{align}
  |P^0_{ x_0, 2^m}(\nabla u)| \le   c 2^{ \delta k}+ | P^0_{ x_0, 2^k}(\nabla u)|+  c| \nabla u|_{ BMO}\log (1+ | u|_{ \mathscr{L}^{ 1+\delta }_{ 1 (p, 1)}}).  
 \label{7.39}
 \end{align}
 \end{lem}

{\bf Proof}:  1. Let $k  \in \Z$, and let $ l\in \N$, specified below. Using H\"older's inequality and Poincar\'e's inequality,  we easily get 
\begin{align*}
&\sum_{j=-\infty}^{ k} 2^{ -j}\osc_{p, 1} (u; x_0, 2^j)
\\
&= \sum_{j=-\infty}^{ k-l} 2^{ -j}\osc_{p, 1} (u; x_0, 2^j)+ 
\sum_{j=k-l+1}^{ k} 2^{ -j}\osc_{p, 1} (u; x_0, 2^j)
\\
& \le  \frac{2^{ \delta (k-l)}}{1-2^\delta }| u|_{ \mathscr{L}^{ 1+\delta }_{ 1 (p, 1)}}+ 
cl | \nabla u|_{ BMO}. 
\end{align*}
Choosing $ l =  \Big\lfloor \frac{1}{\log 2}\log (1+ | u|_{ \mathscr{L}^{ 1+\delta }_{ 1 (p, 1)}} )\Big\rfloor+1$, we infer from the above estimate  
 \begin{align*}
 &\sum_{j=-\infty}^{ k } 2^{ -j}\osc_{p, 1} (u; x_0, 2^j)
 \le c 2^{ \delta k}+ c| \nabla u|_{ BMO}\log (1+ | u|_{ \mathscr{L}^{ 1+\delta }_{ 1 (p, 1)}}). 
 \end{align*}
Whence, \eqref{7.38}.  
 
 \vspace{0.2cm}
 2. Let $ m, k \in \Z, m<k$. Arguing as in the proof of Theorem\,\ref{thm5.4}, and using \eqref{7.38},   we estimate  
 \begin{align}
| P^0_{ x_0, 2^m}(\nabla u)| &\le | P^0_{ x_0, 2^k}(\nabla u)| + \sum_{j=m}^{ k}  2^{ -j}\osc_{p, 1} (u; x_0, 2^j) 
 \cr
 &\le 
  c 2^{ \delta k}+| P^0_{ x_0, 2^k}(\nabla u)|+  c| \nabla u|_{ BMO}\log (1+ | u|_{ \mathscr{L}^{ 1+\delta }_{ 1 (p, 1)}}).  
 \label{7.39a}
 \end{align}
 This completes the proof of  \eqref{7.39}.  
 \hfill \Beweisende  
 
 \vspace{0.3cm}
{\bf  Proof of Theorem\,\ref{thm4}}:  
 1. First applying the known Calderon-Zygmund estimate in BMO to the Biot-Savart formula, we get the estimate 
\begin{equation}
| \nabla  v(\tau )|_{ BMO}  \le c | \omega(\tau ) |_{ BMO}\quad \forall \tau \in [0,T_{ \ast}). 
\label{7.51a}
\end{equation}

2. Let $ x_0\in \R^{n}$ be fixed. 
Let $ k\in \Z$ be appropriately chosen,   which will be specified below. Our  aim is to provide an uniform bound for $ \sup_{ j \ge k}|P^0_{ x_0, 2^j}(\nabla v)|$.
Let $ t\in (0,T_{ \ast})$ be fixed. Let $ \xi \in C^{ 1,1}([0, T_{ \ast}])$ be a characteristic such that 
\[
\dot{\xi}(t) = P^0_{ x_0, 2^k}(v(\cdot +\xi (t), t))\quad  \forall\,t\in (0, T_{ \ast}),\quad  \xi (t)=0. 
\]   
We set 
\[
V(x,t) = v(x+ \xi (t), t) - \dot{\xi}(t),\quad  \Pi (x,t) = \pi (x+ \xi (t), t) + \ddot{\xi}(t)x,\quad  (x,t)\in Q_{ T_{ \ast}}.  
\]
Clearly,
\[
P^0_{ x_0, 2^k} (V(t))=0\quad  \forall\,t\in (0, T_{ \ast}), 
\]
and $ (V, \Pi )$ solves the Euler equations. 
\begin{equation}
\begin{cases}
\nabla \cdot V =0\quad  \text{ in}\quad Q_{ T_{ \ast}},  
\\[0.3cm]
\partial _t V + (V\cdot \nabla) V = - \nabla \Pi\quad  \text{ in}\quad Q_{ T_{ \ast}}.   
\end{cases}
\label{7.53}
\end{equation}

In view  \eqref{2.4a}, putting $ f= V, \pi = \Pi $ and $ r= 2^{ j+2}$ therein, we find 
 \begin{align}
& \osc_{ 2, 1} (v(t); x_0,  2^j)
\le   c \osc_{2, 1} (v_0; x_0,  2^{ j+1}) +c 2^{- j} \intl_{0}^{t} \|   V(\tau )\|_{ L^\infty(B(x_0,2^{ j+1}))}  
\osc_{ 2, 1} (V(\tau ); x_0,  2^{ j+2})d\tau  
\cr
& \quad +  c  \intl_{0}^{t}  \osc_{ 2, 1} (V(\tau ); x_0, 2^{ j+1})  | P^{ 0}_{ x_0, 2^{ j+1}}(\nabla V(\tau )| d\tau 
\cr
&\quad +  c\intl_{0}^{t} \osc_{ \frac{2n}{n+2}, 1} (\nabla \Pi (\tau ); x_0,  r)  d\tau.
  \label{7.53e}
 \end{align}     
 Using  Poincare's inequality along with \eqref{7.51a}, we find for all $ i \ge j$  
\begin{equation} 
| P^0_{ x_0, 2^i}(V(\tau )) - P^0_{ x_0, 2^{ i+1}}(V(\tau )) | \le c  2^i |\omega (\tau )|_{ BMO}+ c2^i\sup_{ i \ge j} 
| P^{ 0}_{ x_0, 2^i} (\nabla V(\tau ))|).
\label{7.54b}
 \end{equation} 
Hence,  recalling that $ P^0_{ x_0, 2^k}(V(\tau ))=0 $, using triangle inequality, we estimate for all $ j \ge k$
\begin{align}
| P^0_{ x_0, 2^j}(V(\tau ))|& \le | P^0_{ x_0, 2^k}(V(\tau ))| + c  \sum_{i=k}^{j}2^i \Big(| \omega (\tau )|_{ BMO}+ c\sup_{ i \ge k} | P^{ 0}_{ x_0, 2^i} (\nabla V(\tau ))|)\Big)
\cr
&\le   c 2^j\Big(| \omega(\tau ) |_{ BMO}+ c\sup_{ i \ge k} | P^{ 0}_{ x_0, 2^i} (\nabla V(\tau ))|)\Big).
\label{7.54a}
\end{align}
By virtue of Sobolev-Poincar\'e's inequality together with  \eqref{7.54a} and  \eqref{7.51a} we estimate for all $ \tau \in (0,T)$ 
\begin{align*}
& 2^{- j}\|   V(\tau )\|_{ L^\infty(B(x_0, 2^{ j+1}))} 
\\
&\le c 2^{- \frac{j}{2}}
\|   \nabla V(\tau )- \nabla P^1_{ x_0, 2^{ j+1}}(V(\tau ))\|_{ L^{ 2n}(B(x_0, 2^{ j+1}))} 
+c 2^{ -j}|P^1_{ x_0, 2^{ j+1}}(V(\tau ))| 
 \\
 &\le c |\nabla v(\tau )|_{ BMO} +    c 2^{ -j} |P^0_{ x_0, 2^{ j+1}}(V(\tau ))| + c |P^0_{ x_0, 2^{ j+1}}
 (\nabla V(\tau ))|
 \\
 &\le c |\omega (\tau )|_{ BMO} +     c |P^0_{ x_0, 2^{ j+1}}
 (\nabla V(\tau ))|.
\end{align*}
We also need to estimate the pressure.  Noting that 
\[
\osc_{ \frac{2n}{n+2}, 1}(\nabla \pi(\tau ) ; x_0, 2^j)  = 
\osc_{ \frac{2n}{n+2}, 1}(\nabla \Pi(V(\tau ), V(\tau ) ) ; x_0, 2^j)
\]
 consulting  \eqref{6.53} with $ r= \frac{2n}{n+2}$ and 
$ p=2$, and applying  \eqref{7.51a},  we get   
for all $\tau \in [0,T_{ \ast}),  x_0\in  \R^{n}$ and $ j\in \Z$, 
\begin{equation}
\osc_{ \frac{2n}{n+2}, 1}(\nabla \Pi(\tau ) ; x_0, 2^j) \le c  \Big(|\omega(\tau ) |_{ BMO} +
\sup_{ i \ge j} | P^{ 0}_{ x_0, 2^i} (\nabla V(\tau ))|\Big) S_{ 2, 1}(\osc_{2, 1}(V(\tau ); x_0))_j.
\label{7.52}
\end{equation} 
Inserting the above estimates into the right-hand side of  \eqref{7.53e} along with  \eqref{7.51a} and \eqref{7.52},  we get for  all $ j\in \Z$,  
\begin{align}
&\osc_{2, 1}(v(t); x_0, 2^j) 
 \le c\osc_{2, 1}(V(0); x_0, 2^j) 
 \cr
 & +  c\intl_{0}^{t} \Big(| \omega (\tau )|_{ BMO}
  +\sup_{ i \ge j} 
| P^{ 0}_{ x_0, 2^i} (\nabla V(\tau ))| \Big)S_{ 2, 1}(\osc_{2, 1}(V(\tau ); x_0))_j   d\tau. 
\label{7.54}
\end{align}
Applying $ S_{ 1,1}$ to both sides of the above inequality, and using Lemma\,2.1, we get, 
\begin{align}
& S_{ 1,1}(\osc_{2, 1}(v(t); x_0))_k
\cr
&  \le cS_{ 1,1}(\osc_{2, 1}(V(0); x_0))_k +  c \intl_{0}^{t} (| \omega (\tau )|_{ BMO} +\sup_{ i \ge k} 
| P^{ 0}_{ x_0, 2^i} (\nabla V(\tau ))|) S_{1, 1}(\osc_{2, 1}(V(\tau ); x_0))_k   d\tau.
\label{7.55}
\end{align}
On the other hand, estimating for all $ i \in \Z$
\[
 | P^{ 0}_{ x_0, 2^i} (\nabla V(\tau ))- P^{ 0}_{ x_0, 2^{ i+1}} (\nabla V(\tau ))| \le c 2^{ -i} \osc_{2, 1}(V(\tau ); x_0, 2^{ i+1}),
\]
using triangle inequality, we find for $ l\in \Z$
\begin{equation}
| P^{ 0}_{ x_0, 2^l} (\nabla V(\tau )) | \le  c \sum_{i=l}^{\infty} 2^{ -i}\osc_{2, 1}(V(\tau ); x_0, 2^i) + c | P^0_\infty(\nabla v(\tau ))|.  
\label{7.55b}
 \end{equation} 
In particular, 
\begin{equation}
\sup_{ i \ge k}  | P^{ 0}_{ x_0, 2^i} (\nabla V(\tau )) | \le  c \sum_{i=k}^{\infty} 2^{ -i}\osc_{2, 1}(V(\tau ); x_0, 2^i) + c | P^0_\infty(\nabla v(\tau ))|.  
\label{7.55a}
\end{equation}
Combining this estimate with \eqref{7.55},  multiplying the resultant inequality by $ 2^k$,    and taking the supremum over all $ x_0\in \R^{n}$, we arrive at  
\begin{align}
&\beta  _k (t)
\le c_0 \beta _k (0) +  c_0\intl_{0}^{t} \alpha (\tau )
\beta _k (\tau ) + \beta _k (\tau )^2  d\tau,
\label{7.56}
\end{align}
where 
\begin{align*}
\alpha(\tau ) &= | \omega (\tau )|_{ BMO} +  | P^0_\infty(\nabla v(\tau ))|, 
\\
\beta  _k(\tau ) &=  \sup_{ x_0\in \R^{n}}\sum_{i=k}^{\infty} 2^{ -i}\osc_{2, 1}(v(\tau )); x_0, 2^i),\quad  \tau \in [0,T]. 
\end{align*}
According to our assumption \eqref{1.37} we have $ \alpha \in L^1(0, T_{ \ast})$.  
We define 
\[
\var = \frac{1}{ 2c_0 c_1 e T_{ \ast}},\quad  \text{ where}\quad  c_1:= c_0 e^{  \intl_{0}^{T_{ \ast}}  \alpha (\tau ) d\tau }.
\]
Observing \eqref{1.34a}, we may choose $ k\in \Z$ such that $ \beta _k(0) \le \var $.  
Applying Gronwall's lemma,  we deduce    from \eqref{7.56} for all $ t\in (0,T_{ \ast})$
\begin{equation}
\beta_k (t) \le c_1  \var  e^{ c_0 t\supl_{ \tau \le t} \beta _k(\tau )  }.
\label{7.57}
\end{equation}
Without loss  of generality we may assume that $ c_1 \ge 1$.  Clearly, $ \beta (0) < \frac{1}{c_0 T_{ \ast}}$. 
Assume there exists $ t\in [0, T_{ \ast}]$ such that $ \beta _k(t) = \frac{1}{c_0 T_{ \ast}} $ and $ \supl_{ \tau \le t} 
\beta _k(\tau )
\le \frac{1}{c_0 T_{ \ast}}$. Then \eqref{7.57} would imply that 
\[
\beta_k (t)  = \frac{1}{c_0 T_{ \ast}} \le c_1  \var  e = \frac{1}{2 c_0 T_{ \ast}}  ,
\]
 which is a contradiction. Consequently, 
 \begin{equation}
\beta _k (t) \le \frac{1}{c_0 T_{ \ast}}\quad \forall\,t\in [0, T_{ \ast}].   
 \label{7.58}
 \end{equation}

3.  We verify that $ v_0 \in \mathscr{L}^{ 1+ \frac{\delta }{4}}_{1 (2, 1)}(\R^{n} )$. In fact, since $ \nabla v_0\in L^\infty(\R^{n} )$ by the help of H\"older's inequality and Poincar\'e's inequality we easily get 
\[
\osc_{ 2,1}(v_0; x_0, 2^j) \le c 2^{ \frac{j}{2}}\osc_{p,1}(v_0; x_0, 2^j)^{ \frac{1}{2}} \|\nabla v_0\|_{ \infty}^{ \frac{1}{2}}
\]   
By means of H\"older's inequality we find 
\begin{align*}
 \sum_{j\in \Z} 2^{ - j (1+ \frac{\delta }{4} ) }\osc_{ 2,1}(v_0; x_0, 2^j) 
 &\le \sum_{j=-\infty}^0 2^{ - j (1+ \frac{\delta }{4} ) }\osc_{ 2,1}(v_0; x_0, 2^j)  + c\|\nabla v_0\|_{ \infty}
 \\
 &\le 2 \sum_{j=-\infty}^0 2^{ j\frac{\delta }{4}}2^{ -j \frac{1+\delta }{2}}\osc_{p,1}(v_0; x_0, 2^j)^{ \frac{1}{2}} \|\nabla v_0\|_{ \infty}^{ \frac{1}{2}}+ c\|\nabla v_0\|_{ \infty}
 \\
 &\le c |v_0|_{ \mathscr{L}^{ 1+ \delta }_{q (p, 1)}}^{ \frac{1}{2}}\|\nabla v_0\|_{ \infty}^{ \frac{1}{2}}+ c\|\nabla v_0\|_{ \infty}.   
\end{align*}

4. Let $ k\in \Z$  chosen such that $ \beta _k(0) \le \var  $.   Let $ x_0\in \R^{n} $.  Let $ j\in \Z$.  
By $ (V, P)$ we denote a centered solution in $ x_0$  to  \eqref{euler}, which is equivalent to $ (v, p)$.  Since 
$ V(x_0, \tau )=0$ for all $ \tau \in (0,T)$, this yields 
\[
\lim_{i \to -\infty} P^0_{ x_0, 2^i} (V(\tau )) =0. 
\]
Using  \eqref{7.54b} and triangle inequality, we get 
\begin{align}
  2^{ -j}| P^{ 0}_{ x_0, 2^j} (V(\tau ))| &\le  c 2^{ -j}\sum_{i=-\infty}^{j} 2^i \Big(|\omega |_{ BMO} +
  \sup_{ i \in \Z} | P^0_{ x_0, 2 ^i}(\nabla V(\tau ))| \Big) 
  \cr
 &\le c\Big(|\omega |_{ BMO} +
  \sup_{ i \in \Z} | P^0_{ x_0, 2 ^i}(\nabla V(\tau ))|  \Big).
\label{7.54c}
\end{align}
Inserting \eqref{7.54c} into the right-hand side of \eqref{7.54}, we obtain 
\begin{align}
&\osc_{2, 1}(v(t); x_0, 2^j) 
 \le c\osc_{2, 1}(V(0); x_0, 2^j) 
\cr
&+  \intl_{0}^{t} \Big(|\omega (\tau )|_{ BMO} +
\sup_{ i \in \Z} 
P^{ 0}_{ x_0, 2^i} (\nabla V(\tau ))\Big) S_{ 2, 1}(\osc_{2, 1}(V(\tau ); x_0))_j   d\tau. 
\label{7.60}
\end{align}  

\vspace{0.5cm}  
We proceed with the estimation of $ \sup_{ i \in \Z} | P^0_{ x_0, 2 ^i}(\nabla V(\tau ))|$. Clearly, by 
 \eqref{7.55b} we see that for all $ i\in \Z$
\begin{align*}
 | P^0_{ x_0, 2 ^i}(\nabla V(\tau ))|
&\le  c \sum_{m\in \Z} 2^{ -m}\osc_{2, 1}(V(\tau ); x_0, 2^m) + c | P^0_\infty(\nabla v(\tau ))|
\\
&\le  c \sum_{m=-\infty}^{ k-1} 2^{ -m}\osc_{2, 1}(V(\tau ); x_0, 2^m) + c \beta _k(\tau )+ c | P^0_\infty(\nabla v(\tau ))|.  
\end{align*}  

 By the aid of \eqref{7.38} with $ p=2$ and $ \frac{\delta }{4}$ in place of $ \delta $ (cf. Lemma\,\ref{lem7.3}) 
 together with \eqref{7.51a} we find 
\[
 \sup_{ i \in \Z} | P^0_{ x_0, 2 ^i}(\nabla V(\tau ))| \le c \Big(2^{k \frac{\delta }{4}} + \beta _k(\tau )  +| \omega (\tau )|_{ BMO}\log (1+ | v(\tau )|_{ \mathscr{L}^{ 1+ \frac{\delta }{4} }_{ 1 (2, 1)}}) + |P^0_\infty(\nabla V(\tau )) |\Big).
\]

  Inserting this estimate into the right-hand side of \eqref{7.60}, and applying $ S_{1+ \frac{\delta }{4} ,1 }$ to both sides, using Lemma\,\ref{lem10.1},  we are led to 
  \begin{align}
& S_{ 1+\frac{\delta }{4} ,1 }(\osc_{2, 1}(v(t); x_0))_j 
 \le cS_{ 1+\frac{\delta }{4} ,1}(\osc_{2, 1}(V(0); x_0))_j 
\cr
&+  c\intl_{0}^{t} \Big(\alpha (\tau )   \log (1+ | v(\tau )|_{ \mathscr{L}^{ 1+\frac{\delta }{4}}_{ 1 (2, 1)}})
+ 2^{ \delta k} + \beta _k(\tau )\Big) S_{ 1+ \frac{\delta }{4}, 1}(\osc_{2, 1}(V(\tau ); x_0))_j   d\tau. 
\label{7.61}
\end{align}  
Multiplying both sides of \eqref{7.61} by $ 2^{ - j (1+ \frac{\delta }{4} )}$, taking the supremum over all $ x_0\in \R^{n}$ after summing over $j\in \Bbb Z$, and observing \eqref{7.58},  we  
deduce that 
  \begin{align}
& | v(t)|_{ \mathscr{L}^{ 1+\frac{\delta }{4} }_{ 1 (2, 1)}}
\cr
& \le c| v_0|_{ \mathscr{L}^{ 1+\frac{\delta }{4} }_{ 1 (2, 1)}}
+ c \intl_{0}^{t} (1+\alpha (\tau ))   \log (e+ | v(\tau )|_{ \mathscr{L}^{ 1+\frac{\delta }{4} }_{ 1 (2, 1)}})| v(\tau )|_{ \mathscr{L}^{ 1+\frac{\delta }{4} }_{ 1 (2, 1)}}   d\tau,
\label{7.62}
\end{align}   
 for a constant $ c>0$ independent of $ t$. Applying Gronwall's lemma,  we obtain from \eqref{7.62} that 
 \[
| v(t)|_{ \mathscr{L}^{ 1+\frac{\delta }{4} }_{ 1 (2, 1)}} \le \exp \Bigg[c | v_0|_{ \mathscr{L}^{ 1+\frac{\delta }{4} }_{ 1 (2, 1)}} 
\exp \bigg( \intl_{0}^{T_{ \ast}} (1+ \alpha (\tau )) d\tau \bigg)\Bigg]. 
\]  
 
  Accordingly, $ v\in L^\infty(0, T_{ \ast}; \mathscr{L}^{ 1+\frac{\delta }{4}}_{ 1 (2, 1)}(\R^{n}))$. Taking into account  \eqref{7.58}, we see that  $  v\in L^\infty(0, T_{ \ast}; \mathscr{L}^{ 1 }_{ 1 (2, 1)}(\R^{n}))$. 
In particular, $ \nabla v$ is bounded. Repeating the above argument and recalling $ v_0\in 
\mathscr{L}^{ 1+\delta  }_{ q (p, 1)}(\R^{n})$, we obtain $ v\in L^\infty(0, T_{ \ast}; \mathscr{L}^{ 1+\delta }_{ q (p, 1)}(\R^{n}))$, which completes the proof of the theorem.  \hfill \Beweisende 
  
hspace{0.5cm}
$$\mbox{\bf Acknowledgements}$$
Chae was partially supported by NRF grants 2016R1A2B3011647, while Wolf has been supported 
supported by NRF grants 2017R1E1A1A01074536.
The authors declare that they have no conflict of interest.

 \end{document}